\documentclass[12pt]{article}
\usepackage{amsmath,amsfonts,amssymb,mathabx,txfonts}
\usepackage{mathrsfs,upgreek}
\usepackage{dsfont}
\usepackage[utf8]{inputenc}
\usepackage[T1]{fontenc}
\usepackage{lmodern} \usepackage{ucs}
\usepackage{fullpage}

\usepackage{graphicx}
\usepackage{caption}
\usepackage{subcaption}

\usepackage{array}
\usepackage{multicol}
\usepackage{multirow}
\usepackage{color}
\usepackage{authblk}
\usepackage{esint}
\usepackage{enumitem,hyperref}

\usepackage{accents}

\usepackage{tikz}
\usetikzlibrary{decorations}
\usetikzlibrary{decorations.pathreplacing}
\usetikzlibrary{arrows}

\newtheorem{lemma}{Lemma}[section]
\newtheorem{theo}[lemma]{Theorem}
\newtheorem{rmk}[lemma]{Remark}
\newtheorem{proposition}[lemma]{Proposition}

\makeatletter
\def\namedlabel#1#2{\begingroup
    #2%
    \def\@currentlabel{#2}%
    \phantomsection\label{#1}\endgroup
}
\makeatother

\makeatletter
\renewcommand*{\eqref}[1]{%
  \hyperref[{#1}]{\textup{\tagform@{\ref*{#1}}}}%
}
\makeatother

\newcommand{\QED}{\mbox{}\hfill \raisebox{-0.2pt}{\rule{5.6pt}{6pt}\rule{0pt}{0pt}} \medskip\par}

\newcommand{\R}{\mathbb{R}}

\newcommand{\U}{\mathcal{U}}

\newcommand{\ds}{\displaystyle}
\newcommand{\ud}{\, {\mathrm{d}}}

\title{A simple testbed  for stability analysis of  quantum dissipative  systems}

\author{Thierry~Goudon\thanks{ {\tt thierry.goudon@inria.fr}}}
\author{Simona Rota Nodari\thanks{ {\tt simona.rotanodari@univ-cotedazur.fr}}}
\date{}

\affil[1]{\small Universit\'e C\^ote d'Azur, Inria,  CNRS, LJAD,

Parc Valrose, F-06108 Nice, France}

\begin{document}
\maketitle

\begin{abstract}
We study a two-state quantum system 
with a non linearity intended to describe interactions with a complex environment, arising
 through a non local coupling term.
We study the stability of particular solutions, obtained as  constrained extrema of the energy functional of the system.
The simplicity of the model allows us to justify a complete stability analysis.
This is the opportunity to review in details the techniques to investigate the stability issue.
We also bring out the limitations of perturbative approaches based on simpler asymptotic models.  
\end{abstract}

\vspace*{.5cm}
{\small
 Particles interacting with a vibrational field.   Orbital stability.
\\[.4cm]

\noindent{\bf Math.~Subject Classification.} 
35Q40 
35Q51 

\section{Introduction}

In this work, we consider a simple quantum system
characterized by a single degree of freedom which can take only two values, hereafter referred to as 0 and 1.
The quantum system interacts with its environment, 
the description of which is embodied into a vibrational field, oscillating in some abstract direction $z\in \mathbb R^n$.
 Therefore the evolution of the system is governed by the ODE system 
\begin{equation}\label{S1}
\begin{array}{l}
i \ds\frac
{\ud}{\ud t} u_0(t)=u_0(t) -u_1(t)+u_0(t)\ds\int_{\mathbb R^n} \sigma(z)\psi_0(t,z)\ud z,
\\
i \ds\frac
{\ud}{\ud t} u_1(t)=u_1(t) -u_0(t)+u_1(t)\ds\int_{\mathbb R^n} \sigma(z)\psi_1(t,z)\ud z,
\end{array}
\end{equation}
coupled to the wave equations
\begin{equation}\label{W1}
\begin{array}{l}
\left(\ds\frac{1}{c^2}\partial_{tt}^2-\Delta \right)\psi_0(t,z)=-\sigma(z) |u_0(t)|^2,
\\
\left(\ds\frac{1}{c^2}\partial_{tt}^2-\Delta \right)\psi_1(t,z)=-\sigma(z) |u_1(t)|^2.
\end{array}
\end{equation}
These equations are completed by initial data
\begin{equation}\label{id}
(u_0,u_1,\psi_0,\partial_t \psi_0,\psi_1,\partial_t \psi_1)\big|_{t=0}=
(u_{0,\mathrm{init}},u_{1,\mathrm{init}},\psi_{0,\mathrm{init}},\varpi_{0,\mathrm{init}},\psi_{1,\mathrm{init}},\varpi_{1,\mathrm{init}}).\end{equation}
The coupling is embodied into the form function $z\in\mathbb R^n\mapsto \sigma(z)$, which, throughout the paper is assumed to be non-negative, smooth, with fast enough decay (say compactly supported to fix ideas).
The free problem ($\sigma=0$) reduces 
to
\begin{equation}\label{0cou}
\ds\frac
{\ud}{\ud t}
\begin{pmatrix} u_0\\ u_1\end{pmatrix} = \ds\frac1i
\begin{pmatrix} 1& -1 \\
  -1& 
 1 \end{pmatrix}\begin{pmatrix} u_0\\ u_1\end{pmatrix} .
\end{equation}
We infer that the system oscillates with frequency 2 around  a constant state: the solutions of \eqref{0cou} read
$$
\begin{pmatrix}
u_0(t)\\ u_1(t)
\end{pmatrix}=\frac12
\begin{pmatrix} 1 & 1
\\
1 & -1\end{pmatrix}
\begin{pmatrix}
(u_{0,\mathrm{init}}+u_{1,\mathrm{init}})  \\
(u_{0,\mathrm{init}}-u_{1,\mathrm{init}}) e^{-2it} 
\end{pmatrix}.
$$
Hence, we are wondering how the coupling ($\sigma\neq 0$) impacts this simple dynamics.
It is also worth considering the large  speed regime $c\to \infty$ which leads to the following non-linear
ODE system
\begin{equation}\label{Hartree}
\ds\frac
{\ud}{\ud t}\begin{pmatrix} u_0\\ u_1\end{pmatrix} = \ds\frac1i
\begin{pmatrix} 1-\kappa |u_0|^2& -1\\
 -1& 
 1-\kappa|u_1|^2 \end{pmatrix}\begin{pmatrix} u_0\\ u_1\end{pmatrix} 
\end{equation}
where 
\begin{equation}\label{defkappa}
\kappa =
\ds\int_{\mathbb R^n}
 \sigma(z)(-\Delta)^{-1} \sigma(z)\ud z
= \ds\int_{\mathbb R^n}
 \ds\frac{|\widehat \sigma(\xi)|^2}{|\xi|^2}\ds\frac{\ud \xi}{(2\pi)^n}>0.
\end{equation}
It will be interesting to compare the behavior 
of the asymptotic model \eqref{Hartree} with \eqref{S1}-\eqref{W1};
in particular we are going to point out the limitations of a perturbative approach that would try to
deduce properties of \eqref{S1}-\eqref{W1} from the analysis of \eqref{Hartree}.
\\

According to the ideas of quantum mechanics,  the $|u_j|$'s represent the probability of being in the state  labelled by $j$; in turn, the
  total probability  should be one: we always have
  \begin{equation}\label{consL2}
  |u_0(t)|^2 + |u_1(t)|^2=1.
 \end{equation}
 If this property holds initially, we check that it holds forever.
Moreover, the equations describe the energy exchanges between the quantum system and the environment
which traduces into an additional conservation property, namely, we have (detailed computations can be found in Appendix~\ref{App_Cons}, 
 but the result also directly follows from the symplectic form of the problem, exhibited below,  
 combined to  Noether's theorem)
\begin{equation}\label{EnerSW}
\begin{array}{l}
\textrm{for \eqref{S1}-\eqref{W1}:}\qquad \ds\frac
{\ud}{\ud t}
\left(
\ds\frac{|u_0-u_1|^2}{2}
+\ds\frac14
\ds\int_{\mathbb R^n}
\Big(
\ds\frac{1}{c^2}
(|\partial_t \psi_0|^2+|\partial_t \psi_1|^2) +|\nabla\psi_0|^2+|\nabla\psi_1|^2\Big)\ud z
\right.
\\
\hspace*{9cm}
\left.
+\ds\frac12\ds\int_{\mathbb R^n} \sigma(\psi_0|u_0|^2+\psi_1|u_1|^2)\ud z
\right)=0
\end{array}\end{equation}
 which becomes
\begin{equation}\label{EnerHa}
\textrm{for \eqref{Hartree}:}\qquad
\ds\frac
{\ud}{\ud t}
\left(
\ds\frac{|u_0-u_1|^2}{2}
-\ds\frac\kappa 4(|u_0|^4+|u_1|^4)\right)=0
\end{equation}
for the asymptotic model \eqref{Hartree}.
These conservation properties play a central role in the analysis of the equations.
\\

The question we address
  comes from the 
  modeling of \emph{quantum open systems}.
The motivation, inspired from the seminal work of Caldeira and Leggett \cite{CL}, is to understand how the interactions with the environment induce some kind of dissipative   
effects. The intuition is that the quantum system exchanges energy with the vibrational field, 
and the energy is eventually evacuated ``at infinity'' in the $z$-direction; this mechanism can be interpreted as a sort of friction acting 
on the quantum system.  
For the sake of concreteness, 
the energy transfer mechanisms at work between the quantum system and the environment with the model  \eqref{S1}-\eqref{W1}
are illustrated in 
Figure~\ref{figenergytransf}
which show typical evolutions of the different  contributions, wave and particle, to the total energy:   
albeit these curves are suggestive, in fact,  
  they correspond to very different behaviors of the system, as we shall discuss below.
Such an issue has been studied in details for the case of a single classical particle in \cite{BdB}, 
 where the dissipation mechanisms are explicitly exhibited.
 This situation 
 has been further investigated in 
 \cite{AdBLP,dBP,dBPS,dBS, bDPL};
we also refer the reader to \cite{JP2} 
 or \cite{KKSb}
for different, but related,  viewpoints on the dynamic of a classical particle coupled to a complex environment.
Dealing with many classical particles
leads to consider Vlasov-like equations \cite{dBGV,GV,Vi2}, and the dissipation effects can be interpreted as a sort of Landau damping.
The friction mechanisms are intimately related to the  dispersion properties of the wave equation that need to be strong enough, 
an effect driven by the condition $n\geq 3$ on the $z$-direction, that will be 
assumed throughout the paper.
In particular, it can be noticed that it guarantees 
the quantity defined by \eqref{defkappa} to be finite.
We refer the reader to \cite{Vi2} for detailed comments about this assumption.
Coming back to quantum particles, one is led to systems coupling the Schr\"odinger equation with a wave equation: the model
\begin{equation}
\label{SchW}
\begin{array}{l}
i\partial _tu+\ds\frac{\Delta_x}2 u=\Phi u,
\\
\Phi(t,x)=\ds\int \sigma_1(x-y)\sigma(z)\psi(t,y,z)\ud z\ud y,
\\
\Big(\ds\frac1{c^2}\partial_{tt}^2\psi-\Delta_z \psi\Big)(t,x,z)=-\sigma(z)\ds\int\sigma_1(x-y)|u|^2(t,y)\ud y,\end{array}\end{equation}
is the quantum analog of the equation introduced in \cite{BdB} (other quantum frameworks are discussed for instance in \cite{Fau,DS,JP1}).
The equation is analysed when the variable $x$ lies in $\mathbb R^d$ in \cite{Vi3} and ground states can be identified by variational approaches.
However, the stability analysis of the ground states is delicate because of 
the non local definition of the self-consistent potential, and the arguments developped 
for NLS ($\Phi=-|u|^2$ in the first equation of \eqref{SchW}) or Schr\"odinger-Newton ($\Phi=\frac{1}{|\cdot|}\star|u|^2$ with $d=3$) do not adapt directly (note at least that here the coupling has a more dynamical nature).
The attempt in this direction presented in \cite{Vi3}, completed by the numerical investigation in  \cite{Vi1}, relies on  
a perturbative approach, inspired from \cite{Lenzmann}.
However, it induces some restrictions which are not completely satisfying.
In order to understand this difficulty, we have studied the simpler framework of plane waves  ($x$ lies in the torus $\mathbb T^d$) in 
\cite{SRN3}, where the Hamiltonian structure is further exploited, in the spirit of the pioneering work \cite{GSS}, see also the recent overview \cite{SRN1}. 
It allows us to identify fundamental differences between \eqref{SchW} and its asymptotic counterpart as $c\to \infty$; in particular, the coupling with the wave equation induces spectral difficulties
which make perturbative arguments inoperative.
We wish to explore in further details these issues by considering the simpler systems \eqref{S1}-\eqref{W1}and \eqref{Hartree}.

\begin{figure}[!htbp]
  \begin{subfigure}{0.32\textwidth}
 \includegraphics[height=4cm]{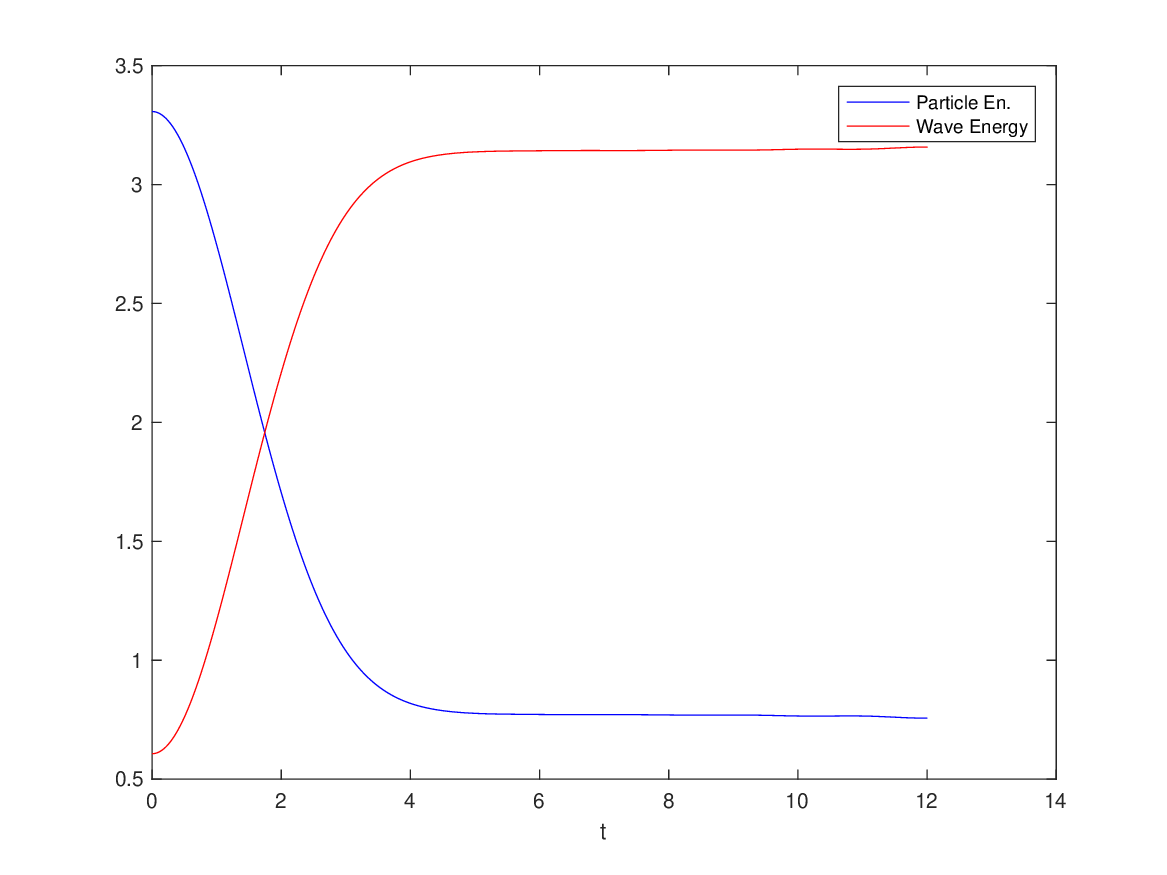}
 \caption{}
   \end{subfigure}
 \begin{subfigure}{0.32\textwidth}
  \includegraphics[height=4cm]{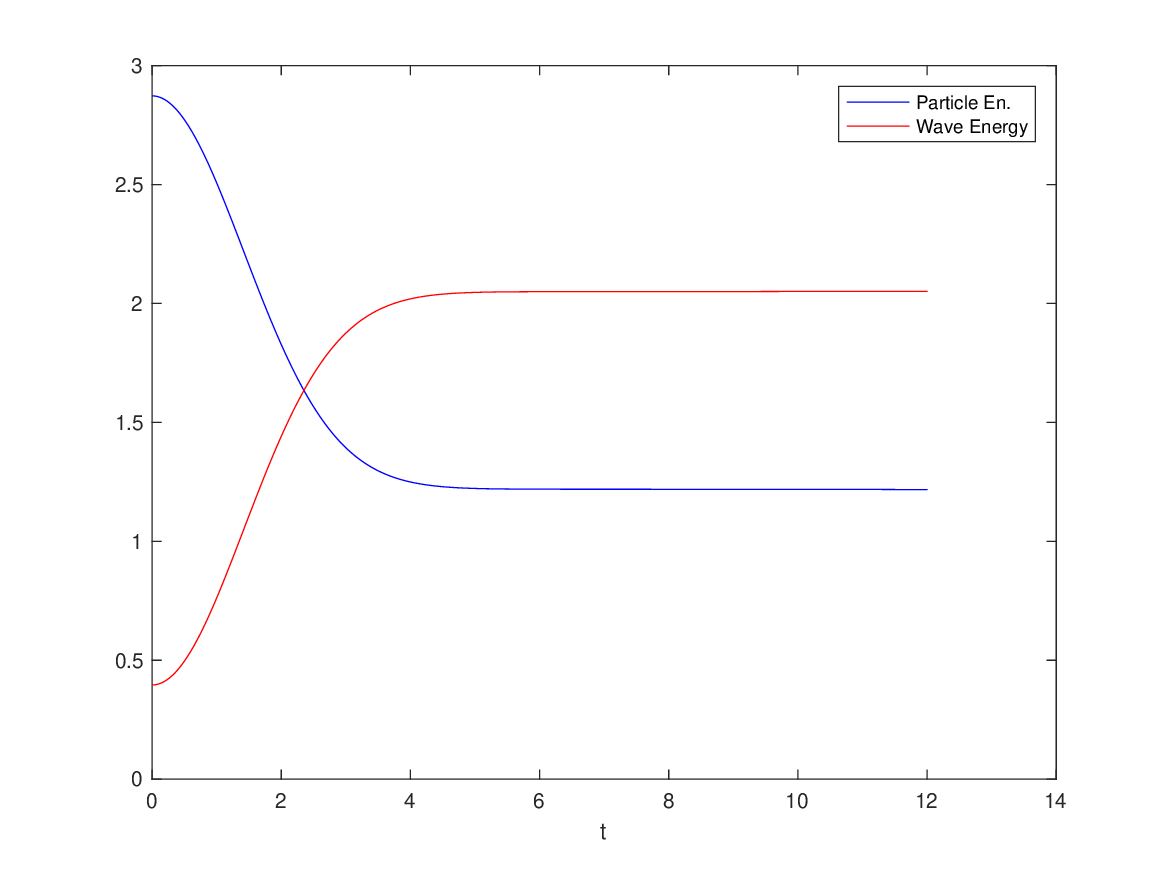}
  \caption{}
    \end{subfigure}
  \begin{subfigure}{0.32\textwidth}
    \includegraphics[height=4cm]{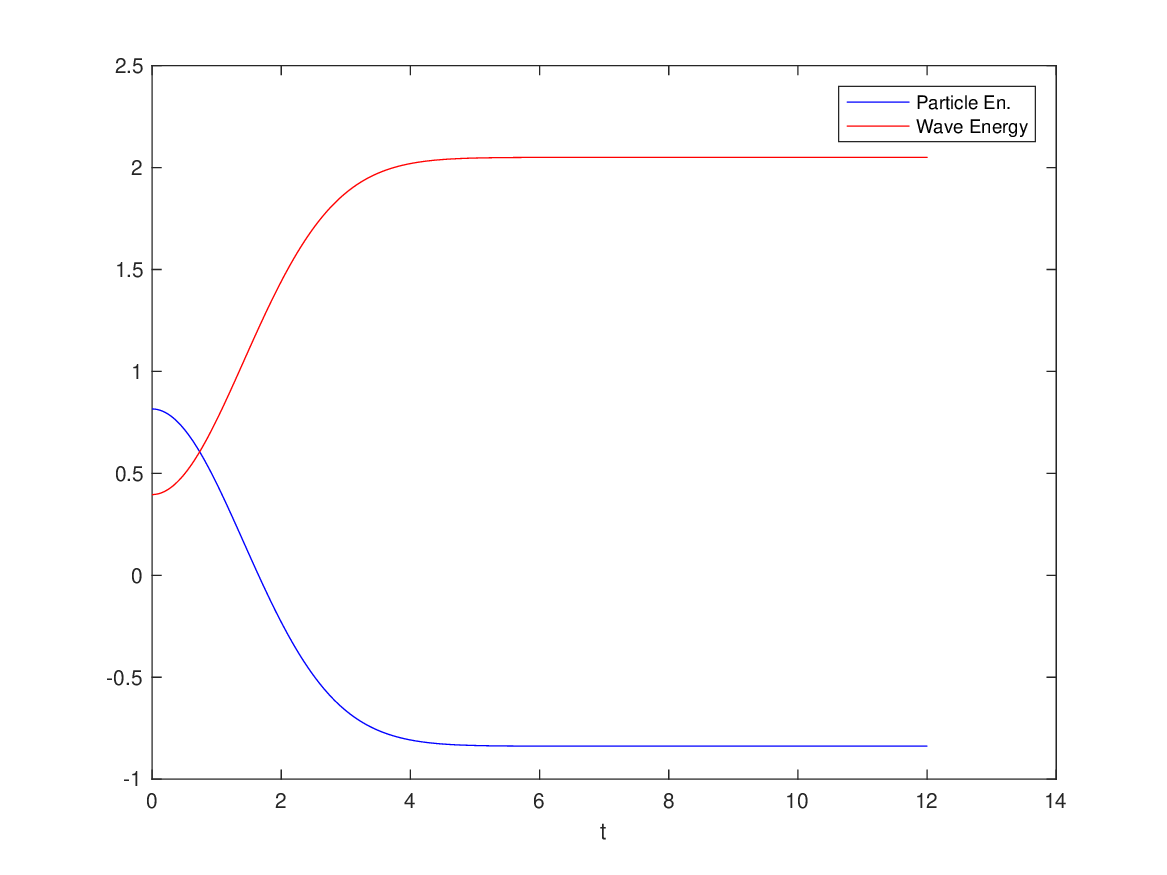}
    \caption{}
    \end{subfigure}
  
 \caption{  Evolution of the ``Wave contribution'' 
 $\frac14\sum_{j=0}^1
\int_{\mathbb R^n}
\big(
\frac{|\partial_t \psi_j|^2}{c^2}+|\nabla\psi_j|^2\big)\ud z$
and
the ``Particle contribution''
 $\frac{|u_0-u_1|^2}{2}
+\frac12\sum_{j=0}^1\int_{\mathbb R^n} \sigma \psi_j|u_j|^2\ud z
$
 to the Total Energy \eqref{EnerSW} associated to \eqref{S1}-\eqref{W1}.
 The simulations correspond to various cases that will be discussed in details below:
 (a) $\tau=+1$ and large $\kappa$,
 (b) $\tau=-1$, (c): $\tau=+1$ and small $\kappa$ 
 }\label{figenergytransf} 
 \end{figure}

In what follows, we pay attention to solutions  of \eqref{S1}-\eqref{W1} or \eqref{Hartree}, where 
the quantum particles distribution has the specific form $e^{i\omega t}(U_{*0},U_{*1})$, with $U_{*0},U_{*1}$, fixed complex numbers. These solutions can be classified in terms of extrema of the energy.
  The question we address is about the stability of these specific solutions.
  At first sight, 
  the problem under consideration can be seen as a discrete version of the non linear Schr\"odinger equation:
  we roughly interpret
  $u_0-u_1$ and $u_1-u_0$ as the discrete laplacian $(\Delta^du)_0=\frac{-u_{-1}+2u_0-u_1}{2}$,
  $(\Delta^du)_1=\frac{-u_{0}+2u_1-u_2}{2}$
   endowed with periodic conditions $u_{-1}=u_1$, $u_0=u_2$~!
  Stability analysis relies on the properties of
  the energy functional which can be used as a Lyapounov functional, and establishing coercivity properties
  is key for proving the orbital stability of the ground state, see
  \cite{Weinstein1, Weinstein2} and the recent review \cite{Tao}.
  A quite general framework has been set up in 
  \cite{GSS,GSS2}, see also  \cite{SRN1,SRN2}, intended to cover the 
  analysis of a wide class of Hamiltonian systems.
  However, the coupling with a vibrational environment lead to difficulties of a different nature, which are not covered by the abstract framework of \cite{GSS,GSS2} since 
 the nonlinearity has the form $\Phi u$, where the potential $\Phi$ 
  is non local both ``in space'' (here it means that it mixes the two states 0 and 1) and in time, with some kind of memory effects, so that the arguments of \cite{GSS,GSS2} do not apply.
  \\

 The interest of the two-level model is to be both  simple enough to allow us to perform many explicit computations, and rich enough to exhibit interesting phenomena; in turn
 \begin{itemize}
 \item we are able to provide a complete stability analysis for the models \eqref{S1}-\eqref{W1} and \eqref{Hartree};
 \item we review in full details the techniques for investigating such systems, and explain how they can be adapted to handle the non local coupling;
 \item it allows us to clarify where are the main difficulties and it provides valuable hints to 
 study more complex models. We expect this work to provide useful ideas to go back to the more challenging problem \eqref{SchW}.
\end{itemize}

The paper is organized as follows.
In Section~\ref{Ham}, we discuss the Hamiltonian structure of the problem 
  and make the connection appear between extrema of the energy
  functional
  and specific solutions with the  form $u(t)=e^{i\omega t}(U_{*0},U_{*1})$.
  Section~\ref{AnHar} is devoted to the analysis of the asymptotic system \eqref{Hartree}, which is a mere ODE system.
  In Section~\ref{AnSW}, we discuss the system \eqref{S1}-\eqref{W1}. 
Throughout the paper, numerical simulations 
illustrate the obtained statements.

\section{Hamiltonian formulation, extrema  of the energy and traveling-wave-like solutions}
\label{Ham}
 
 Throughout the paper, we  split a complex number  $u=q+ip$, where $q,p$ are real valued. 
Coming back to the unknown describing the quantum state, it makes the following correspondance appear
\begin{equation}\label{notation}
U=\begin{pmatrix} u_0\\u_1\end{pmatrix}\in \mathbb C^2
\longleftrightarrow X=\begin{pmatrix} q_0\\p_0\\ q_1\\p_1 \end{pmatrix}\in \mathbb R^4.\end{equation}

\subsection{Analysis of the asymptotic model}

We start with the simpler system \eqref{Hartree}.
Let us introduce the  function
\[
\mathscr H: (u_0,u_1)\in \mathbb C^2\longmapsto 
\ds\frac{|u_0-u_1|^2}{2}
-\ds\frac\kappa 4(|u_0|^4+|u_1|^4).\]
We have observed that $t\mapsto \mathscr H (u_0(t),u_1(t))$ is conserved 
by the differential system \eqref{Hartree}.
This property can be interpreted as a consequence of the following reformulation of the problem, in terms 
of the real valued quantities defined by \eqref{notation}.
The conserved quantity becomes
\begin{equation}\label{EnerHaBis}
\mathscr H (X)=\ds\frac{|q_0-q_1|^2}{2}+\ds\frac{|p_0-p_1|^2}{2}
-\ds\frac\kappa 4(|q_0|^2+|p_0|^2)^2-\ds\frac\kappa 4(|q_1|^2+|p_1|^2)^2,
\end{equation}
and \eqref{Hartree} can be cast in the \emph{symplectic form}
\begin{equation}\label{Hartreebis}
\partial_t X=\mathscr J\nabla_X\mathscr H(X),
\end{equation}
with $\mathscr J$ the  skew-symmetric matrix
\[\mathscr J=\begin{pmatrix}
  0 & 1 & 0 & 0
  \\
    -1 & 0 & 0& 0
    \\
     0 & 0& 0 & 1
     \\
     0 & 0& - 1 &0
     \end{pmatrix}.\]
We are  interested in specific solutions of \eqref{Hartree}, 
having the special form $(e^{i\omega t} U_{*0},e^{i\omega t} U_{*1})$ where 
 $U_{*0}=Q_{*0}+iP_{*0}$, and  $U_{*1}=Q_{*1}+iP_{*1}$ are fixed complex numbers.
 We are led to the relation 
 \begin{equation}\label{EulerLag}
 -\omega X_*=\nabla_X\mathscr H(X_*)
 =
 \begin{pmatrix}
 Q_{*0}-Q_{*1}-\kappa(Q_{*0}^2+P_{*0}^2)Q_{*0}
 \\
 P_{*0}-P_{*1}-\kappa(Q_{*0}^2+P_{*0}^2)P_{*0}
 \\
 Q_{*1}-Q_{*0}-\kappa(Q_{*1}^2+P_{*1}^2)Q_{*1}
 \\
 P_{*1}-P_{*0}-\kappa(Q_{*1}^2+P_{*1}^2)P_{*1}
 \end{pmatrix}
 ,\end{equation}
 which similarly arises when searching for the extrema of $\mathscr H$ under the constraint of fixed $L^2$ norm $|U_{*0}|^2+|U_{*1}|^2=1$,
  $\omega$ being interpreted as the associated Lagrange multiplier. We thus focus on this optimization viewpoint.

We write $u_j=r_je^{i\theta_j}=q_j+ip_j$, with $r_j\geq 0$ and $\theta_j\in [0,2\pi)$ and we realize that all terms in
$\mathscr H(u_0,u_1)$ do not depend on the angles $\theta_j$, 
but 
\[
|u_0-u_1|^2=r_0^2+r_1^2-2r_0r_1\cos(\theta_1-\theta_0)
 \]
 so that 
 \[
 (r_0-r_1)^2\leq |u_0-u_1|^2\leq (r_0+r_1)^2
 \]
 holds.
 The inequalities are saturated when $\theta_1=\theta_0\ \mathrm{ mod}(2\pi)$ (left) or  $\theta_1=\theta_0\ \mathrm{ mod}(\pi)$ (right).
 If $(u_0,u_1)$ minimizes $\mathscr H$ 
 over the unit sphere of $\mathbb C^2$, we deduce from
 $\mathscr H(r_0,r_1)\leq \mathscr H(u_0,u_1)
$ and $r_0^2+r_1^2=1$, that 
$(r_0,r_1)\in [0,1]\times
 [0,1]$ equally reaches the minimum.
 Conversely, if $(q_0,q_1)$ minimizes  $\mathscr H$ 
 over the unit sphere of $\mathbb R^2$, then, for any $u_j=r_je^{i\theta_j}$, with $r_0^2+r_1^2=1$, we get 
  $\mathscr H(q_0,q_1)\leq \mathscr H(r_0,r_1)\leq \mathscr H(u_0,u_1)$
  so that $(q_0,q_1)$ minimizes $\mathscr H$ 
 over the unit sphere of $\mathbb C^2$. A similar equivalence holds for maximizing $\mathscr H$.
 
 Therefore, all extrema can be described by restricting first to the case $p_0=p_1=0$, and then, from the obtained (real valued)
 optima $(q_0,q_1)$, by setting $u_0=e^{i\theta_0} q_0$, $u_1= e^{i\theta_0} q_1$, $\theta_0\in [0,2\pi)$.
 Moreover, we should also bear in mind the conservation of the $L^2$ norm, so that we are actually interested in  
 extrema over the sphere $\{(q_0,q_1)\in \mathbb R^2,\ |q_0|^2+|q_1|^2=1\}$.
 Accordingly, we can reinterpret the problem as to optimize  a function  
 of a mere scalar variable 
 \[
\theta\in [0,2\pi)\longmapsto  \mathscr H_1(\theta)=\ds\frac{(\cos(\theta)-\sin(\theta))^2}{2}-\ds\frac\kappa 4
 (\cos^4(\theta)+\sin^4(\theta)).\]
 For reader's convenience, graphs of $\theta\mapsto \mathscr H_1(\theta)$ are plotted for several values of $\kappa$ in Fig.~\ref{GrH1}.
 We have
 \[
 \mathscr H_1'(\theta)=-(\cos^2(\theta)-\sin^2(\theta))+\kappa\cos(\theta)\sin(\theta)(\cos^2(\theta)-\sin^2(\theta))
 =
 \ds\frac\kappa 2\cos(2\theta)\left(\sin(2\theta)-\ds\frac2\kappa\right).\]
 It vanishes when $\theta=\frac\pi4$ which yields the solution $q_0=1/\sqrt2$, $q_1=1/\sqrt2$,  or $\theta=\frac{3\pi}4$,
 which yields the solution 
  $q_0=1/\sqrt2$, $q_1=-1/\sqrt2$.
  If the smallness condition $0<\kappa<2$ holds, this completely describes the extrema of the function $\mathscr H_1$.
  When $\kappa>2$, 
  we can find other solutions by setting $\theta=\frac{\mathrm{arcsin}(2/\kappa)}{2}\in (0,\pi/4)$
  and $\theta=\frac{\pi}{2}-\frac{\mathrm{arcsin}(2/\kappa)}{2}\in (\pi/4/,\pi/2)$.
 We have 
 \[
 \mathscr H_1''(\theta)=
\kappa\left( \cos^2(2\theta)
-\sin(2\theta)\left(\sin(2\theta)-\ds\frac2\kappa\right)
\right).
 \]
Therefore, we distinguish the following cases:
\begin{itemize}
\item if $0<\kappa<2$, $\theta=\pi/4$ minimizes the energy 
 ($\mathscr H_1''(\pi/4)=\frac\kappa2(\frac2\kappa-1)>0$)
and  $\theta=3\pi/4$ maximizes the energy ($\mathscr H_1''(3\pi/4)=-\frac\kappa2(\frac2\kappa+1)<0$): 
we have $\mathscr H_1(\pi/4)=-\frac\kappa8\leq \mathscr H_1(\theta)\leq \mathscr H_1(3\pi/4)=1-\frac\kappa8$;
\item if $\kappa>2$, 
$\theta^+_\kappa=\frac{\mathrm{arcsin}(2/\kappa)}{2}$, 
$\theta^-_\kappa=\frac{\pi}{2}-\frac{\mathrm{arcsin}(2/\kappa)}{2}$ minimize the energy 
($\mathscr H_1''(\theta_\kappa^\pm)=\kappa\cos^2(2\theta\kappa^\pm)>0$),
$\theta=\pi/4$ is a local maximum of the energy ($\mathscr H_1''(\pi/4)=\frac\kappa2(\frac2\kappa-1)<0$)
and  $\theta=3\pi/4$ maximizes the energy ($\mathscr H_1''(3\pi/4)=-\frac\kappa2(\frac2\kappa+1)<0$); 
we have $\mathscr H_1(\theta_\kappa^\pm)=\frac12(1-\kappa/2-1/\kappa))\leq \mathscr H_1(\theta)\leq \mathscr H_1(3\pi/4)=1-\frac\kappa8$ and
$\mathscr H_1(\pi/4)=-\frac\kappa8\in (\mathscr H_1(\theta_\kappa^ \pm),\mathscr H_1(3\pi/4))$.
\end{itemize} 

\begin{figure}[!h]
\begin{center}
\includegraphics[height=8cm]{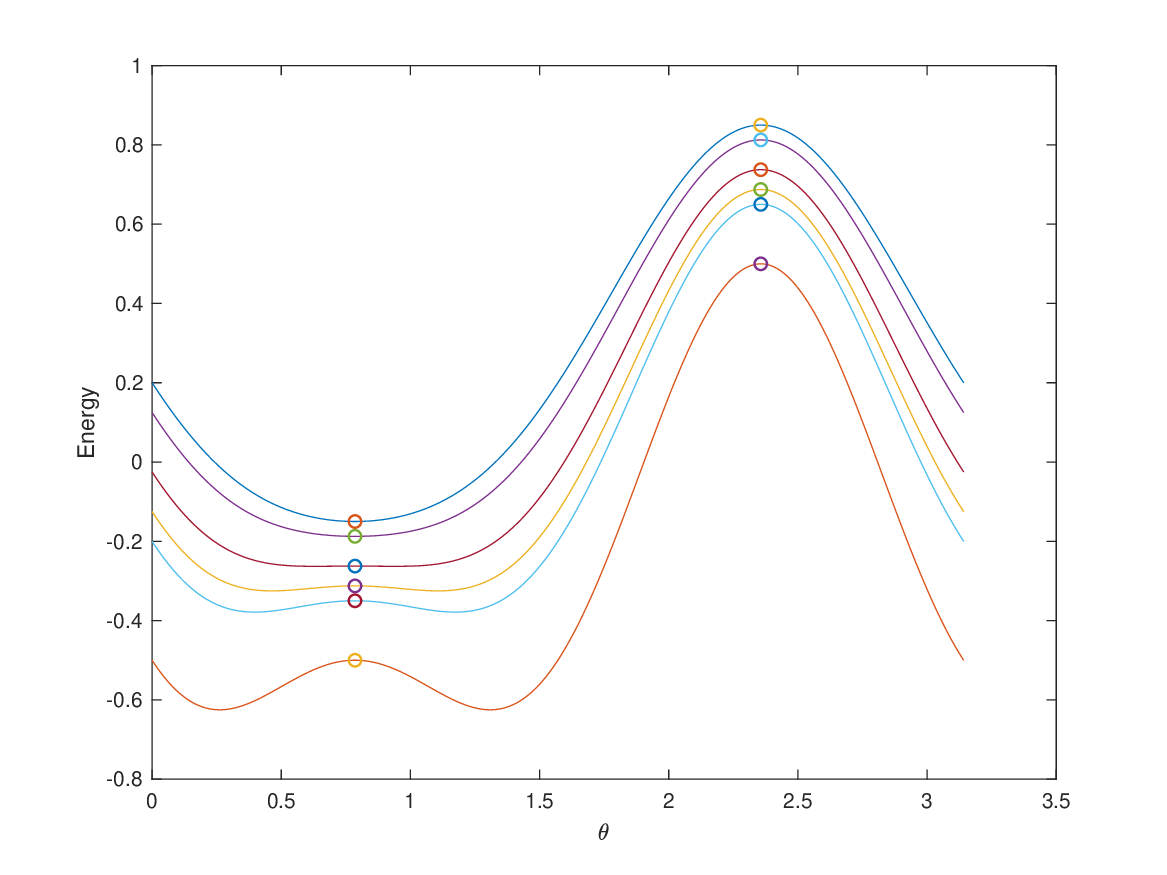}
\end{center}
\caption{Graph of $\theta\mapsto \mathscr H_1(\theta)$ for several values of $\kappa$
($\kappa\in\{ 1.2, 1.5, 2.1, 2.5, 2.8, 4\}) $. The circles corresponds to $(\pi/4,\mathscr H(\pi/4)$ and 
$(3\pi/4,\mathscr H(3\pi/4)$:
$\kappa=2$ is the threshold at which the convexity at $\pi/4$ changes 
}
\label{GrH1}
\end{figure}

Assuming the smallness condition 
\begin{equation}
  \label{small}
  0<\kappa<2,
  \end{equation}
we thus denote 
\begin{equation}
  \label{solstatH}
e^{i\omega t}U_*=\ds\frac{e^{i\omega t}}{\sqrt2}\begin{pmatrix}
1  \\ \tau
\end{pmatrix},\qquad \tau=\pm1
  \end{equation}
  the obtained  solution of \eqref{Hartree}, with 
$\tau=1$ corresponding to the state of minimal energy, and $\tau=-1$  
corresponding to the state of maximal energy.
Equivalently, we can consider
\begin{equation} \label{solstatH2}
X_*=\ds\frac{1}{\sqrt2}\begin{pmatrix}
1\\0\\
\tau\\0 
\end{pmatrix}\end{equation}
so that, given the extended rotation matrix
\[
R(\theta)=\begin{pmatrix}
\cos(\theta)&  -\sin(\theta)& 0 & 0
\\
\sin(\theta)& \cos(\theta)& 0 & 0
\\
  0 & 0&\cos(\theta)&  -\sin(\theta)
  \\
    0 & 0&\sin(\theta)& \cos(\theta)
\end{pmatrix},\]
 $R(\omega t)X_*$ defines a solution to \eqref{Hartreebis}.

 When $\kappa>2$, \eqref{solstatH} still both define a
  solution of \eqref{Hartree}, but for $\tau=1$ the solution
  does not reach the minimal energy.
   Moreover, in this situation we find two extra solutions 
   \begin{equation}\label{extra}
   \textrm{$e^{i\omega t}U_{*\kappa,\pm}$  with $U_{*\kappa,\pm}=\begin{pmatrix}\sin(\theta_\kappa^\pm)\\
    \cos(\theta_\kappa^\pm)
   \end{pmatrix}$ and 
   $
   \theta_\kappa^+=\ds\frac12\arcsin\Big(\ds\frac2\kappa\Big)$, $
   \theta_\kappa^-=\ds\frac\pi2-\ds\frac12\arcsin\Big(\ds\frac2\kappa\Big)$}.
    \end{equation}
      Since $\frac2\kappa>0$, both $\sin(\theta_\kappa^\pm)$ and 
   $\cos(\theta_\kappa^\pm)$ are positive.
   Using the elementary relation $\sin^2(\theta)=\frac{1-\cos(2\theta)}{2}$, 
   we can write
   \[
   \sin(\theta_\kappa^+)=\sqrt{\ds\frac{1-\cos(\arcsin(2/\kappa))}{2}}
   =\sqrt{\ds\frac{1-\sqrt{1-(2/\kappa)^2}}{2}}
   =\ds\frac12\Big(
   \sqrt{1+2/\kappa}-\sqrt{1-2/\kappa}\Big),
  \]
   and 
   \[
   \cos(\theta_\kappa^+)=\sqrt{1-\sin^2(\theta_\kappa)}=
   \ds\frac12\Big(
   \sqrt{1+2/\kappa}+\sqrt{1-2/\kappa}\Big), \]
   so that, with $\sin(\pi/2-\theta)=\cos(\theta)$, $\cos(\pi/2-\theta)=\sin(\theta)$, we can rewrite the solution $U_{*\kappa,\pm}$ as follows
   \[
   U_{*\kappa,+}=\ds\frac{1}{2\sqrt\kappa}\begin{pmatrix}
   \sqrt{\kappa+2}-\sqrt{\kappa-2}
   \\
\sqrt{\kappa+2}+\sqrt{\kappa-2})
   \end{pmatrix},\qquad
   U_{*\kappa,-}=\ds\frac{1}{2\sqrt\kappa}\begin{pmatrix}
   \sqrt{\kappa+2}+\sqrt{\kappa-2}
   \\
\sqrt{\kappa+2}-\sqrt{\kappa-2})
   \end{pmatrix}.
   \]
    The corresponding solution for  \eqref{Hartreebis} reads
\begin{equation}\label{solk2}
R(\omega t)X_*,\qquad 
X_*= \begin{pmatrix}
\sin(\theta_\kappa^\pm)
\\
0
\\
\cos(\theta_\kappa^\pm)
\\
0
\end{pmatrix}
=\ds\frac{1}{2\sqrt\kappa}
\begin{pmatrix}
\sqrt{\kappa+2}- \tau\sqrt{\kappa-2}
\\
0
\\
\sqrt{\kappa+2}+\tau\sqrt{\kappa-2}
\\
0
\end{pmatrix},\qquad \tau=\pm 1.
\end{equation}

   Going back to \eqref{EulerLag}, we find the Lagrange multiplier $\omega$ to be associated to all  these solutions. Namely we get 
    $$    -\omega U_{*0}=U_{*0}-U_{*1}-\kappa| U_{*0}|^2 U_{*0},\quad
 -\omega U_{*1}=U_{*1}-U_{*0}-\kappa| U_{*1}|^2 U_{*1}.$$
 Adding these relations and using $|U_{*0}|^2+|U_{*1}|^2=1$, we are led to 
 $$2(\omega+1)-\kappa = \ds\frac{U_{*1}}{U_{*0}}+\ds\frac{U_{*0}}{U_{*1}}=
\ds\frac{1}{U_{*1}U_{*0}}.$$
 Hence, we conclude that
\begin{align}\label{defom}
\textrm{for \eqref{solstatH}} &\quad \omega=\ds\frac\kappa2+\tau-1=\left\{
\begin{array}{ll}
\kappa/2,\quad & \text{if $\tau=+1$},
\\
-2+\kappa/2,\quad & \text{if $\tau=-1$},
\end{array}
\right.
\\
\label{defomextra}
\textrm{for \eqref{extra}} &\quad \omega=\kappa-1.
\end{align}

\subsection{Analysis of the coupled model}\label{sec:solscoupled}

Denoting $u_j=q_j+ip_j$ and $\varpi_j=\frac{\partial_t \psi_j}{2c^2}$, the energy functional \eqref{EnerSW} casts as
\begin{align}\label{EnerSWbis}
\mathscr H(X)
=&\ds\frac{|q_0-q_1|^2+|p_0-p_1|^2}2
+\ds\int_{\mathbb R^n}
\left(\ds 
c^2\ds\left(|\varpi_0|^2+\ds|\varpi_1|^2\right)
+\frac14\left(|\nabla\psi_0|^2+|\nabla\psi_1|^2
\right)\right)\ud z
\nonumber\\
&+\ds\frac12\ds\int_{\mathbb R^n}\sigma
(\psi_0 (|q_0|^2+|p_0|^2)+\psi_1 (|q_1|^2+|p_1|^2))
\ud z,
\end{align}
where $X$ is the shorthand notation for $(q_0,p_0,q_1,p_1,\psi_0,\varpi _0,\psi_1,\varpi_1)$.
Repeating the arguments used for the asymptotic model, we realize 
 that extrema of $\mathscr H$ can be found 
 by considering only the case $p_0=p_1=0$ and, taking into account the constraint of 
 normalized norm, $|u_0|^2+|u_1|^2=1$, we are led to 
 investigate the extrema of 
 \[\begin{array}{lll}
 \mathscr H_1(\theta,\psi_0,\varpi _0,\psi_1,\varpi_1)
&=&\ds\frac{|\cos(\theta)-\sin(\theta)|^2}2
+\ds\int_{\mathbb R^n}
\left(\ds 
c^2\ds\left(|\varpi_0|^2+\ds|\varpi_1|^2\right)
+\frac14\left(|\nabla\psi_0|^2+|\nabla\psi_1|^2
\right)\right)\ud z
\\
&&+\ds\frac12\ds\int_{\mathbb R^n}\sigma
(\psi_0 \cos^2(\theta)+\psi_1 \sin^2(\theta))
\ud z,
\end{array}\]
where $\theta$ lies in $[0,2\pi)$.
At the extrema, we infer that
\[\varpi_0=\varpi_1=0
 \]
 together with 
 \[
 -\Delta\psi_0=-\sigma \cos^2(\theta),\qquad -\Delta\psi_1=-\sigma \sin^2(\theta) 
 .\]
 The latter relation leads to 
 $\int_{\mathbb R^n}\sigma \psi_0\ud z=-\int_{\mathbb R^n}\sigma (-\Delta)^{-1}\sigma\ud z\cos^2(\theta)
 =-\kappa \cos^2(\theta)$, and similarly
 $\int_{\mathbb R^n}\sigma \psi_1\ud z
 =-\kappa \sin^2(\theta)$.
 Eventually, computing $\partial_\theta\mathscr H_1$ 
 yields \[
 -(\cos^2(\theta)-\sin^2(\theta))-
 \ds\frac12\ds\int_{\mathbb R^n}
 \sigma(\psi_0 - \psi_1)\sin(2\theta)\ud z
 .\]
 Therefore, at the extrema we obtain
 \[
 \ds\frac\kappa2\cos(2\theta)\left(\ds\frac2\kappa-\sin(2\theta)\right)=0.
 \]
 Hence, we find the same extrema as for the asymptotic model.

 In particular,  we set
 $Q_{*0}=\frac1{\sqrt 2}$, $P_{*0}=0$, $Q_{*1}=\frac \tau{\sqrt 2}$, $P_{*1}=0$,
 $\Psi_{*0}=\Psi_{*1}=-\frac{(-\Delta)^{-1}\sigma}2$,
 $\varpi_{*0}=\varpi_{*1}=0$,
and  the 
 energy is made minimal (resp. maximal) when $\tau=+1$ with $0<\kappa<2$ (resp. $\tau=-1$ without condition on $\kappa$).
 This analysis provides specific solutions of \eqref{S1}-\eqref{W1}, 
having the special form $(e^{i\omega t} U_{*0},e^{i\omega t} U_{*1},
 \Psi_{*0}, \Psi_{*1})$ where 
 $U_{*0}, U_{*1}$ are fixed complex numbers and 
  $\Psi_{*0}, \Psi_{*1}$ are fixed functions in $L^2(\mathbb R^n)$.
 This leads to the relations
  \[\begin{array}{ll}
  -\omega U_{*0}=U_{*0}-U_{*1}+U_{*0}\ds\int_{\mathbb R^n}\sigma \Psi_{*0}\ud z,\qquad &
    -\omega U_{*1}=U_{*1}-U_{*0}+U_{*1}\ds\int_{\mathbb R^n}\sigma \Psi_{*1}\ud z,
    \\[.4cm]
    -\Delta  \Psi_{*0}=-\sigma |U_{*0}|^2,\qquad &- \Delta  \Psi_{*1}=-\sigma |U_{*1}|^2.
   \end{array} \]
    Let $\Gamma$ denote the solution of $-\Delta \Gamma=\sigma$, which can be alternatively defined 
  by means of Fourier transform
  \[
  \Gamma=\mathscr F_{\xi\to z}^{-1}\left(\ds\frac{\widehat\sigma(\xi)}{|\xi|^2}
  \right).\]
  Hence, we get 
 $$\Psi_{*0}(z)=-  |U_{*0}|^2 \Gamma (z),\qquad \Psi_{*1}(z)=   -  |U_{*1}|^2 \Gamma (z),$$
 so that $U_{*0}, U_{*1}$ are required to satisfy 
 \[
(\omega+1) U_{*0}-U_{*1}-\kappa |U_{*0}|^2U_{*0}=0
 =(\omega+1 )U_{*1}-U_{*0}-\kappa |U_{*1}|^2U_{*1}
,\]
  together with  the physical normalisation
  $$|U_{*0}|^2+|U_{*1}|^2=1.$$
  With the extrema discussed above, we have
  $
  |U_{*0}|=|U_{*1}|=\frac{1}{\sqrt2}
  $ and for the system 
  \[
  \begin{pmatrix}
  \omega+1 -\kappa/2 & -1
  \\
  -1 & \omega+1 -\kappa/2
  \end{pmatrix}\begin{pmatrix}U_{*0}\\U_{*1}
  \end{pmatrix}=0
  \]
  to admit non trivial solutions, the  dispersion relation \eqref{defom}
  should be fulfilled.
  Given this condition, we conclude that  
  \begin{equation}\label{SpecSol}
  u_0(t)=\ds\frac{e^{i\omega t}}{\sqrt 2},  \quad 
  u_1(t)=\tau \ds\frac{e^{i\omega t}}{\sqrt 2},\quad 
  \psi_0(t,z)=-\ds\frac{\Gamma(z)}{2}, 
  \quad\psi_1(t,z)= - \ds\frac{\Gamma(z)}{2}\end{equation}
  satisfies \eqref{S1}-\eqref{W1}.
  \\
  
  If $\kappa>2$, we find two extra solutions which minimize the energy 
  \begin{equation}\label{extraSW}
  \begin{array}{l}
  Q_{*0}=\sin(\theta_\kappa^\pm),\quad Q_{*1}=\cos(\theta_\kappa^\pm),\quad P_{*0}=P_{*1}=0,\\
  U_{*}=\begin{pmatrix} Q_{*0}\\Q_{*1}\end{pmatrix},
  \Psi_{*0}=-|Q_{*0}|^2\Gamma,\quad\Psi_{*1}=-|Q_{*1}|^2\Gamma.
  \end{array}
  \end{equation}
With $\omega$  still given by \eqref{defomextra}, we conclude that 
 \begin{equation}\label{SpecSolextra}
  u_0(t)=e^{i\omega t}Q_{*0},  \quad 
  u_1(t)= e^{i\omega t}Q_{*1},\quad 
  \psi_0(t,z)=\Psi_{*0}, 
  \quad\psi_1(t,z)= \Psi_{*1}\end{equation}
  satisfies \eqref{S1}-\eqref{W1}.  
  
  Finally, we observe that the system can be expressed in the Hamiltonian formulation
  \[\partial_t X=
  \begin{pmatrix} 
  \mathscr J & 0
  \\
  0 & \mathscr J
  \end{pmatrix}
  \nabla_X\mathscr H(X),\qquad
  \mathscr J=\begin{pmatrix}
  0 & 1 & 0 & 0
  \\
  -1 & 0 & 0& 0
  \\
   0 & 0 & 0 & 1
  \\
0 & 0 &   -1& 0
  \end{pmatrix}.\]
  We shall see later on a more adapted  formulation, more convenient for the
    stability analysis.
    For the time being, this formulation makes
    the same connections between different viewpoints appear, as we did for the asymptotic equation.
    
    \subsection{Statement of the results}

    Let us collect here the main statements that will be obtained (definitions of the notions of stability will be made precise later on).
    
 \begin{theo}[Stability analysis for \eqref{Hartree}]
 Let us assume one of the following cases:
 \begin{itemize}
 \item[i)] 
 $\tau=-1$,
 \item[ii)]   $\tau=+1$ with $0<\kappa<2$, 
 \item[iii)]  $\kappa>2$.
\end{itemize}
We consider the reference solution of   \eqref{Hartree} given by 
 \eqref{solstatH} for i) and ii) or by  \eqref{extra} for iii).
 Then, the reference solution is spectrally and orbitally stable.
 \end{theo}

  \begin{theo}[Instability result for \eqref{Hartree}]\label{ThUnstHa}
  Let $\kappa>2$.  Then, the state $e^{i\omega t}(1/\sqrt2,1/\sqrt2)$ is a spectrally and  orbitally
unstable solution of \eqref{Hartree}.
  \end{theo}
  
     \begin{theo}[Stability analysis for \eqref{S1}-\eqref{W1}]
 Let  $\tau=+1$ with $0<\kappa<2$. Then, the reference solution  \eqref{SpecSol} of \eqref{S1}-\eqref{W1} is spectrally and orbitally stable.
 Let $\kappa>2$. Then, the reference solution \eqref{extraSW}-\eqref{SpecSolextra} is spectrally and orbitally stable.
 \end{theo}
 
\begin{theo}[Instability result for \eqref{S1}-\eqref{W1}]
Let $\tau=1$ with $\kappa>2$ or $\tau=-1$. Then the reference solution \eqref{SpecSol} is spectrally and orbitally unstable. 
\end{theo}
    
    These results are in line with the  analysis performed in \cite{SRN3} for plane waves solutions for the PDE system \eqref{SchW}
    and its asymptotic Hartree-like counterpart.
    It confirms 
    that the asymptotic model has more stable solutions than the original model, and that the dynamic coupling \eqref{W1}
    induces intricate and rich selection mechanisms.
    We expect this study will provide fruitful ideas to come back to  \eqref{SchW}
    set for  $x\in\mathbb R^d$, and will allow us to fill a gap in the understanding of open quantum systems.
    
  \section{Stability analysis of the asymptotic model \eqref{Hartree}}
  \label{AnHar}

   \subsection{Spectral and linearized stability}
   
   We start by linearizing \eqref{Hartree} about the solutions \eqref{solstatH}.
   We search for solutions of \eqref{Hartree} on the form 
  $$u_j=e^{i\omega t}(U_{*j}+v_j).$$
  Using $|u+h|^2=|u|^2+2\mathrm{Re}(\overline u h)+|h|^2$,  the dispersion relation \eqref{defom}, and 
  neglecting the non linear terms,  one is led to the following linearized system
  \begin{equation}\label{linHa}
  i\ds\frac{\ud}{\ud t}v_0=\tau v_0-v_1-\kappa \mathrm {Re}(v_0),\qquad
  i\ds\frac{\ud}{\ud t}v_1=\tau v_1-v_0-  \kappa \mathrm {Re}(v_1).
  \end{equation}
  We write $v_j=q_j+ip_j$, with $q_j,p_j$ real-valued. The unknown is now represented by the 
  vector $X=(q_0,p_0,q_1,p_1)$;  we get
  \[
\ds\frac{\ud}{\ud t} X=
  \mathbb LX,\qquad
  \mathbb L=
  \begin{pmatrix}
  0&  \tau & 0& -1
  \\
   \kappa-\tau & 0 & 1& 0
  \\
  0 & -1 & 0 & \tau
  \\
   1& 0 & \kappa- \tau & 0 
  \end{pmatrix}.
  \]
  The stability of this ODE system is related to the spectral analysis of the matrix $\mathbb L$: spectral stability means that the 
  real part of the eigenvalues of $\mathbb L$ are all non positive; 
  linearized stability means that any solution of this linear system remains uniformly bounded for any $t\geq 0$.
  
     \begin{proposition}\label{staLinHa}
   If $\tau=-1$, the system \eqref{linHa} is spectrally stable; if 
   $\tau=+1$, the system \eqref{linHa} is spectrally stable under the condition \eqref{small}.
   In these situations, if, moreover, $\mathrm{Re}(v_0+\tau v_1)\big|_{t=0}=0$, then, the solution 
 of   \eqref{linHa} remains uniformly bounded for any $t\geq 0$. 
 If  $\tau=+1$ with $\kappa>2$, the system is spectrally unstable.
   \end{proposition}
  
\noindent
{\bf Proof.}
    We observe that 0 is an eigenvalue of $\mathbb L$.
    Indeed, $\mathbb LX=0$ leads to the independent relations
    $$
    \left\{\begin{array}{l}
    \tau p_0=p_1,
    \\
     \tau p_1=p_0\end{array}\right.
     \textrm{ and } 
     \left\{\begin{array}{l}
    (\kappa-\tau) q_0=-q_1,
    \\
     (\kappa-\tau)  q_1=-q_0.\end{array}\right.
    $$
    Since $\tau^2 =1$, the former yields a non trivial solution, while
    the latter in general ($(\kappa-\tau)^2-1=\kappa(\kappa-2\tau)\neq 0$)
    has only the solution $q_0=q_1=0$.
    Hence we find  
    the eigenspace 
  $\mathrm{Ker}(\mathbb L)=\mathrm{Span}\{(0, 1, 0, \tau)\}$.
  Note however that $\mathbb L$ has a Jordan block associated to the eigenvalue 0, since 
  the kernel of $$\mathbb L^2
  =\begin{pmatrix}
  \kappa \tau -2 & 0 & 2\tau-\kappa & 0
  \\
  0 & \kappa \tau -2 & 0 & 2\tau-\kappa 
  \\
 2\tau -  \kappa & 0 & \tau\kappa -2 & 0
 \\
 0 &2 \tau - \kappa & 0 & \tau\kappa-2 
  \end{pmatrix}$$ is spanned by $\{(0, 1, 0, \tau), (1, 0, \tau, 0)\}$.
   This leads to solutions of \eqref{linHa}  the norm of which  can grow linearly.
  Next, let $\lambda\neq 0$, $X\neq 0$ satisfy
  $\mathbb LX=\lambda X$.
  Since $\tau^2=1$, we observe that $\tau q_0=-q_1$.
  Therefore, we obtain 
  $\lambda p_0=q_1-\tau q_0+\kappa q_0=(-2\tau+\kappa) q_0$,
  together with 
  $\lambda p_1= q_0+ (-\tau+\kappa) q_1
  =(2-\tau\kappa)q_0$.
  It yields
  $\lambda q_0=\tau p_0-p_1=- \big(\tau\frac{2\tau-\kappa}{\lambda}
  +\frac{2-\tau \kappa}{\lambda}\big)q_0$.
A non trivial solution $q_0$ exists provided $\lambda $ satisfies 
\[
\lambda^2  = -4 +2\tau\kappa.\]
  If $\tau=-1$, we find $\lambda=\pm 2i\sqrt{1+\kappa/2}$.
  If $\tau=1$, we find 
  $\lambda =\pm 2i\sqrt{1-\kappa/2}$, assuming the smallness condition 
  \eqref{small}; otherwise, $\lambda=\pm 2\sqrt{\kappa/2-1}$ and the system admits  a positive eigenvalue.   
  
  In fact, the problem \eqref{linHa} can be easily solved by hand. On the one hand,
  we have 
  $$\ds\frac{\ud}{\ud t}
  (q_0+\tau q_1)=0,\qquad 
  \ds\frac{\ud}{\ud t}
  (p_0+\tau p_1)=\kappa(q_0+\tau q_1)$$
  so that
  \[(q_0+\tau q_1)(t)=C_1,\qquad (p_0+\tau p_1)(t)=C_2+C_1\kappa t.\]
  On the other hand, the pair $(q_0-\tau q_1)$ and $(p_0-\tau p_1)$ solves a linear system
  associated to the matrix
  \[
  \begin{pmatrix}
  0 & 2\tau
  \\
  \kappa-2\tau & 0\end{pmatrix}
   \]
   which is diagonalizable with eigenvalues satisfying $\lambda^2=-4(1-\tau\kappa/2)<0$.
   The analysis of the linearized system is therefore complete.
   \QED

Similar computations 
can be performed with the solutions \eqref{extra}.
The linearized system now reads 
\begin{equation}
\label{linHaex}
i\ds\frac{\ud}{\ud t} v_0
=
(1+\omega-\kappa\alpha^2) v_0
-v_1  -2\kappa \alpha^2\mathrm{Re}(v_0)
,\qquad
i\ds\frac{\ud}{\ud t} v_1
=
(1+\omega-\kappa\beta^2) v_1
-v_0  -2\kappa \beta^2\mathrm{Re}(v_1),
\end{equation}
with 
\[
U_{*\kappa,\pm}=\begin{pmatrix}\alpha\\\beta\end{pmatrix}
,\qquad
\alpha=\ds\frac{\sqrt{\kappa+2}-\tau \sqrt{\kappa-2}}{2\sqrt\kappa}
=\sin(\theta_\kappa^\tau),\qquad
\beta=\ds\frac{\sqrt{\kappa+2}+\tau  \sqrt{\kappa-2}}{2\sqrt\kappa}=\cos(\theta_\kappa^\tau).\]
Let us set 
\[
A=
1+\omega-\kappa\alpha^2,\qquad
B=
1+\omega-\kappa\beta^2.
\]
Elementary manipulations lead to
\begin{equation}\label{notation3}
A=\ds\frac{\kappa}{2}+\tau \ds\frac{\sqrt{\kappa^2-4}}{2},\quad 
B=\ds\frac{\kappa}{2}-\tau \ds\frac{\sqrt{\kappa^2-4}}{2},\quad
AB=1,
\quad \kappa\alpha^2=B,\quad \kappa\beta^2=A.\end{equation}
The matrix associated 
to the linearized system thus reads
\[
\mathbb L=\begin{pmatrix}
0 & A& 0 & -1
\\
-A+2B&  0 & 1& 0
\\
0 & -1& 0 & B
\\1& 0& -B+2A& 0
\end{pmatrix}.\]
In turn, it can be checked that 
$$\mathrm{Ker}(\mathbb L)=\mathrm{Span}\{(0,1,0,A)\}.$$ 
Next, let $(\lambda,X)$ be an eigenpair of $\mathbb L$, with $\lambda\neq 0$.
We observe that $\lambda Aq_1=A(Bp_1-p_0)=-\lambda q_0$, which implies 
$Aq_1+q_0=0$. It follows that
\begin{equation*}
\lambda p_0=(-A+2B)q_0+q_1=(-A+2B)(-Aq_1)+q_1= A(A-B)q_1
\end{equation*}
and
\begin{equation*}
  \lambda p_1=(-B+2A)q_1+q_0=(-B+2A)q_1
-Aq_1= (A-B)q_1,
\end{equation*}
which lead to 
\[
\lambda q_1=Bp_1-p_0=B\ds\frac{ A-B}{\lambda}q_1
-\ds\frac{A(A-B)} {\lambda}q_1=-\frac{q_1}{\lambda}(A-B)^2.\]
Therefore, we obtain
\begin{align*}
  \lambda^2= -(A-B)^2=-(\kappa^2-4)=-\kappa^2+4<0.
\end{align*}
 We deduce that $\lambda\in i\mathbb R
$.

  \begin{proposition}\label{staLinHaEx}
   The system \eqref{linHaex} is spectrally stable.
   If, moreover, $\mathrm{Re}(v_0+A v_1)\big|_{t=0}=0$, then, the solution 
 of   \eqref{linHaex} remains uniformly bounded for any $t\geq 0$. 
   \end{proposition}

\noindent
{\bf Proof.} 
The spectral stability has just been established above, all eigenvalues of $\mathbb L$ being with a non positive real part.
Next, we introduce the vectors
\[\Psi=(1,0,A,0), \qquad
\Psi_1=\Big(0,\ds\frac{-\tau}{2\sqrt{\kappa^2-4}},0,
\frac{ \kappa\tau +\sqrt{\kappa^2-4}}{4\sqrt{\kappa^2-4}}\Big).\]
They satisfy 
\[
\mathbb L^\intercal\Psi=0,\qquad  \mathbb L^\intercal\Psi_1=\Psi.\]
Let $X$ satisfy $\frac{\ud}{\ud t}X=\mathbb LX$.
We observe that $\frac{\ud}{\ud t}X\cdot \Psi=
\frac{\ud}{\ud t}(q_0+Aq_1)=X\cdot \mathbb L^\intercal \Psi=0$, and
$\frac{\ud}{\ud t}X\cdot \Psi_1=X\cdot \mathbb L^\intercal \Psi_1=X\cdot\Psi
$. Hence $X(t)\cdot \Psi=X_{\mathrm{init}}\cdot \Psi$ is conserved and 
$X(t)\cdot \Psi_1=X_{\mathrm{init}}\cdot \Psi_1+tX_{\mathrm{init}}\cdot \Psi$
grows at most linearly.
Assuming $X_{\mathrm{init}}\cdot\Psi=0$ prevents the linear growth.
Finally, the pair $(Ap_0-p_1,Aq_0+q_1)$ satisfies
the $2\times 2 $ system governed by the matrix
\[
\begin{pmatrix}0& (B-A)
\\
(A-B)& 0\end{pmatrix}\]
the eigenvalues of which are clearly purely imaginary.
These observations completely characterize the solution of the linear system \eqref{linHaex}.
\QED

Propositions~\ref{staLinHa} and \ref{staLinHaEx} are illustrated in Fig.~\ref{stalin}
where we perform simulations of the different scenario:
the stable case ((a)-(b)) requires a condition on both the coefficients ($\tau$, $\kappa$) 
and the data;  when the    orthogonality condition of   Proposition~\ref{staLinHa} is violated, one observes a linear 
growth of the $L^2$ norm ((c)-(d)); when the condition on the data is not fulfilled, one observes an exponential blow up ((e)-(f)).
\\

System \eqref{Hartree} is a mere finite dimensional differential system.
As far as one is concerned with the stability of equilibrium solution of differential systems in finite dimension, spectral stability 
implies non linear stability, see e.~g.~\cite[Prop.~1.41]{Tao2}, \cite[Th. 1.1 \& 1.2]{Str}.
Here, we are dealing with the notion of \emph{orbital stability}, and the reference solutions remains time-dependent which induces some subtleties. 
We shall detail approaches which do not use properties specific to the finite dimensional framework, having in mind more complicated couplings.

     \begin{figure}[!htpb]
     \begin{subfigure}{0.48\textwidth}
  \includegraphics[height=5cm]{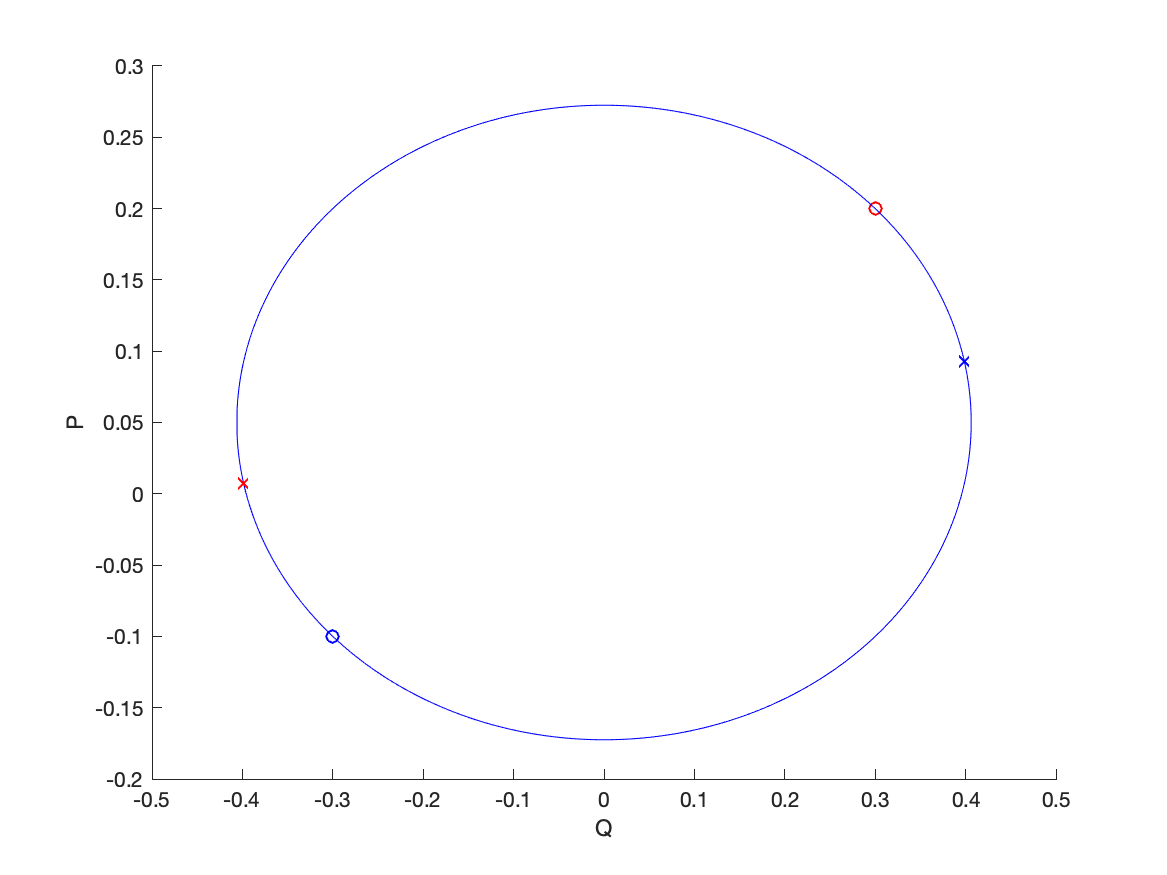}
  \caption{}
  \end{subfigure}
   \begin{subfigure}{0.48\textwidth}
   \includegraphics[height=5cm]{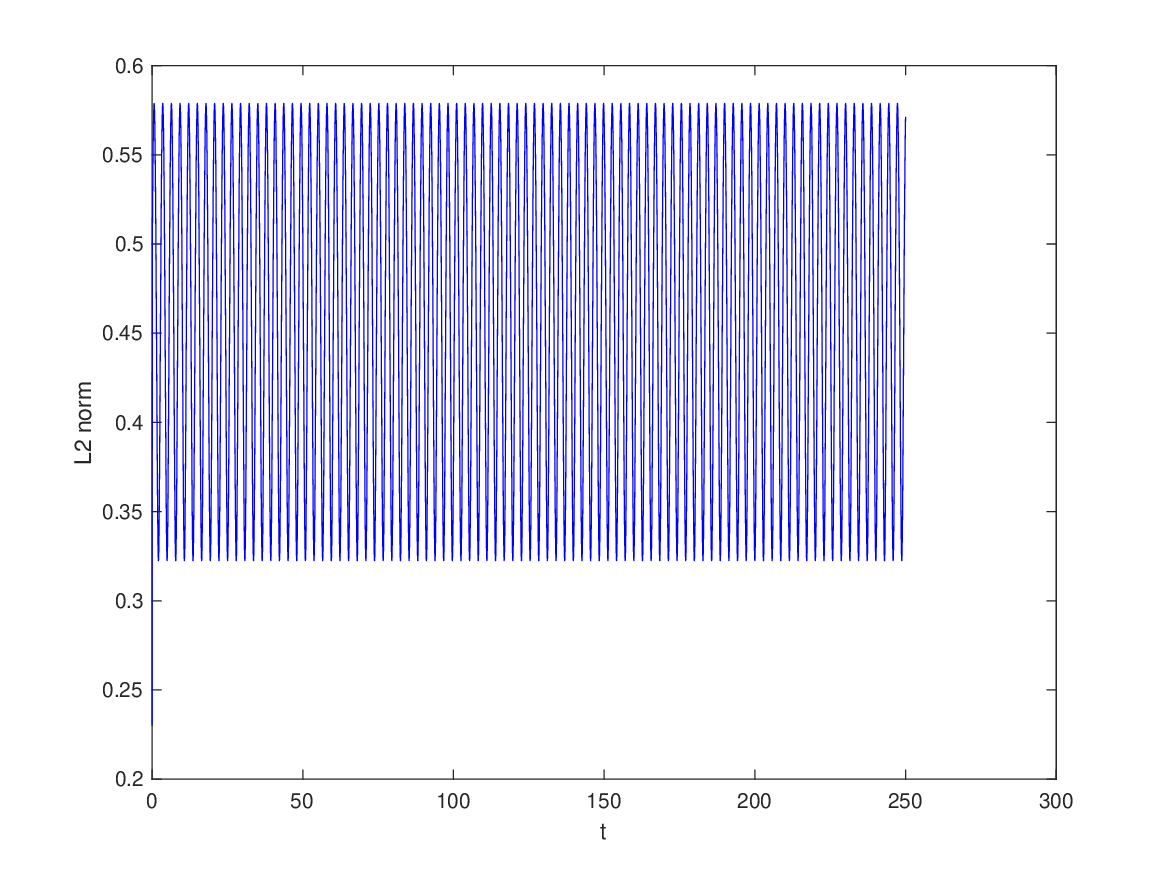}
     \caption{}
    \end{subfigure}
   \begin{subfigure}{0.48\textwidth}
   \includegraphics[height=5cm]{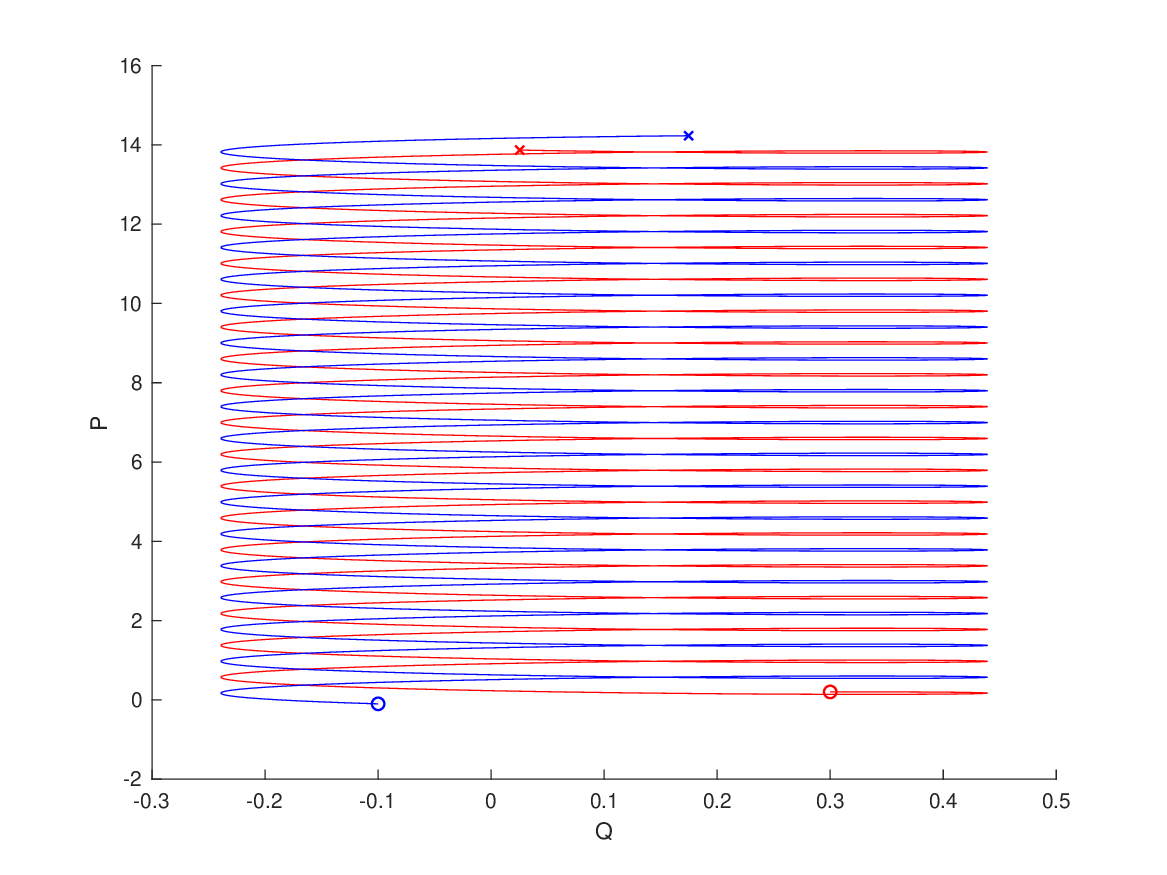}
     \caption{}
    \end{subfigure}
   \begin{subfigure}{0.48\textwidth}
   \includegraphics[height=5cm]{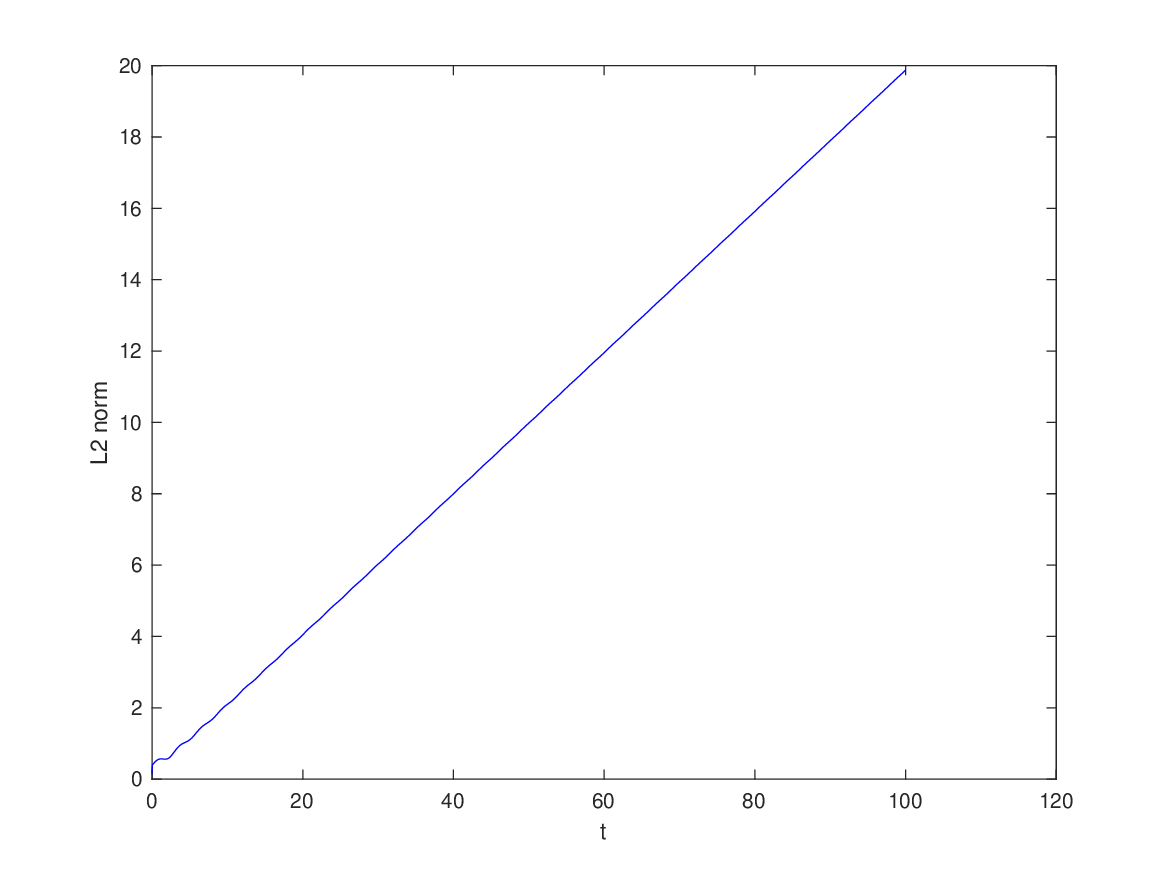}
     \caption{}
 \end{subfigure}
   \begin{subfigure}{0.52\textwidth}
    \includegraphics[height=5cm]{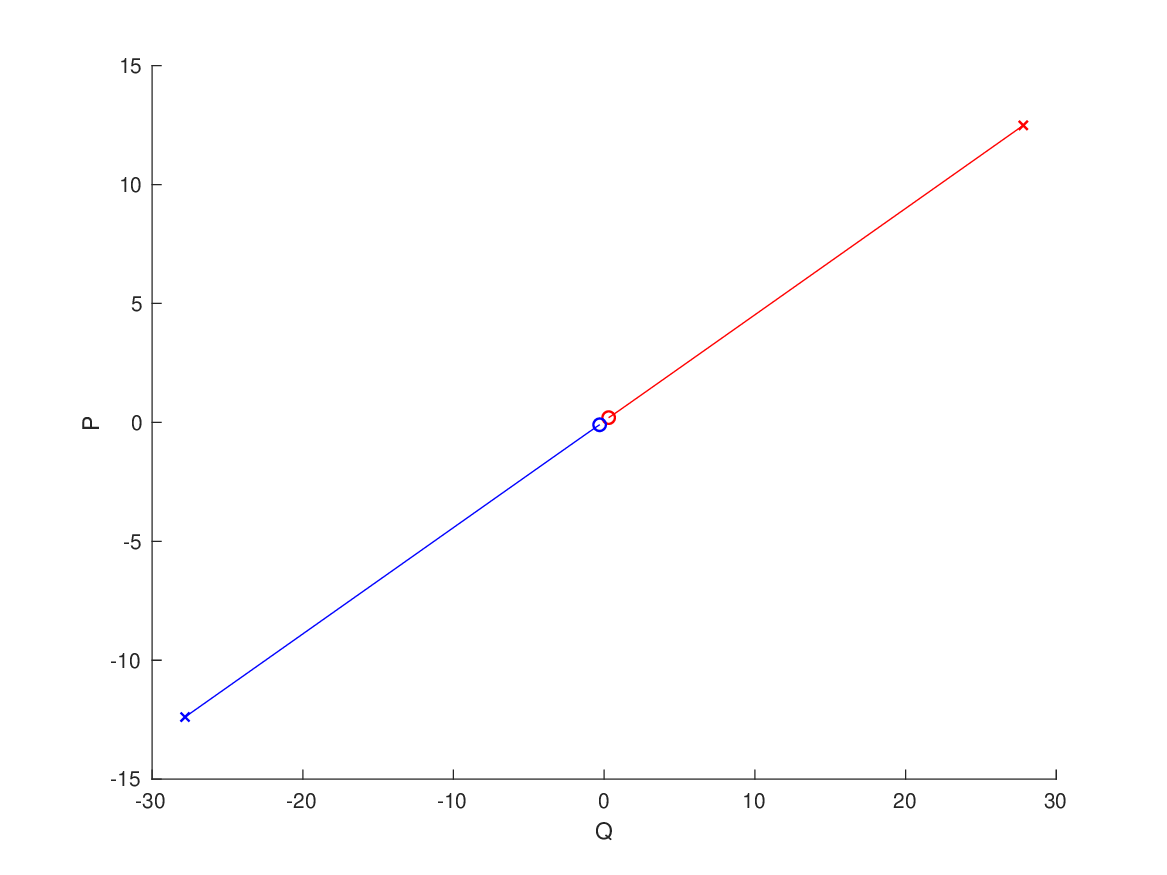}
      \caption{}
     \end{subfigure}
   \begin{subfigure}{0.48\textwidth}
   \includegraphics[height=5cm]{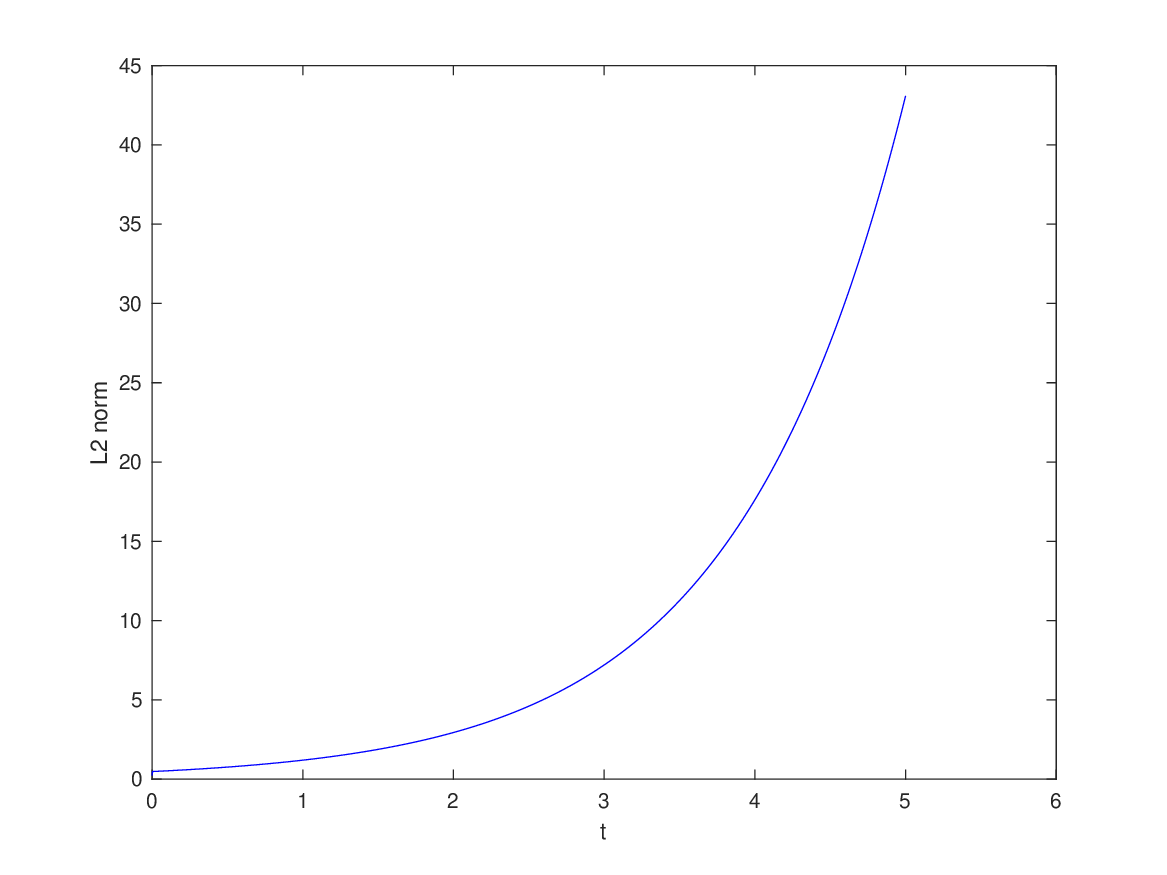}
     \caption{}
    \end{subfigure}
    \caption{Simulation of the  linearized asymptotic model \eqref{linHa}.
  The circled points indicate the initial state, the cross indicate the final state.
 (a)-(b): stable case $\kappa=1.4$ and $\tau=+1$; phase portrait at $T=250$ (a) and evolution of the $L^2$ norm (b) for a well prepared data. 
     The solution remains in a bounded domain. Similar results can be obtained 
   when $\tau=-1$ or, with $\kappa>2$, for  the linearized problem \eqref{linHaex}.
   (c)-(d):  $\kappa=1.4$ and $\tau=+1$ with ill prepared data; phase portrait at $T=100$ (c) and evolution of the $L^2$ norm (d);
   the $L^2$  norm of the solution grows linearly.
   (e)-(f): instable case  $\kappa=2.4$ and $\tau=+1$; phase portrait at $T=50$ (e) and evolution of the $L^2$ norm (f)
  } \label{stalin}
  \end{figure}
  
%
%
%
%
\subsection{Orbital stability}

  Let us  set $F(X)=\frac{|X|^2}{2}=\frac{Q_0^2+Q_1^2+P_0^2+P_1^2}{2}$ and 
  introduce the functional
  \[\mathscr E (X)=\mathscr H(X)+\omega F(X)\]
  with $\mathscr H$ defined by \eqref{EnerHaBis}. This quantity is thus conserved by the dynamical system \eqref{Hartreebis}.
  We observe that \eqref{EulerLag} can be reformulated as 
  \begin{equation}\label{0gradE}\nabla\mathscr E(X_*)=0\end{equation}
  and $\mathscr L$ corresponds to the Hessian of $\mathscr E$ evaluated at $X_*$.
 Inspired by the strategy described in \cite{SRN1}, we introduce the level set of solution if \eqref{Hartreebis} associated to $X_*$,
  \[\mathscr S=\{X\in \mathbb R^4,\ F(X)=F(X_*)=1/2\}.
  \]
  We wish to establish a coercivity estimate, on a certain subspace,
  for the quadratic form $ X\mapsto \mathscr LX\cdot X$. This is a crucial property for establishing the orbital stability, an idea that dates back to  \cite{Weinstein1,Weinstein2} for Schr\"odinger equations, see \cite{ SRN1,GSS,Tao}.
   \\
  
  With $X_*$ given by \eqref{solstatH2}, 
   the tangent set to the level set  is given by
  \[
  T\mathscr S=\{X\in \mathbb R^4, \ \nabla F(X_*)\cdot X=0\}=
  \{(q_0,p_0,q_1,p_1)\in \mathbb R^4, \ q_0+\tau q_1=0\}.
  \]
  The orbit associated to $X_*$ is given by
  \[\mathscr O=\Big\{\frac1{\sqrt 2}(\cos(\theta), \sin(\theta), \tau\cos(\theta),\tau \sin(\theta)),\  \theta\in \mathbb R\Big\}.\]
  and 
  we get  \[(T\mathscr O)^\perp=\{(q_0,p_0,q_1,p_1)\in \mathbb R^4, \ p_0+\tau p_1=0\}.\]
  The reference solution associated to $X_*$ is said orbitally stable 
  if, for any $\epsilon>0$, there exists $\delta>0$, such that, for any solution $t\mapsto Y(t)$ of \eqref{Hartreebis},   
   $|Y(0)-X_*|\leq \delta$ implies that $\mathrm{dist}(Y(t)-\mathscr O)\leq \epsilon$ holds for any $t\geq 0$.
 
 \begin{rmk}\label{phase_inv}
 Bearing in mind the transformation \eqref{notation}, 
 multiplying the components of $U\in \mathbb C^2$ by $e^{i\theta}$ is equivalent to 
 apply the (extended) rotation $R(\theta)$ to $X\in \mathbb R^4$, with leaves the energy $\mathscr H(X)$, as well as $\mathscr E(X)$, invariant.
 The identity  $\mathscr H(R(\theta)X)=\mathscr H(X)$ yields $R(\theta)^\intercal \nabla\mathscr H(R(\theta)X)=\nabla\mathscr H(X)$
 and we observe that $R(\theta)^{-1}R'(\theta)=-\mathscr J$.
 These observations allow us to derive directly the linearized system: with $\partial_t X=\mathscr J\nabla \mathscr H(X)$ and $X(t)=
 R(\omega t)(X_*+\tilde X(t))$, we get
 \[\partial_t \tilde X=\omega \mathscr J(X_*+\tilde X)+ \mathscr J\nabla \mathscr H(X_*+\tilde X)
 .\]
 Assuming the perturbation to be small, at leading order the right hand side reads
 \[ \mathscr J(\omega X_* +\nabla \mathscr H(X_*))+ \mathscr J
 (\omega\tilde X+ 
 D^2 \mathscr H(X_*)\tilde X)
 =0+\mathscr J\mathscr L\tilde X=\mathbb L\tilde X.
 \]
 \end{rmk}

  In order to investigate the  orbital stability of the system, we recast 
the linearized system by using  the symplectic form
  \[
  \mathbb L=
  \underbrace{\begin{pmatrix}
  0 & 1 & 0 & 0
  \\
    -1 & 0 & 0& 0
    \\
     0 & 0& 0 & 1
     \\
     0 & 0& - 1 &0
     \end{pmatrix}}_{=
  \mathscr J}
  \underbrace{\begin{pmatrix}
  \tau-\kappa & 0 & -1 & 0
  \\
    0 & \tau & 0& -1
    \\
     -1 & 0& \tau-\kappa & 0
     \\
     0 & -1&  0& \tau
     \end{pmatrix}}
  _{\mathscr L},
  \]
  with $\mathscr L=D^2\mathscr E(X_*)$ symmetric. 
  
  \begin{lemma}\label{specL}
  The spectrum of the matrix $\mathscr L$ is 
$\sigma(\mathscr L)=\{0,-\kappa,2\tau,2\tau-\kappa\}$
with  eigenspaces spanned respectively by 
\[\begin{array}{ll}
X_0=(0,1,0,\tau),\qquad&
X_{-\kappa}=(1,0,\tau,0),\\
X_{2\tau-\kappa}=(1,0,-\tau,0),\qquad & 
X_{2\tau}=(0,1,0,-\tau)
.\end{array}
.\]
\end{lemma}
  
  Hence, we get
  \[\mathscr LX\cdot X=(\tau-\kappa)(Q_0^2+Q_1^2)
  -2 Q_1Q_0
  +\tau (P_0^2+P_1^2)-2P_1P_0.
  \]
  As a matter of fact, when $\tau=1$, it recasts as
   \[
    \mathscr LX\cdot X=|P_0-P_1|^2+|Q_0-Q_1|^2
    -\kappa(Q_0^2+Q_1^2).\]
    Restricting to the subspace $T\mathscr S\cap (T\mathscr O)^\perp$, we have 
    $Q_0=-\tau Q_1$ and $P_0=-\tau P_1$, so that, still for $\tau=1$, we get
    \[
    \mathscr LX\cdot X=4|P_0|^2+ 2(2-\kappa)|Q_0|^2\geq \ds(2-\kappa)|X|^2
    .\]
%
This coercivity estimate 
is key in establishing the orbital stability.
Surprisingly, the case  $\tau=-1$ is simpler.
We now work with 
$$\mathscr E(X)=-\mathscr H(X)-\omega F(X).$$
 We still have $\nabla \mathscr E(X_*)=0$ and $D^2\mathscr E(X_*)=-\mathscr L$. The spectral decomposition of  $\mathscr L$ implies that $-\mathscr L$ is coercive on $(\mathrm{Ker}(\mathscr L))^\perp=(T\mathscr O)^\perp$.
 This allows us to justify the orbital stability.
 \\

We turn to the case where $\kappa>2$  and $X_*=(\alpha,0,\beta,0)$ is given by \eqref{solk2}. 
Now, we  look at 
\[
\mathscr  L=
\begin{pmatrix}
0 & -1& 0 & 0
\\
1& 0& 0 & 0
\\
 0 & 0& 0 & -1
 \\
 0 & 0&
1& 0
\end{pmatrix}\mathbb L=
\begin{pmatrix}
A-2B&  0 & -1& 0
\\
0 & A& 0 & -1
\\-1& 0& B-2A& 0
\\
0 & -1& 0 & B
\end{pmatrix}.
\]
The equations for the eigenpairs uncouple since we get
\[\begin{array}{ll}
(A-\lambda)p_0=p_1,\qquad &(B-\lambda)p_1=p_0,
\\
(A-2B-\lambda)q_0=q_1,\qquad &(B-2A-\lambda)q_1=q_0.
\end{array}
\]
The former leads to 
\[
\lambda(\lambda-(A+B))=\lambda(\lambda-\kappa)=0,
\]
and the latter gives
\[
(B-2A-\lambda)(A-2B-\lambda)-1
=\lambda^2+\lambda(A+B)+(A-2B)(B-2A)-1
=
\lambda^2+\lambda\kappa-2(\kappa^2-4)=0.
\]
This gives the eigenelements of $\mathscr L$.

\begin{lemma}\label{specL2}
We have
$$\sigma(\mathscr L)=\Big\{0,\kappa,\ds\frac{-\kappa+\sqrt{9\kappa^2-32}}{2},\ds\frac{-\kappa-\sqrt{9\kappa^2-32}}{2}\Big\},$$
where only the last value is negative,
with eigenspaces spanned respectively by 
$$\begin{array}{ll}
X_0=\Big(0,1,0,\ds\frac\kappa2+\tau \frac{\sqrt{\kappa^2-4}}{2}\Big),\qquad & X_\kappa=\Big(0,1,0,-\ds\frac\kappa2+\tau \frac{\sqrt{\kappa^2-4}}{2}\Big),
\\[.4cm]
X_+=\Big(1,0,\tau \ds\frac32\sqrt{\kappa^2-4}-\ds\frac12\sqrt{9\kappa^2-32},0\Big),\qquad &
X_-=\Big(1,0,\tau\ds\frac32\sqrt{\kappa^2-4}+\ds\frac12\sqrt{9\kappa^2-32},0\Big).
\end{array}$$
\end{lemma}
Establishing the orbital stability amounts to check the coercivity of $\mathscr L$ on 
$T\mathscr S\cap(T \mathscr O)^\perp$, where, now, 
\[T\mathscr S=\big\{X=(q_0,p_0,q_1,p_1)\in \mathbb R^4,\ X\cdot X_{*}=\alpha q_0+\beta q_1=0\big\},\]
and
\[(T\mathscr O)^\perp=\big\{X=(q_0,p_0,q_1,p_1)\in \mathbb R^4,\ \alpha p_0+\beta p_1=0\big\}.\]
We have
\[
\mathscr LX\cdot X=(A-2B)q_0^2-2q_0q_1+
(B-2A)q_1^2+Ap_0^2-2p_0p_1+Bp_1^2.\]
Since $AB=1$ and $\frac\alpha\beta=B$, on $T\mathscr S\cap(T \mathscr O)^\perp$, it reduces to
\[
\mathscr LX\cdot X\big|_{T\mathscr S\cap(T \mathscr O)^\perp}
=\big(A+(B-2A)B^2\big) q^2_0+
\big(A+B^3+2B\big)p_0^2
.\]
A tedious, but elementary,  computation yields
\[
\mathscr LX\cdot X\big|_{T\mathscr S\cap(T \mathscr O)^\perp}
=\ds\frac{\kappa-\tau  \sqrt{\kappa^2-4}}{2}((\kappa^2-4)q_0^2
+\kappa^2 p_0^2), 
\] 
hence the desired coercivity estimate holds.

\subsection{Symplectic formulation and further comments about spectral stability}

Let 
 us keep focused on the spectral stability issue.
 For the problem \eqref{Hartree}, the spectrum of $\mathbb L=\mathscr J\mathscr L$ is completely determined, as seen above, and we have directly a full understanding of the linearized problem. 
 However, for more intricate system, like \eqref{S1}-\eqref{W1}, we do not have a direct access to the spectrum of $\mathbb L$. The strategy is to deduce information about stable/instable modes 
  from the study of $\mathscr L$ which could be easier (in particular because $\mathscr L$ is symmetric).
  To this end, let us introduce the auxilliary operators 
\[\mathscr M=-\mathscr J\mathscr L\mathscr J,\qquad
\mathbb A=\mathscr P\mathscr M\mathscr P,\]
where $\mathscr P$ is the orthogonal projection on $(\mathrm{Ker}(\mathscr L))^\perp$.
We  also introduce 
 $$\mathbb K=\mathscr P\mathscr L^{-1}\mathscr P.$$
The counting of the eigenvalues of $\mathbb L$ is based on the following considerations. 
We are interested in the coupled system 
\begin{equation}\label{GenEi}\mathscr MX
=-\lambda\tilde X,\qquad \mathscr L\tilde X=\lambda X.\end{equation}
It turns out that this problem \eqref{GenEi} admits non trivial solutions iff 
$\pm\lambda$  are eigenvalues of $\mathbb L$.
Next, \eqref{GenEi} admits non trivial solution with $\lambda\neq 0$, iff the generalized eigenvalue problem 
\begin{equation}\label{GenEi2}
\mathbb AX
=\mu \mathbb K X\end{equation}
(which recasts as $\mathscr M X=\mu \tilde X$, $\mathscr L\tilde X=X$, with $X\in (\mathrm{Ker}(\mathscr L))^\perp$)
admits non trivial solutions with $\mu=-\lambda^2$.
The spectral stability means that  the spectrum of $\mathbb L$ is contained in $i\R$.  
This can be reformulated as saying
that all the eigenvalues of the generalized eigenproblem \eqref{GenEi2} are real and positive. 
In order  to count the eigenvalues $\mu$ of the generalized eigenvalue problem, we define the following quantities: 
\begin{itemize}
\item $N^-_n$, the number of negative eigenvalues,
\item $N^0_n$, the number of  eigenvalues zero,
\item $N^+_n$, the number of positive eigenvalues,
\end{itemize}
counted with their algebraic multiplicity, the eigenvectors of which are associated to 
non-positive 
values of the the quadratic form $X\mapsto (\mathbb KX|X)=(\mathscr L^{-1}\mathscr PX|\mathscr PX)$.
 Moreover, let
$N_{C^+}$ be the number of  generalized eigenvalues $\mu\in \mathbb C$ of \eqref{GenEi2} with $\mathrm{Im}(\mu)>0$.
As said above, the eigenvalues counted by $N^-_n$ and $N_{C^+}$ correspond to cases of instabilities for the linearized problem. 
We now use the counting argument of \cite[Theorem~1]{ChouPel} (see also the review \cite{LZ}) which asserts that
\[N^-_n+N^0_n+N^+_n+N_{C^+}=n(\mathscr L),\]
    the number of negative eigenvalues of $\mathscr L$.
    Let us check how this counting machinery works for \eqref{Hartree}.
  \\

  Let us begin with the case where $X_*$ is given by \eqref{solstatH2}.
We use the notation of Lemma~\ref{specL}.
For further purposes, we remark that 
\[
\mathscr J X_{-\kappa}=-X_0,\qquad \mathscr JX_{2\tau-\kappa}=-X_{2\tau}.\]
In particular, for $\tau=-1$, $\mathscr L$ has three negative eigenvalues;
for $\tau=+1$ and assuming \eqref{small}, 
there are two  positive eigenvalues and one negative eigenvalue
but if  $\tau=+1$ and  \eqref{small} is violated, 
there are one  positive eigenvalue and two negative eigenvalues.
Note that 
\begin{itemize}
\item [e1)] the eigenvectors  $X_0,X_{-\kappa}, X_{2\tau-\kappa}, X_{2\tau}$ form a  orthogonal basis of $\mathbb R^4$;
\item [e2)] with $X_*=\frac{1}{\sqrt2}(1, 0 , \tau, 0)=\frac{X_{-\kappa}}{\sqrt2}$ the reference solution, we have
 \[X_*\cdot  X_0=X_*\cdot X_{2\tau-\kappa}= X_*\cdot  X_{2\tau}=0;\]
 \item [e3)] and $X_*\cdot X_{-\kappa}=\sqrt 2>0.$
\end{itemize}

We start by computing $N_n^0=1$. We have seen that $\mathrm{Ker}(\mathscr L)$ is spanned by $X_0=(0,1,0,\tau)$. Hence, we have to solve $\mathscr L\tilde X_0=Y_0 $ with $Y_0=-\mathscr JX_0=(-1, 0, -\tau, 0)$ and $\tilde X_0\in (\mathrm{Ker}(\mathscr L))^\perp$. This leads to  
$\tilde X_0=\frac1\kappa(1,0,\tau,0)$ which yields 
 $\mathbb K Y_0\cdot Y_0=\mathscr L^{-1} Y_0\cdot Y_0 = \tilde X_0\cdot Y_0=-\frac2\kappa<0$
and thus $N_n^0=1$.

Next, solving the generalized eigenvalue problem amounts to solve
\[\begin{array}{ll}
-\tilde q_1+\tau \tilde q_0-\kappa \tilde q_0=q_0,
\qquad \qquad&
\tau q_0-q_1=\mu \tilde q_0,
\\
-\tilde q_0+\tau \tilde q_1-\kappa \tilde q_1=q_1,
\qquad\qquad &
\tau q_1-q_0=\mu \tilde q_1,
\\
\tau  p_0- \kappa p_0-p_1=\mu \tilde p_0,
\qquad\qquad&
\tau \tilde p_0-\tilde p_1=p_0,
\\
-p_0+\tau  p_1- \kappa p_1=\mu \tilde p_1,
\qquad\qquad&
\tau \tilde p_1-\tilde p_0=p_1,
\end{array}\]
with  $X=(q_0,p_0,q_1,p_1), \tilde X=(\tilde q_0,\tilde p_0,\tilde q_1,\tilde p_1)\in  (\mathrm{Ker}(\mathscr L))^\perp$.
We set
\begin{equation}\label{defmk}
M_\kappa=\begin{pmatrix}
\tau -\kappa & -1
\\
-1 & \tau -\kappa\end{pmatrix}.\end{equation}
The $q$ and $p$ equations decouple and we have, on the one hand
\[
M_\kappa\begin{pmatrix}\tilde q_0\\\tilde q_1\end{pmatrix}=
\begin{pmatrix}q_0\\q_1\end{pmatrix},\qquad
M_0\begin{pmatrix} q_0\\q_1\end{pmatrix}=\mu
\begin{pmatrix}\tilde q_0\\\tilde q_1\end{pmatrix},\]
and, on the other hand
\[M_\kappa\begin{pmatrix}p_0\\ p_1\end{pmatrix}=\mu
\begin{pmatrix}\tilde p_0\\ \tilde p_1\end{pmatrix},\qquad
M_0\begin{pmatrix} \tilde p_0\\ \tilde p_1\end{pmatrix}=
\begin{pmatrix}p_0\\p_1\end{pmatrix}.\]
It amounts to say that $(\tilde q_0,\tilde q_1)$ and 
$(\tilde p_0,\tilde p_1)$ are eigenvectors for $\mu$ of $M_0 M_\kappa$ and 
$M_\kappa M_0$, respectively. Here, we get
\[
M_\kappa M_0=
\begin{pmatrix}
2-\tau \kappa & \kappa-2\tau
\\
\kappa-2\tau & 2-\tau \kappa\end{pmatrix}=
M_0M_\kappa,\]
the eigenvalues of which being $0$ and $4(1-\tau \kappa/2)$.
We thus obtain the solutions $\tilde X_1=(1,0,-\tau,0)$ and $\tilde X_2=(0,1,0,-\tau)$, associated to 
$X_1=\mathscr L\tilde X_1=(2\tau -\kappa, 0 , \kappa\tau-2, 0)$, $X_2=\mathscr L\tilde X_2=(0,2\tau ,0,-2)$
which both belong to $(\mathrm{Ker}(\mathscr L))^\perp$.
We compute $\mathscr L^{-1}X_1 \cdot X_1=\tilde X_1\cdot X_1=
 2(2\tau -\kappa)$, 
which  is negative when $\tau =-1$ and has the sign of $2-\kappa$ when $\tau=+1$, and  
$\mathscr L^{-1}X_2\cdot X_2=\tilde X_2\cdot X_2=4\tau$.
Therefore, we can verify the counting formula in the following three cases
\begin{itemize}
\item $\tau=-1$: $n(\mathscr L)=3$ and
$N^0_n=1$, $N^+_n=2$, $N^-_n=0$, which yields $N_{C^+}=0$
and indeed we found that $\mathbb L$ has two purely imaginary eigenvalues, 
there is no exponentially unstable solution to the linearized system;
\item $\tau=1$ and $\kappa>2$: $n(\mathscr L)=2$
and
$N^0_n=1$, $N^+_n=0$, $N^-_n=1$, which yields $N_{C^+}=0$ 
and indeed we found that $\mathbb L$ has two real eigenvalues,
we can find exponentially unstable solutions to the linearized system;
\item $\tau=1$ and $0<\kappa<2$: $n(\mathscr L)=1$
and
$N^0_n=1$, $N^+_n=0$, $N^-_n=0$, which yields $N_{C^+}=0$ and indeed we found that $\mathbb L$ has two purely imaginary eigenvalues, there is no exponentially unstable solution to the linearized system.
\end{itemize}

We can perform similar computations for the solution \eqref{solk2}. We now use the notation of Lemma~\ref{specL2}.
We have seen that $\mathrm{Ker}(\mathscr L)$ is spanned by $X_0=(0,1,0,A)$.
We start by solving $\mathscr L\tilde X_0=Y_0 $ with $Y_0=-\mathscr JX_0=(-1, 0, -A, 0)$ so that 
$\tilde X_0=\frac{1}{2(A-B)}(-1,0,A,0)$ which yields $\tilde X_0\cdot Y_0=
\mathbb K Y_0\cdot Y_0=
- \frac A2<0
$, and thus $N_n^0=1$.
Solving the generalized eigenvalue problem amounts to solve
\[\begin{array}{l}
M\begin{pmatrix}\tilde q_0\\\tilde q_1\end{pmatrix}=\mu
\begin{pmatrix}\tilde q_0\\\tilde q_1\end{pmatrix},\qquad
M^\intercal \begin{pmatrix}\tilde p_0\\\tilde p_1\end{pmatrix}=\mu
\begin{pmatrix}\tilde p_0\\\tilde p_1\end{pmatrix},
\\
M=\begin{pmatrix}A & -1 \\
-1 & B\end{pmatrix}
\begin{pmatrix}A-2B & -1 \\
-1 & B-2A\end{pmatrix}=
\begin{pmatrix}A^2-1 & A-B \\
B-A & B^2-1\end{pmatrix},
\end{array}\]
the eigenvalues of $M$ being $0$ and $\kappa^2-4>0$ (thus $N_n^-=0$).
We thus obtain the solutions $\tilde X_1=(1,0,-B,0)$
%
and $\tilde X_2=(0,1,0,B)$.
Accordingly, we get $X_1=\mathscr L\tilde X_1
=(A-B,0,1-B^2,0)$,
and $X_2=\mathscr L\tilde X_2=(0,A-B,0,B^2-1)$,
so that
$\tilde X_1\cdot X_1=\tilde X_2\cdot X_2
=A-2B+B^3=\frac{(B^2-1)^2}{B}>0$ and $N_n^+=0$.
Since we found $n(\mathscr L)=1$, we conclude that $N_{C^+}=0$:
and there is no exponentially unstable solution to the linearized system
(which is indeed consistent with the fact that $\mathbb L$ has two purely imaginary eigenvalues).

   \begin{figure}[!htbp]
   \begin{subfigure}{0.32\textwidth}
  \includegraphics[height=4cm]{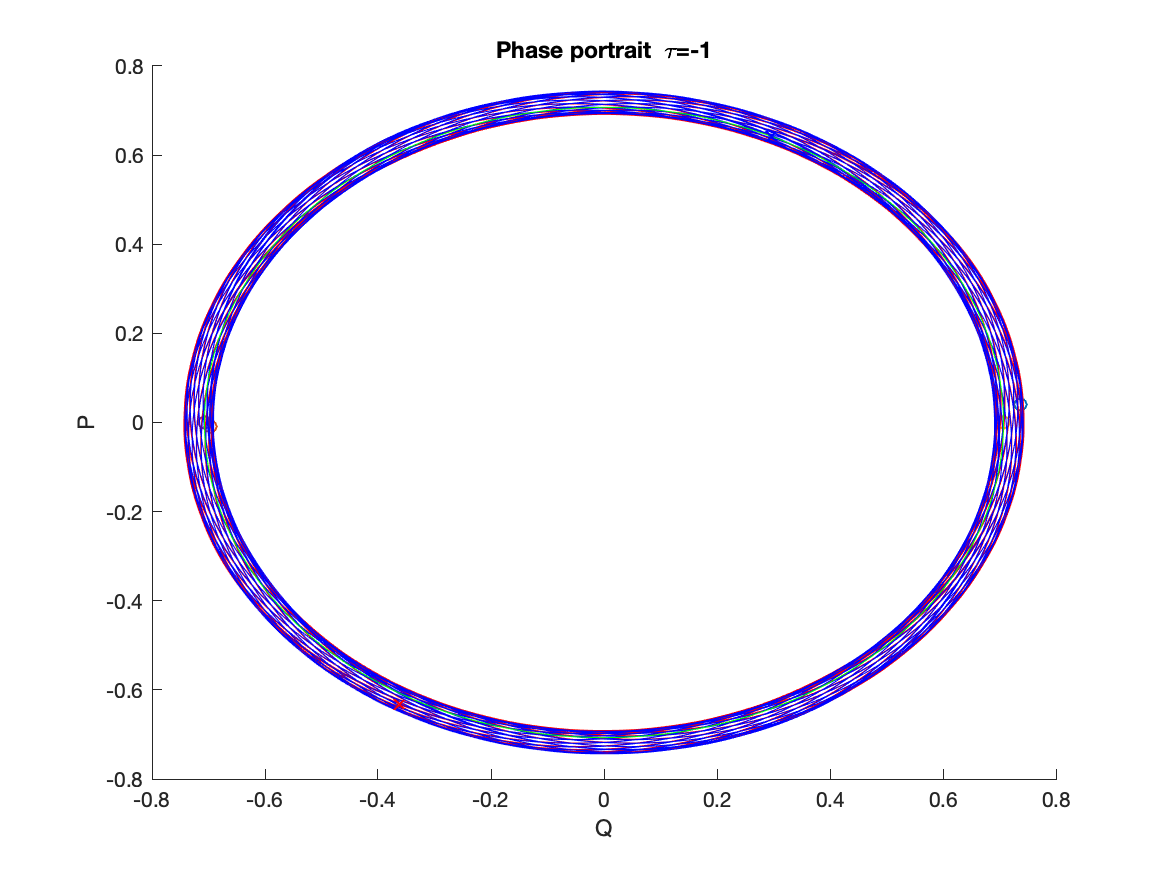}
  \caption{}
    \end{subfigure}
  \begin{subfigure}{0.32\textwidth}
   \includegraphics[height=4cm]{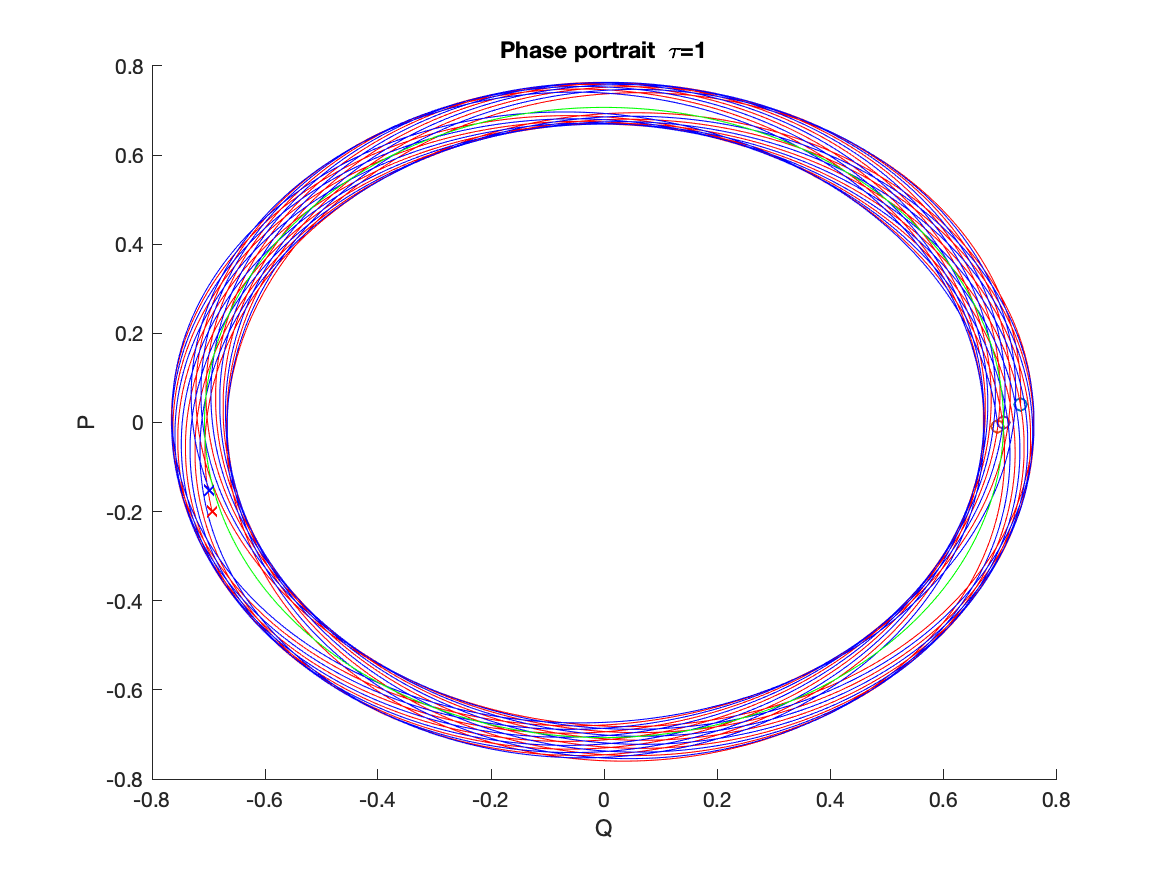}
   \caption{}
     \end{subfigure}
   \begin{subfigure}{0.32\textwidth}
     \includegraphics[height=4cm]{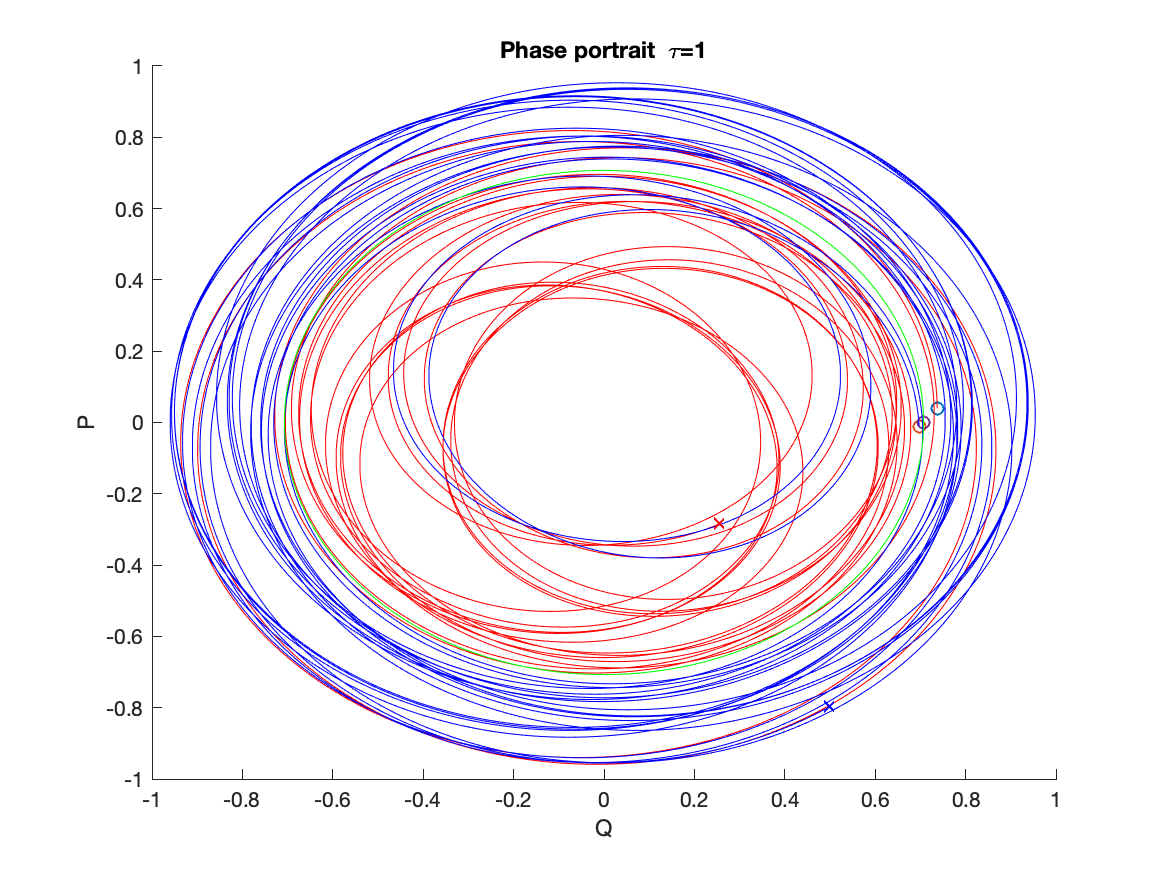}
     \caption{}
     \end{subfigure}
  \caption{Simulation of the non linear asymptotic model: phase portrait at $T=100$,
   with $\kappa=1.4$ and $\tau=-1$ (a),
   with $\kappa=1.4$ and $\tau=1$ (b),
   with $\kappa=2.4$ and $\tau=1$ (c).
    The circled points indicate the initial state, the cross indicate the final state
  } 
  \end{figure}

\subsection{Instability}\label{insta_As}

For $\tau=+1$, the status of the solution $X_*$  given by \eqref{solstatH2} changes as $\kappa$ overtakes the threshold 2: 
being a minimizer of the energy when $0<\kappa<2$, it becomes a local maximum when $\kappa>2$.   
  We have also seen that the Morse index of $\mathscr L$ switches from 1 to 2. 
  In this case, we can adapt the arguments presented in \cite{GSS,Maeda} to justify the instability of the reference solution when $\kappa
  >2$.
    To prove this statement, we need a series of preparation lemma, which exploit the algebraic properties of $\mathscr L$ and its spectral decomposition.

  \begin{lemma}\label{coerc1}
  We can find a constant $c>0$ such that for any $X\in \mathbb R^4$ verifying 
  $X\cdot X_*=X\cdot  X_{2-\kappa}=X\cdot X_0=0$, we have 
  $\mathscr LX\cdot X\geq c|X|^2$.
   \end{lemma}

\noindent
  {\bf Proof.}
  Since $(X_0,X_{-\kappa},X_{2-\kappa},X_2)$ forms an orthogonal basis of $\mathbb R^4$ and $X_*=X_{-\kappa}/\sqrt 2$, the vector we are considering is in fact proportional to $X_2$: from  $X=aX_2$, we deduce that 
  \[\mathscr LX\cdot X=a^2 \mathscr LX_2\cdot X_2=2a^2|X_2|^2=2|X|^2.\]
  \QED

  It is convenient to split $X_*=(X_{*0},X_{*1})$, with $X_{*0}=X_{*1}=\frac{1}{\sqrt 2}(1,0)$ and 
  to consider the rotation matrix
  in the plane
  \[R(\theta)=\begin{pmatrix}\cos(\theta)& -\sin(\theta)
  \\
  \sin(\theta)& \cos(\theta)\end{pmatrix}.\]
 We shall use the same notation for $V=(V_0,V_1)\in \mathbb R^2\times\mathbb R^2$, 
 $R(\theta)V=(R(\theta)V_0,R(\theta)V_1)$.
  
    \begin{lemma}\label{l1}
  Let $\epsilon>0$ and set 
  \[
  \mathscr U_\epsilon=\Big\{V=(V_0,V_1)\in\mathbb R^4,\ \inf_{\theta}
 |R(\theta) V- X_{*}|^2\leq\epsilon  \Big\}.\]
  For any $V\in \mathscr U_\epsilon$, there exists $\theta_*(V)\in [0,2\pi)$ such that 
  \[
   \inf_{\theta} |R(\theta) V- X_{*}|^2
   =
   |R(\theta_*(V)) V- X_{*}|^2.\]
    Moreover, the following relations hold 
    \[
   (i)\qquad  \theta_*(R(\theta')V)=\theta_*(V)-\theta',\qquad\qquad (ii)\qquad
    \nabla_{V_j}\theta_*(V)=\ds\frac{R'(\theta_*(V))^{^\intercal} X_{*j}}{R(\theta_*(V))^{^\intercal} X_{*j}\cdot V_j}.
    \]
     \end{lemma}

\noindent
  {\bf Proof.}
  The standard argument \cite{GSS,Maeda} relies on an application of the implicit function theorem.
  Here the construction can be made fully explicit.
  Indeed, given $V\in \mathbb R^4$, the $2\pi$-periodic function 
  $$\theta\longmapsto F(\theta)=  |R(\theta) V- X_{*}|^2=|R(\theta) V_0- X_{*0}|^2+|R(\theta) V_1- X_{*1}|^2$$
  admits a minimizor on $[0,2\pi]$, characterized by 
  \[
  F'(\theta)=2(R(\theta) V_0- X_{*0})\cdot R'(\theta)V_0+2(R(\theta) V_1- X_{*1})\cdot R'(\theta)V_1=0,\]
  where 
  \[
  R'(\theta)=\begin{pmatrix}-\sin(\theta)& -\cos(\theta)
  \\
  \cos(\theta)& -\sin(\theta)\end{pmatrix}.\]
  Since $$(R'(\theta))^{^\intercal} R(\theta)=
  \begin{pmatrix}0 & 1
  \\-1 & 0\end{pmatrix}
  $$
the relation becomes 
\begin{equation}\label{derF}
F'(\theta)=-2X_{*0}\cdot R'(\theta)V_0- 2X_{*1}\cdot R'(\theta)V_1=0.\end{equation}
Let $V_j=(Q_j,P_j)$. Using the specific expression of $X_{*j}$, we obtain
\[\sin(\theta)(Q_0+Q_1)
+\cos(\theta)(P_0+P_1)=0,
\]
which eventually determines the minimizer by
\[\mathrm{tan}(\theta_*(V))=-\ds\frac{P_0+P_1}{Q_0+Q_1}.\]
Differentiating \eqref{derF} with respect to $V_j$ and using $R''(\theta)=-R(\theta)$ yield
\[
(R'(\theta_*(V)))^{^\intercal} X_{*j} -(R(\theta_*(V)))^{^\intercal} X_{*j}\cdot V_j \nabla_{V_j}\theta_*(V)=0\]
and thus 
$$
\nabla_{V_j}\theta_*(V)=\ds\frac{(R'(\theta_*(V)))^{^\intercal} X_{*j}}
{(R(\theta_*(V)))^{^\intercal} X_{*j}\cdot V_j}
.$$
Finally, from $R(\theta+\theta')=R(\theta)R(\theta')$, we infer, for any $\theta,\theta'$,
\[
|R(\theta_*(V)-\theta')R(\theta')V_j-X_{*j}|=|R(\theta_*(V))V_j-X_{*j}|
\leq |R(\theta+\theta')V_j-X_{*j}|= |R(\theta)R(\theta')V_j-X_{*j}|\]
which means $\theta_*(V)-\theta'=\theta_*(R(\theta')V)$.
   \QED

\noindent
We observe that we can move from $X_*$ in a specific direction so that the energy decreases.
  
  \begin{lemma}\label{defci} Let $\kappa >2$
 and   set $V_s:s\in (-1/\sqrt 2,1/\sqrt2)\mapsto V_s=\sqrt{1-2s^2}  X_*+sX_{2-\kappa}$.
  Then, there exists $0<s_*<1/\sqrt 2$ such that for any $s\in [-s_*,s_*]$, we have
  $|V_s|=1 $ and $\mathscr E(V_s)<\mathscr E(X_*)$.
  \end{lemma}
  
  \noindent
  {\bf Proof.}
  We compute 
  \[|V_s|^2=(1-2s^2)|X_*|^2+s^2|X_{2-\kappa}|^2+s\sqrt{1-2s^2}  X_*\cdot X_{2-\kappa}
  =1-2s^2+2s^2+0=1.\]
  Next, owing to \eqref{0gradE}, we get the following  Taylor expansion
  \[
  \mathscr E(V_s)=\mathscr E\big(X_*+ sX_{2-\kappa}+(\sqrt{1-2s^2}-1)X_*\big)
  =\mathscr E(X_*)+\ds\frac{s^2}{2}\mathscr LX_{2-\kappa}\cdot X_{2-\kappa}+s^2\epsilon(s)
,\]
where $\lim_{s\to 0}\epsilon(s)=0$.
The conclusion follows from the fact that 
\[\mathscr LX_{2-\kappa}\cdot X_{2-\kappa}=(2-\kappa)|X_{2-\kappa}|^2<0.\]
  \QED

We are going to use the specific directions identified in Lemma~\ref{defci} 
to construct unstable solutions.
The instability will be characterized by 
working on a suitable functional framework which is adapted to the structure of the dynamical system. 
Let us now consider the functional
  \[A:V\in \mathscr U_\epsilon\longmapsto -X_{2}\cdot R(\theta_*(V))V
  =(V_1-V_0)\cdot \begin{pmatrix}\sin(\theta_*(V))\\  \cos(\theta_*(V))\end{pmatrix},\]
  bearing in mind $X_{2}=-\mathscr JX_{2-\kappa}$.
  By using Lemma~\ref{l1}-(ii), we get $R(\theta_*(R(\theta)V))R(\theta)V=R(\theta_*(V)-\theta)R(\theta)V
  =R(\theta_*(V))V$ so that 
$  A(R(\theta)V)=A(V)$.
Next, we get
\[\nabla_VA(V)=-R(\theta_*(V))^{^\intercal} X_{2}- (X_{2}\cdot R'(\theta_*(V))V) \nabla_V\theta_*(V).\]
For $V=X_*$, we have $\theta_*(X_*)=0$ and thus $X_{2}\cdot R'(\theta_*(X_*))X_*=\frac{X_{2}\cdot X_0}{\sqrt2}=0$ and
\begin{equation}\label{usef}
\nabla_VA(X_*)=-X_{2},\qquad
\mathscr J\nabla_VA(X_*)=-\mathscr JX_{2}=-X_{2-\kappa}.\end{equation}
Eventually, since  $\begin{pmatrix}
  R(\theta) &0\\ 0 & R(\theta)
\end{pmatrix}\mathscr J=- \begin{pmatrix}
  R'(\theta) &0\\ 0 & R'(\theta)
\end{pmatrix}$ and $\mathscr J^2=-\mathbb I$,  we observe that 
\[\begin{array}{lll}
\nabla_VA(V)\cdot \mathscr J V&=&-R(\theta_*(V))^{^\intercal} X_{2}\cdot \mathscr J V
- (X_{2}\cdot R'(\theta_*(V))V )(\nabla_V\theta_*(V)\cdot \mathscr J V)
\\[.3cm]
&=&
-R(\theta_*(V))^{^\intercal} X_{2}\cdot \mathscr J V
+ (X_{2}\cdot R'(\theta_*(V))V )\ds\frac{-R'(\theta_*(V))^{^\intercal} X_*\cdot \mathscr J V}{
R(\theta_*(V))^{^\intercal} X_*\cdot V}
\\[.3cm]
&=&
 X_{2}\cdot R'(\theta_*(V))V
 +
 (X_{2}\cdot R'(\theta_*(V))V)
 \ds\frac{-
 X_*\cdot R(\theta_*(V)) V}
{
 X_*\cdot R(\theta_*(V)) V}
=0
.
\end{array}\]
This estimate of Lemma~\ref{defci} can be strengthened as follows.

  \begin{lemma}\label{defci2}  Let $\kappa>2$, set
  \[\mathscr P(V)=\nabla_V A(V)\cdot \mathscr J\nabla_V\mathscr E(V)\]
  and let $V_s$ be defined as in Lemma \ref{defci}.
  Then, there exists $0<s_*<1/\sqrt 2$ such that for any $s\in [-s_*,s_*]$, we have
   $$0<\mathscr E(X_*)-\mathscr E(V_s)<-s\mathscr P(V_s).$$
  \end{lemma}

 \noindent
  {\bf Proof.}
  The proof is again based on Taylor expansions.
  In what follows we denote by $\varrho(s)$ a reminder, the expression of which might change from a line to another, but such that $\lim_{s\to 0}\varrho (s)=0$.
  Since $V_s$ looks like $X_*+sX_{2-\kappa}$,  we get, by virtue of \eqref{0gradE} and \eqref{usef}, 
  \[\begin{array}{lll}
  \mathscr P(V_s)&=&
  s(\nabla_V A(X_*)+ sD^2_VA(X_*)X_{2-\kappa})\cdot  \mathscr J  D^2_V\mathscr E(X_*)X_{2-\kappa})+s\varrho(s)
  \\
  &=&
  -sX_{2}\cdot  \mathscr J \mathscr LX_{2-\kappa}+s\varrho(s)
  =s \mathscr LX_{2-\kappa}\cdot X_{2-\kappa}+s\varrho(s)
  =s(2-\kappa)|X_{2-\kappa}|^2+s\varrho(s)
  .
  \end{array}\]
  Accordingly, we obtain
  \[
 \mathscr E(X_*)- \mathscr E(V_s)+s\mathscr P(s)
 = \frac{s^2}{2}\big((2-\kappa)+\varrho(s)\big)
  \]
  which thus remains negative for $s$ small enough.
  \QED

Note that $\mathscr P(V)=\nabla_V A(V)\cdot \mathscr J\nabla_V\mathscr H(V)$ since $\nabla_V F(V)=V$ 
and $\nabla_V A(V)\cdot \mathscr J V=0$ for all $V\in \R^4$.
The motivation for introducing the functional $A$ and $\mathscr P$ comes from  the fact that, for $X$ solution of \eqref{Hartreebis},
we have 
\begin{equation}\label{diffAP}\ds\frac{\ud}{\ud t}A(X(t))
=\nabla_U A(X(t))\cdot \ds\frac{\ud}{\ud t}X(t)
=\nabla_U A(X(t))\cdot \mathscr J\nabla\mathscr H (X(t))
=\mathscr P(X(t)).\end{equation}

\begin{lemma}  Let $\kappa>2$ and
$\epsilon>0$ be sufficiently small.
Let $V\in \mathscr U_\epsilon$ be such that 
$|V|=|X_*|$ and $\mathscr E(X_*)-\mathscr E(V)>0$.
Then, we actually have
\[
\mathscr E(X_*)-\mathscr E(V)<-\Lambda (V)\mathscr P(V)\]
where
\begin{equation}\label{defLamb}
\Lambda (V)=\ds\frac{R(\theta_*(V))V\cdot X_{2-\kappa}}{|X_{2-\kappa}|^2}.\end{equation}
\end{lemma}

 \noindent
  {\bf Proof.}
For $V\in \mathbb R^4$, set
\begin{equation}\label{defM}
M(V)=R(\theta_*(V))V-X_*-\Lambda(V)X_{2-\kappa},\end{equation}
 so that 
$M(V)\cdot X_{2-\kappa}=0$.
Moreover, we have
$$M(V)\cdot X_0=R(\theta_*(V))V\cdot X_0
=\sqrt 2
(-\mathscr J R'(\theta_*(V)V)\cdot(-\mathscr J X_*))
=\sqrt 2R'(\theta_*(V))V\cdot X_*=0,$$
by definition of $\theta_*(V)$, see \eqref{derF}.
As a consequence, $M(V)$ lies in the orthogonal space of $\mathrm{Span}(X_0,X_{2-\kappa})$ and it can be written
\begin{equation}\label{defatM}
M(V)=a(V)X_* + \tilde M(V), \textrm{ where $\tilde M(V)\in \mathrm{Span}(X_2)$}.\end{equation}
Lemma~\ref{coerc1} tells us that $\mathscr L \tilde M(V)\cdot \tilde M(V)\geq c|\tilde M(V)|^2.$

%
  
 We start by proving 
 \begin{equation}\label{eqP}
 \mathscr P(V)=\mathscr P(R(\theta_*(V))V)
 .\end{equation}
Derivating  $\mathscr H(R(\theta)V)=\mathscr H(V)$ and using Lemma~\ref{l1}-(i), we get
\[
R(\theta)^{^\intercal}\nabla \mathscr H(R(\theta)V)=\nabla\mathscr H(V),\qquad
R(\theta)^{^\intercal}\nabla \theta_*(R(\theta) V)=\nabla \theta_*(V),
\]
while
\[
R(\theta)^{^\intercal}=R(-\theta)=R(\theta)^{-1},\qquad 
\begin{pmatrix}R'(\theta) & 0 \\ 0 & R'(\theta) \end{pmatrix}
=-\mathscr J \begin{pmatrix}R(\theta) & 0 \\ 0 & R(\theta) \end{pmatrix}.\]
Therefore, we obtain 
 \[
 \mathscr P(R(\theta)V)=\nabla_U A(R(\theta)V)\cdot \mathscr J\nabla \mathscr H (R(\theta)V) 
 =
 \nabla_U A(R(\theta)V)\cdot \mathscr JR(\theta)\nabla \mathscr H (V). 
 \]
 where
  $$\begin{array}{lll}
  \nabla A(R(\theta)V)&=&
  -R(\theta_*(R(\theta)V))^{^\intercal}X_2-(X_2\cdot R' (\theta_*(R(\theta)V))R(\theta)V)\ \nabla\theta_*(R(\theta)V)
  \\
  &=&
  -R(\theta_*(V)-\theta)^{^\intercal}X_2-
  (X_2 \cdot R'(\theta_*(V)-\theta)R(\theta)V)\ 
  R(\theta)\nabla\theta_*(V)
  \\
  &=&
  -
  R(\theta) R(\theta_*(V))^{^\intercal}X_2
  +(X_2\cdot \mathscr J R(\theta_*(V))R(-\theta)R(\theta)V)\ 
  R(\theta)\nabla\theta_*(V)
  \\
  &=&
  
  R(\theta) \big[-
  R(\theta_*(V))^{^\intercal}X_2
  -(X_2 \cdot R'(\theta_*(V))V)\ 
\nabla\theta_*(V)
\big]=R(\theta)\nabla A(V).
  \end{array}$$
Hence,  \eqref{eqP} holds.

Let $V\in \mathscr U_\epsilon$.
The definition of $\Lambda(V)$ in \eqref{defLamb}, $M(V)$ in \eqref{defM} and  $a(V), \tilde M(V)$ in \eqref{defatM}
leads to
the estimates
\[\begin{array}{l}
  |\Lambda(V)|^2
  =\ds\frac{|R(\theta_*(V))V\cdot X_{2-\kappa}|^2}{|X_{2-\kappa}|^4}=
  \ds\frac{|(R(\theta_*(V))V-X_*)\cdot X_{2-\kappa}|^2}{|X_{2-\kappa}|^4}
  \leq  \ds\frac{|R(\theta_*(V))V-X_*|^2}{|X_{2-\kappa}|^2}
  \leq \ds\frac{\epsilon^2}{4},
  \\
  [.3cm]
  |M(V)|
  \leq |R(\theta_*(V))V-X_*|+|\Lambda(V)X_{2-\kappa}|
  \leq 2\epsilon ,
  \\[.3cm]
  |a(V)|\leq |M(V)|\leq 2\epsilon,
  \\[.3cm]
  |\tilde M(V)|\leq |M(V)|+|a(V)|\leq 4\epsilon.
  \end{array}
\]
 Now, we perform a Taylor expansion on 
 \[
 \mathscr P(V)=\mathscr P(R(\theta_*(V))V)=\mathscr P(X_*+ \Lambda(V)X_{2-\kappa} +a(V)X_*+\tilde M(V) ),\]
 based on the fact that $\varrho(V)= \Lambda(V)X_{2-\kappa} +a(V)X_*+\tilde M(V)$ is of the order of $\epsilon$.
 Hence, we get 
 \[\begin{array}{lll}
  \mathscr P(V)&=&
  \nabla A(X_*+\varrho(V))\cdot \mathscr J \nabla \mathscr H(X_*+\varrho(V))
 \\
 &=& (\nabla A(X_*)+D^2A(X_*)\varrho(V)) \cdot \mathscr J D^2 \mathscr H(X_*)\varrho(V)
 +\mathscr O(\epsilon^2)
  \\
 &=&
  \nabla A(X_*) \cdot \mathscr J \mathscr L\varrho(V)
 +\mathscr O(\epsilon^2)
 =-  \mathscr L\mathscr J \nabla A(X_*)\cdot\varrho(V)+\mathscr O(\epsilon^2)
  \\
 &=&
  \mathscr L X_{2-\kappa}\cdot (\Lambda(V)X_{2-\kappa} +a(V)X_*+\tilde M(V))
  +\mathscr O(\epsilon^2)\\
  &=&
  (2-\kappa) \Lambda(V)|X_{2-\kappa}|^2+ \mathscr O(\epsilon^2).
 \end{array}\]
 Accordingly, we have \begin{equation}\label{est1}
 -\Lambda(V)\mathscr P(V)=-(2-\kappa) \Lambda(V)^2|X_{2-\kappa}|^2+\mathscr O(\epsilon^3).\end{equation}
 Similarly, we go back to the difference of energies
 \[\begin{array}{l}
 0<\mathscr E(X_*)-\mathscr E(V)=\mathscr E(X_*)-\mathscr E(X_*+\varrho(V))
 =-\ds\frac12\mathscr L\varrho(V)\cdot \varrho(V) +\mathscr O(\epsilon ^3)
 \\[.3cm]
 \hspace*{.5cm}
= -\ds\frac12\mathscr L(\Lambda(V)X_{2-\kappa} +a(V)X_*+\tilde M(V))\cdot (\Lambda(V)X_{2-\kappa} +a(V)X_*+\tilde M(V))+\mathscr O(\epsilon ^3)
 \\[.3cm]
 \hspace*{.5cm}
=- \ds\frac{1}2\big((2-\kappa)\Lambda(V)X_{2-\kappa} -\kappa a(V)X_*+\mathscr L \tilde M(V)\big)\cdot (\Lambda(V)X_{2-\kappa} +a(V)X_*+\tilde M(V))+\mathscr O(\epsilon ^3)
 \\[.3cm]
 \hspace*{.5cm}
= -\ds\frac{2-\kappa}2\Lambda(V)^2|X_{2-\kappa}|^2
+\ds\frac\kappa2|a(V)|^2
-\ds\frac12\mathscr L\tilde M(V)\cdot \tilde M(V)+\mathscr O(\epsilon ^3).
 \end{array}\]
 We now need to refine the estimate on 
 $a(V)=M(V)\cdot X_*=(R(\theta_*(V))V-X_*)\cdot X_*$.
 To this end, we use the elementary relation
 \[\begin{array}{lll}
 0=|V|^2-|X_*|^2&=& |R(\theta_*(V))V|^2-|X_*|^2=
 |(R(\theta_*(V))V-X_*)+X_*|^2-|X_*|^2
 \\
 &=&|R(\theta_*(V))V-X_*|^2 +2(R(\theta_*(V))V-X_*)\cdot X_*
\\
& =&|R(\theta_*(V))V-X_*|^2+2a(V),\end{array}\]
which yields $|a(V)|\leq \frac{\epsilon^2}{2}$.
We are thus led to 
\[\begin{array}{lll}
0<\mathscr E(X_*)-\mathscr E(V)&=&
 -\ds\frac12(2-\kappa)\Lambda(V)^2|X_{2-\kappa}|^2
-\ds\frac12\mathscr L\tilde M(V)\cdot \tilde M(V)+\mathscr O(\epsilon ^3)
\\[.4cm]
&\leq& -\ds\frac{(2-\kappa)}{2}\Lambda(V)^2|X_{2-\kappa}|^2+\mathscr O(\epsilon ^3)
\end{array}\]
since $\mathscr L\tilde M(V)\cdot \tilde M(V)\geq 0$.
In particular, this implies  that $\Lambda(V)$ does not vanish.
  We conclude by going back to \eqref{est1}. 
 \QED
 
We argue by contradiction for establishing Theorem~\ref{ThUnstHa}
We assume that $X_*$  given by \eqref{solstatH2} is an orbitally stable of \eqref{Hartreebis}, meaning that for any $\epsilon>0$, we can find $\delta$ such that 
$X_{\mathrm {init}}\in \mathscr U_\epsilon$ implies 
$X(t)\in \mathscr U_\epsilon$ for any $t\geq0$.
Then, as an initial data we pick $X_{\mathrm {init}}
=V_s$ as defined in Lemma~\ref{defci} with $s<0$ small enough (see Lemma~\ref{defci2}) so that 
$|V_s|=|X_*|$, $\mathscr E(X_*)-\mathscr E(V_s)=\epsilon_*>0$ and $\mathscr P(V_s)>0$.
Let $t\mapsto X(t)$ be the associated solution.
By using the conservation properties of the equation, we obtain
$$0<\epsilon_*=\mathscr E(X_*)-\mathscr E(X(t))<-\Lambda (X(t))\mathscr P(X(t)).$$
Since $\mathscr P(V_s)>0$ and $|\Lambda(X(t))|\leq \frac\epsilon 2$, 
we get $\mathscr P(X(t))\geq C\epsilon_*$ 
for a certain $C>0$. 
We now use \eqref{diffAP}.
Consequently,
there holds 
\[
C\epsilon_*t  \leq 
\ds\int_0^t \mathscr P(r)\ud r=
\ds\int_0^t  \ds\frac{\ud}{\ud t}A(X(r))\ud r
=A(X(t))-A(V_s).
\]
This contradicts the stability assumption $\{X(t),\ t\geq 0\}\subset \mathscr U_\epsilon$ 
which implies that $A(X(t))$ remains bounded. Indeed, $|A(X(t))|\leq |X_2|\ |R(\theta_*(X(t)))X(t)|
\leq |X_2|(|R(\theta_*(X(t)))X(t)-X_*| + |X_*|)\leq 
 |X_2|( \epsilon + |X_*|)$.
%
  
  \section{Stability analysis for the coupled system \eqref{S1}-\eqref{W1}}
  \label{AnSW}
  
  \subsection{Linearized equations}

    \subsubsection{Linearization about the solution \eqref{SpecSol}}}
    
  We search for solutions of \eqref{S1}-\eqref{W1} on the form of a perturbation of \eqref{SpecSol}: 
  $$u_j=e^{i\omega t}(U_{*j}+v_j),\qquad \psi_j=\Psi_{*j}+\phi_j, \qquad \Psi_{*j}=-|U_{*j}|^2(-\Delta)^{-1}\sigma.$$
  Using $|u+h|^2=|u|^2+2\mathrm{Re}(\overline u h)+|h|^2$ and the dispersion relation \eqref{defom}, 
  we arrive at the following  linearized system
  \[\begin{array}{l}
  i\ds\frac{\ud}{\ud t} v_0=\tau v_0 -v_1+\ds\frac{1}{\sqrt 2} \ds\int_{\mathbb R^n}\sigma \phi_0\ud z,
  \\[.4cm]
  i\ds\frac{\ud}{\ud t} v_1=\tau v_1 -v_0+\ds\frac{\tau}{\sqrt 2} \ds\int_{\mathbb R^n}\sigma \phi_1\ud z,
   \\[.4cm]
   \Big(\ds\frac{1}{c^2}\partial^2_{tt} -\Delta\Big)\phi_0=-\sqrt 2 \sigma \mathrm{Re}(v_0),
   \\[.4cm]
     \Big(\ds\frac{1}{c^2}\partial^2_{tt} -\Delta\Big)\phi_1=-\tau\sqrt 2 \sigma \mathrm{Re}(v_1).
  \end{array}
  \] 
 It is convenient to introduce new unknowns. On the one hand, we expand the complex unknown and  consider its real and imaginary parts
$u_j=q_j+ip_j;$
on the other hand, for the wave equation, we set
$$\varphi_j =(-\Delta)^{1/2}\phi_j,\qquad \varpi_j= \ds\frac{\partial_t \phi_j}{c}.$$
We use a block decomposition of the unknown: 
\begin{equation}
  \label{defXSW}
  X=\begin{pmatrix} S\\ W\end{pmatrix},\quad
  W=\begin{pmatrix} W_0\\ W_1\end{pmatrix},\quad
 S=\begin{pmatrix} S_0\\ S_1\end{pmatrix},\quad
  W=\begin{pmatrix} W_0\\ W_1\end{pmatrix},\quad
  S_j=\begin{pmatrix} q_j\\ p_j\end{pmatrix},\quad
  W_j=\begin{pmatrix} \varphi_j\\ \varpi_j\end{pmatrix}.  
\end{equation}
 Therefore, $X$ has  8 components
$(q_0,p_0,q_1,p_1,\phi_0,\varpi_0,\phi_1,\varpi_1)$ and is valued in $\mathbb R^4\times (L^2(\mathbb R^n))^4$.
With these notations, the 
problem casts as
\[
\partial_t X=\mathbb L X,\]
where 
\[
\mathbb L X
=
\begin{pmatrix}
\tau p_0- p_1
\\
-\tau q_0+ q_1- \ds\frac{1}{\sqrt 2}\ds\int_{\mathbb R^n} \sigma (-\Delta)^{-1/2}\varphi_0\ud z
\\
-p_0+\tau p_1
\\
q_0-\tau q_1- \ds\frac{\tau}{\sqrt 2}\ds\int_{\mathbb R^n} \sigma (-\Delta)^{-1/2}\varphi_1\ud z
\\
 c(-\Delta)^{1/2}\varpi_0
\\
- c(-\Delta)^{1/2}\varphi_0- c\sqrt 2\sigma q_0
\\
 c(-\Delta)^{1/2}\varpi_1\\
- c(-\Delta)^{1/2}\varphi_1- c\sqrt 2\tau \sigma q_1
\end{pmatrix}.
\]
The following statements bring out the basic spectral properties of $\mathbb L$ and 
makes the symplectic structure appear. In terms of stability analysis, it implies that 
the linearized system is stable provided $\sigma(\mathbb L)\subset i\mathbb R$.
However, the identification of the eigenvalues of $\mathbb L$ is now not so direct than for the asymptotic problem.
The symplectic structure will be crucial to decide whether or 
not the equation is spectrally stable.

\begin{proposition}\label{Spec0LSW}
Let us denote by $\check X$ the vector constructed from $X$ by changing the components $p_j$ and $\varpi_j$ into 
$-p_j$ and $-\varpi_j$.
Let $(\lambda,X)$ be an eigenpair of $\mathbb L$. Then, $(-\lambda, \check X)$, $(\overline \lambda, \overline X) $
and $(-\overline \lambda, \overline {\check X})$ are equally eigenpairs of $\mathbb L$.

Moreover, we can write $\mathbb L=\mathscr J\mathscr L$ with $\mathscr J$ a skew-symmetric operator and $\mathscr L$ a self-adjoint operator.
\end{proposition}

\noindent {\bf Proof.}
The first part of the claim follows by direct inspection and using the fact that $\mathbb L$ has real coefficients.
Next, we introduce the following block-wise operator $\mathscr J$ and its formal inverse $\tilde {\mathscr J}$ 
\begin{equation}\label{defJ}
\mathscr J=
\begin{pmatrix} 
\mathscr J_S & 0 & 0 & 0
\\
  0 &  \mathscr J_S & 0 & 0
  \\
0 & 0 &  \mathscr J_W  & 0
\\
0 & 0 & 0&  \mathscr J_W \end{pmatrix}
,\qquad
\tilde{\mathscr J}=
\begin{pmatrix} 
\tilde{\mathscr J_S} & 0 & 0 & 0
\\
  0 &  \tilde{\mathscr J}_S & 0 & 0
  \\
0 & 0 &  \tilde {\mathscr J}_W  & 0
\\
0 & 0 & 0&  \tilde {\mathscr J}_W\end{pmatrix}
\end{equation}
where
\[\begin{array}{ll}
\mathscr J_S=
\begin{pmatrix}
0 & 1
\\-1&0 
\end{pmatrix},
\quad &
\tilde{\mathscr J}_S=
\begin{pmatrix}
0 & -1
\\1&0 
\end{pmatrix},
\\ 
\mathscr J_W=  2c\begin{pmatrix}
0 & (-\Delta)^{1/2} 
\\
-(-\Delta)^{1/2}  & 0 
\end{pmatrix},\quad&
\tilde{\mathscr J}_W= \ds\frac1{2c}
\begin{pmatrix}
0 & -(-\Delta)^{-1/2} 
\\
(-\Delta)^{-1/2}  & 0 
\end{pmatrix}.
\end{array}\]
We obtain 
 \begin{equation}\label{defL}
\mathscr LX=\tilde {\mathscr J}\mathbb LX
=
\begin{pmatrix}
-q_1+\tau q_0 +\ds\frac{1}{\sqrt 2}\ds\int_{\mathbb R^n} \sigma (-\Delta)^{-1/2}\varphi_0\ud z
\\
\tau p_0-p_1
\\
-q_0+\tau q_1 +\ds\frac{\tau}{\sqrt 2} \ds\int_{\mathbb R^n} \sigma (-\Delta)^{-1/2}\varphi_1\ud z
\\
-p_0+\tau  p_1
\\
\ds\frac1{ 2} \varphi_0 +\ds\frac{1}{\sqrt 2}(-\Delta)^{-1/2}\sigma q_0
\\
\ds\frac1{2} \varpi_0
\\
\ds\frac1{ 2} \varphi_1 +\ds\frac{\tau}{\sqrt 2}(-\Delta)^{-1/2}\sigma q_1
\\
\ds\frac1{ 2} \varpi_1
\end{pmatrix}.\end{equation}
We readily check that $(\mathscr L X|X')=(X|\mathscr L X')$ holds for the inner product
$(X|X')=\sum_{j=0}^1 q_jq'_j+p_jp'_j+\int_{\mathbb R^n} (\varphi_j\varphi'_j+\varpi_j\varpi'_j)\ud z$.

The change of unknowns is boiled down to ensure that $\mathscr L$ is self-adjoint and, moreover, that the product $(\mathbb L X|X)$ 
 does not involve derivatives of $\varphi_j$ or $\varpi_j$,
a property that will be useful later on (see Section~\ref{Sec:insta}).
\QED

A natural attempt to localize the eigenvalues of $\mathbb L$ would rely on a
asymptotic argument from the simplified problem \eqref{linHa}. However this program faces severe difficulties.
We have seen that the eigenvalues of the asymptotic problem 
lie in $i\mathbb R$; we would like to decide whether  they get stuck on the imaginary axis or  they split into 
branches with non zero real parts as the wave speed $c$ becomes finite.
The coupling with the wave equation induces obstructions to develop the asymptotic arguments (as for instance in \cite{HarG}) that can be described as follows.
 Let us introduce the function
\[\epsilon\geq 0\longmapsto  \kappa_{\epsilon}=\ds\int_{\mathbb R^n} \ds\frac{|\widehat \sigma(\xi)|^2}{\epsilon+|\xi|^2}\ds\frac{\ud \xi}{(2\pi)^n}.
\]
We have $0<\kappa_\epsilon\leq \kappa$, and
by applying the Lebesgue theorem, we can check the continuity of $\epsilon\mapsto \kappa_\epsilon$. However 
it fails to be derivable in general since 
$\frac{\ud}{\ud\epsilon}
\frac{|\widehat\sigma(\xi)|^2}{\epsilon+|\xi|^2}
=
-\frac{|\widehat\sigma(\xi)|^2}{(\epsilon+|\xi|^2)^2}
$
is not integrable for $\epsilon=0$ without introducing further restriction on the dimension $n$
(as $\xi\to 0$ it behaves like $\frac{(\int \sigma(x)\ud x)^2}{|\xi|^4}$).
This explains that developments of the eigenvalues as power series of $1/c$ 
are 
misleading. 
Let us go back to  the function 
\[\lambda=a+ib\in \mathbb C\longmapsto  \kappa_{\lambda}=\ds\int_{\mathbb R^n} \ds\frac{|\widehat \sigma(\xi)|^2}{\lambda^2+|\xi|^2}\ds\frac{\ud \xi}{(2\pi)^n}
=
\ds\int_{\mathbb R^n} \ds\frac{|\widehat \sigma(\xi)|^2}{2iab+ a^2-b^2+|\xi|^2}\ds\frac{\ud \xi}{(2\pi)^n}.
\]
that we now define on the complex plane.
The definition makes sense, but on the imaginary axis $a=0$. 
Let us set $A=a^2-b^2$ and $B=2ab$. Since $\sigma$ is radially symmetric, 
 we are led to consider the function
 \[
 P(A,B)=\ds\int_{0}^\infty \ds\frac{\Sigma(r)}{iB+ A+r^2}\ud r,
\]
with $\Sigma(r)=|\widehat \sigma(r)|^2r^{n-1}$.
It is well defined for $B=0$, and $A\geq 0$, and for any $B\neq 0$, $A\in \mathbb R$; 
the difficulty is to deal with $B=0$ and $A=-\mu<0$.
The lack of continuity near the imaginary axis is illustrated by the following Plemelj like formula: for $A<0$ fixed, 
the limits $B\to 0^\pm$ do not coincide. It reflects the jump discontinuity in the resolvent function of $-\Delta$ at the spectrum.

\begin{lemma}\label{Plem}
Let $\mu>0$.
Then, we have
\[\ds\lim_{B\to 0^\pm}
P(-\mu,B)=\mathrm{P.V.}\ds\int_{0}^{+\infty} \ds\frac{\Sigma(r)}{(r-\sqrt\mu)(r+\sqrt\mu)}\ud r
\mp i \ds\frac{\pi\Sigma(\sqrt\mu)}{2\sqrt \mu}.\]
\end{lemma}

For the sake of completeness, the detailed proof is provided in Appendix~\ref{App:Plem}.
The statement  can be expressed by means of  the limited absorption principle for the wave equation.
This difficulty we are facing can indeed be explained by coming  back to the  
 the wave equation, which has an
essential spectrum lying all along  the imaginary axis. 
As we shall 
detail below, 
we need to discuss Helmholtz type equation $(\lambda-\Delta)u=f$.
The equation perfectly makes sense provided $\lambda \in \mathbb C\setminus(-\infty,0]$.
For negative $\lambda$, in dimension $d=3$, this leads to consider 
$u_\pm(x)=\int\frac{e^{\pm i\sqrt{-\lambda}|x-y|}}{|x-y|}f(y)\ud y$
which both define solutions of the  Helmholtz  equation, with a different behavior at infinity.
These solutions can be obtained as the limits of $(\lambda\pm i\epsilon-\Delta)^{-1}f$ as $\epsilon\to 0$.
Hence the resolvent operator is not well defined, and the functional integrals that one would like to apply 
as in \cite{HarG}
are misleading.
\\ 

\begin{figure}[!h]
\begin{center}
\includegraphics[height=7cm]{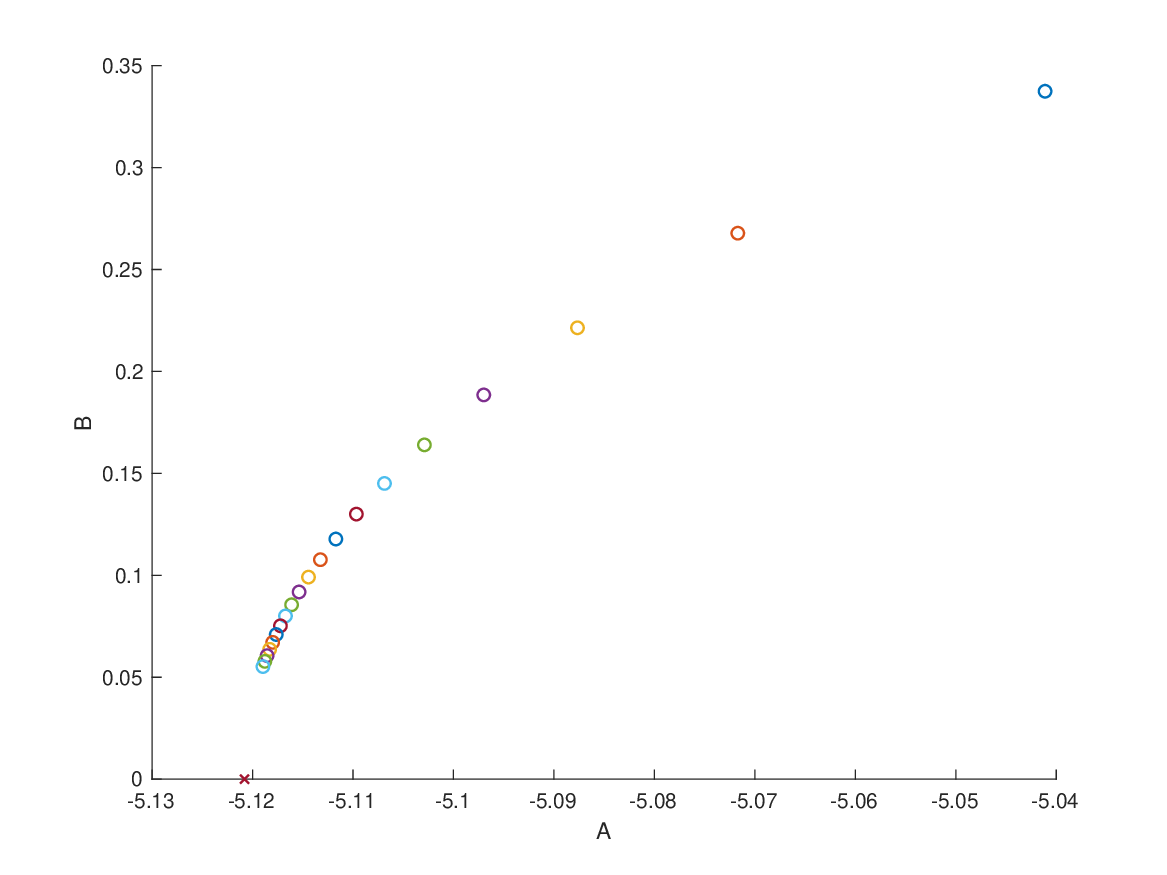}
\end{center}
\caption{Numerical identification of the zeros of $F$ for several values of the wave speed $c$ ($\kappa=0.5604$ and $\tau=-1$). The cross on the horizontal axis indicates the coordinates corresponding to the eigenvalue of the asymptotic problem.
}\label{EvolZero}
\end{figure}

Let us further illustrate how the difficulty shows up.
Searching for eigenvalues of $\mathbb L$, we are led to the following non linear equation for $\lambda\in \mathbb C$
(see the detailed computations in \eqref{resolvent} below)
\begin{equation}\label{nllam}
\lambda^2+4-2\tau \kappa_{\lambda^2/c^2}=0.\qquad
\end{equation}
We wonder whether or not there exists a solution $\lambda=a+ib$ 
with positive real and imaginary parts. Hence we set $A=a^2-b^2$ and $B=2ab$. The latter is supposed to be $\neq 0$ 
and we are thus led to investigate the zeros of the function
\[
F:(A,B)\in \mathbb R^2\longmapsto
\begin{pmatrix}
 A +4 -2\tau c^2\ds\int_{\mathbb R^n} \ds\frac{( A+c^2\xi^2)|\widehat \sigma(\xi)|^2}{(A+c^2\xi^2) ^2+ B^2}
\ds\frac{\ud \xi}{(2\pi)^n}
 \\[.4cm]
 1+2\tau \ds\int_{\mathbb R^n} \ds\frac{|\widehat \sigma(\xi)|^2}{(A+c^2\xi^2)^2 + B^2}\ds\frac{\ud \xi}{(2\pi)^n}

 \end{pmatrix}
.\]
We do not find explicit solutions for the relation $F(A,B)=0$, but the problem can be investigated numerically, 
based on the Newton algorithm.
Note however that the Jacobian matrix $\nabla F(A,B)$ becomes singular as $B$ tends to 0, making the problem stiffer as the solution $\lambda$ is getting close to the imaginary axis.
Fig.~\ref{EvolZero} displays the zeros of $F$ in the $(A,B)$-plane, for several values of the wave speed $c$.
As $c$ becomes large, we see that the zeros tends to the eigenvalue of the asymptotic problem, which lies on the horizontal axis. It conforts the intuition 
that the eigenvalues of $\mathbb L$ for the coupled problem  
do have a real part, thus 
leading to instability, and they should converge as $c\to \infty$ to the purely  imaginary eigenvalues of the asymptotic problem.
\\

%
 
 For these reasons, we are going to deduce spectral properties on $\mathbb L$ from the spectral analysis of $\mathscr L$, as proposed in 
 \cite{ChouPel}.
Indeed, the spectral analysis of the operator $\mathscr L$ is easier; at least we know that the spectrum embeds into $\mathbb R$ due to the self-adjointness character of $\mathscr L$.
The spectral properties of the operator $\mathscr L$ are summarized in the following statement.
Note that, due to the coupling with the wave equation on the whole $\mathbb R^n$, 
there is a non empty essential spectrum. From now on, we denote by $\mathbf{0}=(0,0,0,0)$.

\begin{theo}\label{ThSpecSW}
Let $\mathscr L$ be the operator defined by \eqref{defL}.
Then, the following assertions hold:
\begin{enumerate}
\item $\mathrm{Ker}(\mathscr L)=\mathrm{Span}(X_0)$, with $X_0=(S_0,\mathbf{0})$,  $S_0=(0,1,0,\tau)$;
\item  $\sigma_{\mathrm{ess}}(\mathscr L)=\{1/2\}$; 
\item If $\tau=+1$ and $0<\kappa<2$, $\mathscr L$ has one negative eigenvalue, 
associated to a one-dimensional eigenspace;  if $\tau=+1$ and $\kappa>2$, $\mathscr L$ has two negative eigenvalues, 
associated to one-dimensional eigenspaces; 
if $\tau =-1$ 
$\mathscr L$ has  three negative eigenvalues associated to one-dimensional eigenspaces;
\item Given $Y_0$  a solution of $\mathscr LY_0=-\mathscr J X_0$,  we have $(-\mathscr JX_0|Y_0)<0$.
\end{enumerate}
\end{theo}

\noindent
{\bf Proof.}
The operator $\mathscr L$ being self-adjoint, its spectrum lies in $\mathbb R$.
Let us study the solutions of $\mathscr L X=\lambda X$.

In particular, we have $\lambda \varpi_j=
\varpi_j/2$.
Hence, when $\lambda=1/ 2
$,  any $X=(\mathbf{0},W)$, with $W=(0,\pi,0,0)$ or $(0,0,0,\pi)$, $\pi\in L^2(\mathbb R^n)$,
 lies in $\mathrm{Ker}(\mathscr L-1/ 2)$. 
 Next, we also have
 \begin{align*}
  \Big(\lambda-\frac12\Big)\varphi_0=\frac{1}{\sqrt{2}}(-\Delta)^{1/2}\sigma q_0,\\
  \Big(\lambda-\frac12\Big)\varphi_1=\frac{\tau}{\sqrt{2}}(-\Delta)^{1/2}\sigma q_1.
 \end{align*}
 Moreover, we can write
 $$\ds\int \sigma (-\Delta)^{-1/2}\varphi_j\ud z=
 \frac{1}{(2\pi)^n}\int \frac{\widehat \sigma \overline{\widehat {\varphi_j}}}{|\xi|}\ud \xi
 =\ds\int (-\Delta)^{-1/2}\sigma \varphi_j\ud z.$$
Hence, when $\lambda=1/ 2$, any $X=(\mathbf{0},W)$, with $W=(\varphi,0,0,0)$ or $(0,0,\varphi,0)$, $\varphi\in L^2(\mathbb R^n)$
 orthogonal to $(-\Delta)^{-1/2}\sigma$,
 lies in $\mathrm{Ker}(\mathscr L-1/ 2)$.
 Therefore, for 
  $\lambda =1/ 2$, 
  the eigenspace is 
   infinite-dimensional.
 Reasoning by a contradiction argument, based on Weyl's criterion, we can show that 
 there is no other values in the essential spectrum of $\mathscr L$, see \cite{SRN3}.
 
From now on, we suppose  $\lambda\neq 1/ 2$.
It allows us to infer $\varpi_0=\varpi_1=0$ and 
$$\varphi_0=\ds\frac{(-\Delta)^{-1/2}\sigma q_0}{\sqrt 2(\lambda-1/2)},\qquad 
\varphi_1=\tau\ds\frac{(-\Delta)^{-1/2}\sigma q_1}{\sqrt 2(\lambda-1/ 2)}.$$
Consequently, 
bearing in mind 
$\int \sigma(-\Delta)^{-1}\sigma\ud z=\kappa$, 
we obtain the following $4\times 4 $ system for $S=(q_0,p_0,q_1,p_1)$,
\[
\lambda S=\begin{pmatrix}
\tau +\ds\frac{ \kappa/2}{\lambda- 1/ 2}& 0 & -1 & 0
\\
0 & \tau & 0 & -1
\\
-1 & 0 & \tau + \ds\frac{ \kappa/2}{\lambda- 1/ 2}& 0
\\
0 & -1 & 0&  \tau
\end{pmatrix}S
.\]
We remark that the relations for $(q_0,q_1)$ and $(p_0,p_1)$ are uncoupled. 
 We start by observing that $\lambda p_0=\tau p_0-p_1$ 
 and $\lambda p_1= -p_0+\tau p_1$
 which admit non trivial solutions provided 
 \[
 (\lambda-\tau)^2-1=\lambda(\lambda -2\tau)=0.\]
 Hence, $0$ and $2\tau$ are eigenvalues for $\mathscr L$ with 
 $\mathrm{Span}(0,1,0,\tau,\mathbf{0})\subset \mathrm{Ker}(\mathscr L)$,
 and $\mathrm{Span}(0,1,0,-\tau,\mathbf{0})\subset \mathrm{Ker}(\mathscr L-2\tau)$, respectively. 
 We turn to the equations for $(q_0,q_1)$
 which admit non trivial solutions provided 
 \[\begin{array}{lll}
 \Big(\lambda - \tau - \ds\frac{ \kappa/2}{\lambda- 1/ 2}\Big)^2-1=
 \Big(\lambda - \tau - \ds\frac{ \kappa/2}{\lambda-1/2}-1\Big)
 \Big(\lambda - \tau - \ds\frac{ \kappa/2}{\lambda- 1/ 2}+1\Big)=0.
 \end{array}\]
 This holds iff
 $(\lambda-  1/2)(\lambda - \tau -1)- \kappa/2=0$ or $(\lambda-  1/ 2)(\lambda - \tau +1)-\kappa/2=0$.
 We distinguish the two cases:
 \begin{itemize}
 \item If $\tau=+1$,
 we get 
 $(\lambda- 1/2)(\lambda - 2)- \kappa/2=0$ or $(\lambda-  1/ 2)\lambda-\kappa/2=0$;
  \item If $\tau=-1$,
 we get 
 $(\lambda-1/ 2)\lambda -\kappa/2=0$ or $(\lambda-  1/2)(\lambda+2)- \kappa/2=0$.
 \end{itemize}
 In both cases, with 
 the second order equation $(\lambda-  1/ 2)\lambda - \kappa/2=\lambda^2-
  \lambda/ 2 - \kappa/2=0$,
 we find the following eigenvalues of opposite signs 
$$\lambda=\ds\frac{1/ 2\pm\sqrt{  1/4+2\kappa}}{2}\in \sigma(\mathscr L).$$
Moreover, 
from $(\lambda- 1/ 2)(\lambda - 2\tau)- \kappa/2=
\lambda^2-( 1/2+2\tau)\lambda+ \tau - \kappa/2
=0$,
we find
\[
\lambda=\ds\frac{  1/2+2\tau\pm \sqrt{(  1/ 2-2\tau)^2+2\kappa}}{2}\in \sigma(\mathscr L).
\]
Hence, when 
$\tau=+1$ with $0<\kappa<2$, this gives two positive eigenvalues; when
$\tau=-1$ or  $\tau=+1$ with $\kappa>2$  we obtain  two eigenvalues of opposite signs.

%
 Finally, $-\mathscr J X_0$ reads $(-1,0,-\tau, 0,\mathbf{0})$. It is orthogonal to $\mathrm{Ker}(\mathscr L)=\mathrm{Span}(X_0,\mathbf{0})$ and it makes sense to consider the equation $\mathscr LY_0=-\mathscr J X_0$.
Imposing $Y_0\in (\mathrm{Ker}(\mathscr L))^\perp$, 
we find 
 \begin{equation*}
  Y_0=\frac1\kappa\left(1, 0, \tau, 0, -\sqrt{2}(-\Delta)^{-1/2}\sigma, 0, -\sqrt{2}(-\Delta)^{-1/2}\sigma, 0\right).
 \end{equation*}
  Accordingly, we get
 $(-\mathscr JX_0|Y_0)=-\frac2\kappa<0$.
 (This product is left unchanged by adding to $Y_0$ any element of $\mathrm{Ker}(\mathscr L)$.)
 \QED
 

  \subsubsection{{Linearization about} the extra solutions when $\kappa>2$}
  
  Let now now assume $\kappa>2$.
We use the same notation as in \eqref{notation3}.
  Considering a perturbation of the solution given by \eqref{extraSW}-\eqref{SpecSolextra}, 
  the linearized equations read
  \[\begin{array}{l}
  i\partial_t v_0=Av_0-v_1+\alpha 
   \ds\int_{\mathbb R^n}\sigma \phi_0 \ud z,
  \\[.3cm]
   i\partial_t v_1=Bv_1-v_0+\beta 
   \ds\int_{\mathbb R^n}\sigma \phi_1 \ud z,
\\[.3cm]
\ds\frac1{c^2}\partial^2_{tt}\phi_0-\Delta \phi_0=-\alpha 
2\sigma\mathrm{Re}(v_0),
\\[.3cm]
\ds\frac1{c^2}\partial^2_{tt}\phi_1-\Delta \phi_1=-\beta
2\sigma\mathrm{Re}(v_1).
\end{array}\]
With the change of variables $$(v_j,=q_j+ip_j,\phi_j)\to \Big(q_j,p_j,\varphi_j=\ds (-\Delta)^{1/2}\phi_j,\varpi_j=\ds\frac{\partial_t \phi_j}{c}
\Big),$$ 
we get
\[
\partial_t X=\mathbb LX
\]
with 
\[\mathbb LX=\begin{pmatrix}
Ap_0-p_1
\\
-Aq_0+q_1-\alpha\ds\int_{\mathbb R^n} (-\Delta)^{-1/2}\sigma \varphi_0\ud z
\\
-p_0+Bp_1
\\
q_0-Bq_1-\beta\ds\int_{\mathbb R^n} (-\Delta)^{-1/2}\sigma \varphi_1\ud z
\\
 c(-\Delta)^{1/2}\varpi_0
\\
- c(-\Delta)^{1/2}\varphi_0-2c\alpha \sigma q_0
\\
 c(-\Delta)^{1/2}\varpi_1
\\
- c(-\Delta)^{1/2}\varphi_1-2c\beta  \sigma q_1
\end{pmatrix}.
\]
We set $\mathbb L=\mathscr J\mathscr L$, with  $\mathscr J$ defined by \eqref{defJ}
and
\begin{equation}\label{defL1}
\mathscr LX=\begin{pmatrix}
Aq_0-q_1+\alpha\ds\int_{\mathbb R^n} (-\Delta)^{-1/2}\sigma \varphi_0\ud z
\\
Ap_0-p_1
\\
Bq_1-q_0+\beta\ds\int_{\mathbb R^n} (-\Delta)^{-1/2}\sigma \varphi_1\ud z
\\
-p_0+Bp_1
\\
\ds\frac{\varphi_0}{2}+\alpha (-\Delta)^{-1/2}\sigma q_0
\\
\ds\frac{\varpi_0}{2}
\\
\ds\frac{\varphi_1}{2}+\beta (-\Delta)^{-1/2}\sigma q_1
\\
\ds\frac{\varpi_1}{2}
\end{pmatrix}.
\end{equation}
We readily obtain the following analog to Proposition~\ref{Spec0LSW}.

\begin{proposition}\label{Spec1LSW}
Let us denote by $\check X$ the vector constructed from $X$ by changing the components $p_j$ and $\varpi_j$ into 
$-p_j$ and $-\varpi_j$.
Let $(\lambda,X)$ be an eigenpair of $\mathbb L$. Then, $(-\lambda, \check X)$, $(\overline \lambda, \overline X) $
and $(-\overline \lambda, \overline {\check X})$ are equally eigenpairs of $\mathbb L$.

Moreover, we can write $\mathbb L=\mathscr J\mathscr L$ with $\mathscr J$ a skew-symmetric operator and $\mathscr L$ a self-adjoint operator.
\end{proposition}
 
 The 
next step consists in studying the spectrum of the self-adjoint operator.

\begin{theo}\label{ThSpec1SW}
Let $\mathscr L$ be the operator defined by \eqref{defL1}.
Then, the following assertions hold:
\begin{enumerate}
\item $\mathrm{Ker}(\mathscr L)=\mathrm{Span}(X_0)$, with $X_0=(S_0,\mathbf{0})$,  $S_0=(0,1,0,A)$;
\item $\sigma_{\mathrm{ess}}(\mathscr L)= \{ 1/2\}$; 
\item $\mathscr L$ has one negative eigenvalue, 
associated to a one-dimensional eigenspace ($n(\mathscr L)=1$);
\item Given $Y_0$  a solution of $\mathscr LY_0=-\mathscr J X_0$,  we have $(-\mathscr JX_0|Y_0)<0$.
\end{enumerate}
\end{theo}

\noindent
{\bf Proof.} The proof of the second item follows exactly the same lines as in Theorem~\ref{ThSpecSW}. 
We also readily check  that $\mathrm{Ker}(\mathscr L)=\mathrm{Span}\{(0,1,0,A,\mathbf{0})\}$.
We have  $-\mathscr J X_0=(-1,0,-A,0)$ and 
solving $\mathscr LY_0=-\mathscr J X_0$ with $Y_0\in (\mathrm{Ker}(\mathscr L))^\perp$ yields 
\begin{equation}
Y_0=\frac{1}{2(A-B)}
(-1,0,A,0, 2\alpha(-\Delta)^{-1/2}\sigma, 0, -2A\beta(-\Delta)^{-1/2}\sigma,0), 
\end{equation}
and thus we get
$(Y_0|-\mathscr JX_0)=-\frac A2<0$.

\begin{figure}[!h]
\begin{center}
\includegraphics[height=7cm]{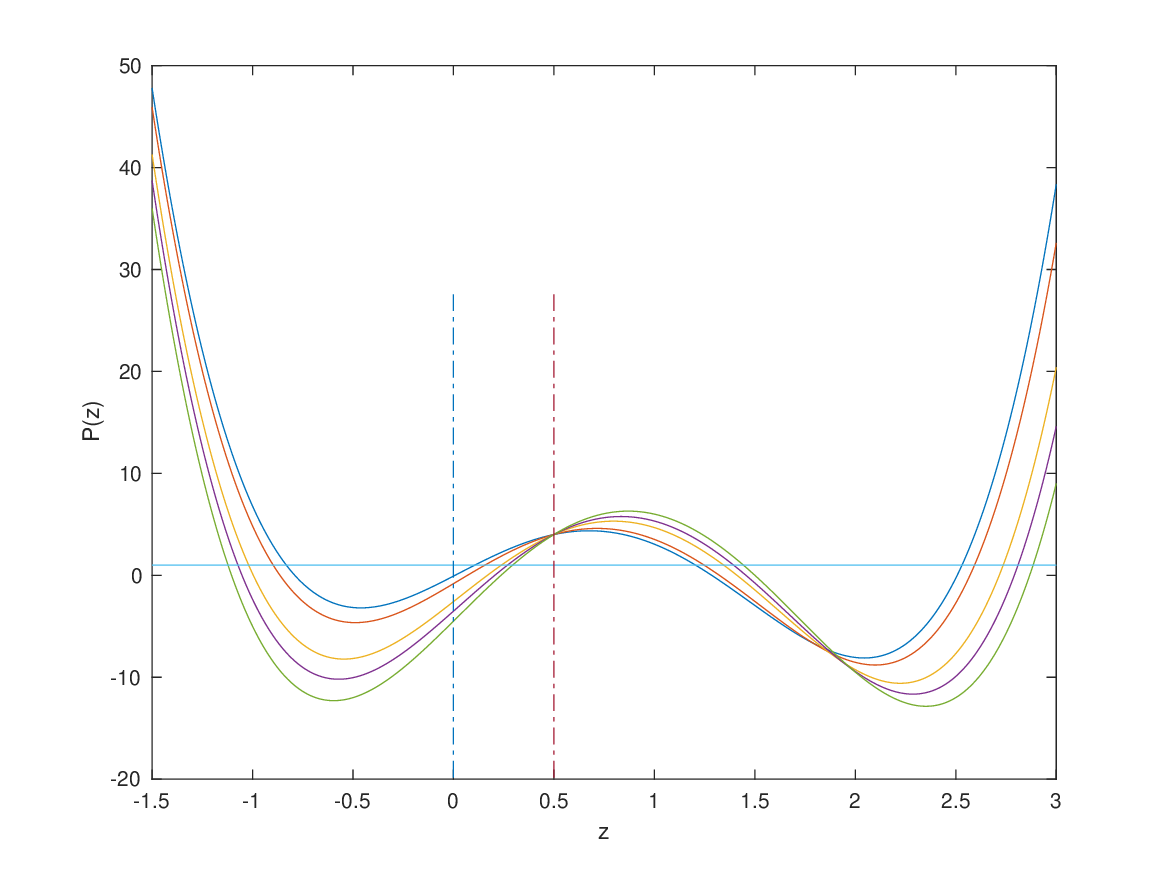}
\caption{Graph of the poynomial function $z\mapsto P(z)$ for several values of $\kappa$
($\kappa\in \{2.01, 2.1, 2.3, 2.4, 2.5\}$}\label{f:p4}
\end{center}
\end{figure}

We now study the eigenvalues   $\lambda\notin\{0,1/2\} 
$ of $\mathscr L$.
We arrive at the matrix system
\[
\begin{pmatrix}
A+\ds\frac{B}{\lambda-1/2} & 0 & -1& 0
\\
0 & A & 0 & -1
\\
-1& 0& B+\ds\frac{A}{\lambda-1/2} & 0
\\
0 & -1& 0& B
\end{pmatrix}\begin{pmatrix}q_0
\\p_0\\q_1\\p_1\end{pmatrix}=\lambda\begin{pmatrix}q_0
\\p_0\\q_1\\p_1\end{pmatrix} .\]
The equations for $(p_0,p_1)$ and $(q_0,q_1)$ uncouple.
The former leads to the relation
\[
\lambda(\lambda-A-B)=0\]
which gives the eigenvalues $0$ and $A+B=\kappa$.
The latter leads to the relation
\begin{align*}
  0=&\left(A+\ds\frac{B}{\lambda-1/2}-\lambda\right)\left(
    B+\ds\frac{A}{\lambda-1/2}-\lambda
    \right)-1\\
    =&\ds\frac{1}{(2\lambda-1)^2}\left( 4\lambda^4 -4(A+B+1)\lambda^3+\lambda^2+(4A^2+4B^2+A+B)\lambda-2(A-B)^2\right)\\
    =&\ds\frac{1}{(2\lambda-1)^2}\left( 4\lambda^4 -4(\kappa+1)\lambda^3+\lambda^2+(4A^2+4B^2+\kappa)\lambda-2(A-B)^2\right)
    =\ds\frac{1}{(2\lambda-1)^2}P(\lambda)
\end{align*}
with $P$ a fourth order polynomial.
Descartes' rule of sign then tells us that 
$P$ has exactly one negative 
root,
see Fig.~\ref{f:p4}.
We have thus proved the third item in Theorem~\ref{ThSpec1SW}.
\QED
 
\subsection{Spectral and linearized stability}

We start with the study of the spectral stability of the solution \eqref{SpecSol} of \eqref{S1}-\eqref{W1}. Let $\mathscr{L}$ be defined by  \eqref{defL}. According to \cite{ChouPel}, we introduce the operator 
\[\mathscr M=-\mathscr J\mathscr L\mathscr J,\qquad
\mathbb A=\mathscr P\mathscr M\mathscr P,\]
where $\mathscr P$ is the orthogonal projection on $(\mathrm{Ker}(\mathscr L))^\perp$, and we set 
 $$\mathbb K=\mathscr P\mathscr L^{-1}\mathscr P.$$
We are interested in the generalized 
eigenvalue problem
\[\mathscr MX
=\mu \tilde X,\qquad \mathscr L\tilde X=X.\]
Recall that $X$ has to belong to $(\mathrm{Ker}(\mathscr{L}))^\perp$ and
we need to compute the product $(\mathbb KX|X)=(\tilde X|X)$, 
which is thus left unchanged by adding to $\tilde X$ an element in $\mathrm{Ker}(\mathscr{L})$.
Hence,  $\tilde X$ can be chosen in $(\mathrm{Ker}(\mathscr{L}))^\perp$.
Solving the generalized eigenvalue problem amounts to solve 
\begin{align*}
  &\tau q_0 -q_1=\mu \tilde q_0 & &\tau p_0 - p_1+\ds c\sqrt2 \ds\int_{\mathbb R^n}
  \sigma \varpi_0\ud z=\mu \tilde p_0,\\
  & -q_0+\tau q_1=\mu \tilde q_1, & &-p_0 +\tau  p_1+\tau\ds c\sqrt2\ds\int_{\mathbb R^n}
  \sigma \varpi_1\ud z=\mu \tilde p_1\\
  &2 c^2(-\Delta \varphi_0)
=\mu \tilde \varphi_0 & & 2c^2(-\Delta \varpi_0)+c\sqrt2\sigma p_0=\mu \tilde \varpi_0\\
& 2 c^2(-\Delta \varphi_1)
=\mu \tilde \varphi_1 & & 2c^2(-\Delta \varpi_1)+\tau c\sqrt2\sigma p_1=\mu \tilde \varpi_1
\end{align*}
coupled to 
\begin{align*}
  &\tau  \tilde q_0-\tilde q_1
  + \ds\frac1{\sqrt 2} \ds\int_{\mathbb R^n} \sigma (-\Delta)^{-1/2} \tilde \varphi_0\ud z=q_0 & & \tau  \tilde p_0 -\tilde p_1=p_0\\
  &- \tilde q_0+\tau  \tilde q_1
  +\ds\frac\tau {\sqrt 2}\ds\int_{\mathbb R^n} \sigma (-\Delta)^{-1/2} \tilde \varphi_1\ud z=q_1& &-p_0+\tau  \tilde p_1=p_1\\
  &\ds\frac{ \tilde \varphi_0}{2}+  \ds\frac1 {\sqrt 2}(-\Delta)^{-1/2}\sigma \tilde q_0=\varphi_0,& &\ds\frac { \tilde\varpi_0}{2}=\varpi_0,\\
  &\ds\frac{ \tilde \varphi_1}{2}+  \ds\frac\tau {\sqrt 2}(-\Delta)^{-1/2}\sigma \tilde q_1=\varphi_1,& &\ds\frac { \tilde\varpi_1}{2}=\varpi_1.
\end{align*}
This leads to the following relations
\begin{align*}
  &\Big(-\ds\frac{\mu}{c^2} -\Delta\Big)\varpi_0=-\ds\frac{1}{c\sqrt2}\sigma p_0,&
   \Big(-\ds\frac{\mu}{c^2} -\Delta\Big)\varpi_1=-\tau\ds\frac{1}{c\sqrt 2} \sigma p_1,
\end{align*}
and
\begin{align*}
  &\Big(-\ds\frac{\mu}{c^2} -\Delta\Big)\tilde \varphi_0=-\sqrt 2(-\Delta)^{1/2}\sigma \tilde q_0,&
\Big(-\ds\frac{\mu}{c^2} -\Delta\Big)\tilde \varphi_1=-\tau\sqrt 2(-\Delta)^{1/2}\sigma \tilde q_1.
\end{align*}

When $\mu <0$, these equations can be solved by means of the Fourier transform and we get
\begin{align*}
  &\widehat\varpi_0(\xi)=-\ds\frac{1}{c\sqrt{2}}\ds\frac{\widehat \sigma(\xi)}{|\mu|/c^2+|\xi|^2}p_0,& & 
\widehat\varpi_1(\xi)=-\ds\frac{\tau}{c\sqrt{2}}\ds\frac{\widehat \sigma(\xi)}{|\mu|/c^2+|\xi|^2}p_1,\\
&\widehat{\tilde \varphi_0}=-\ds\frac{\sqrt2|\xi|\widehat \sigma(\xi)}{|\mu|/c^2+|\xi|^2}\tilde q_0,\qquad
& &
\widehat{\tilde \varphi_1}=-\ds\frac{\tau\sqrt2 |\xi|\widehat \sigma(\xi)}{|\mu|/c^2+|\xi|^2}\tilde q_1.
\end{align*}
It follows that
\begin{align*}
&{\frac{1}{\sqrt2}}\ds\int_{\mathbb R^n} \sigma(-\Delta)^{-1/2}\tilde \varphi_0\ud z=-\tilde q_0
\underbrace{
\ds\int_{\mathbb R^n} \ds\frac{|\widehat \sigma(\xi)|^2}
{|\mu|/c^2+|\xi|^2}\ds\frac{\ud\xi}{(2\pi)^n}}_{=\kappa_{|\mu|/c^2}}
,&&
\frac{\tau}{\sqrt2}\ds\int_{\mathbb R^n} \sigma(-\Delta)^{-1/2}\tilde \phi_1\ud z
=-\tilde q_1\kappa_{|\mu|/c^2},
\\
&c\sqrt{2}\ds\int_{\mathbb R^n} \sigma \varpi_0\ud z=-p_0\kappa_{|\mu|/c^2}
&&
\tau c\sqrt{2}\ds\int_{\mathbb R^n} \sigma \varpi_1\ud z=- p_1\kappa_{|\mu|/c^2}.
\end{align*}
With the matrices defined in \eqref{defmk}, we are thus led to
\[
M_0\begin{pmatrix}
\tilde p_0\\\tilde p_1\end{pmatrix}=\begin{pmatrix}
 p_0\\ p_1\end{pmatrix},\qquad
 M_{\kappa_{|\mu|/c^2}}\begin{pmatrix}
 p_0\\ p_1\end{pmatrix}=\mu \begin{pmatrix}
\tilde p_0\\\tilde p_1\end{pmatrix},
\]
together with
\[
M_{\kappa_{|\mu|/c^2}}\begin{pmatrix}
\tilde q_0\\\tilde q_1\end{pmatrix}=\begin{pmatrix}
 q_0\\ q_1\end{pmatrix},\qquad
 M_0\begin{pmatrix}
 q_0\\ q_1\end{pmatrix}=\mu \begin{pmatrix}
\tilde q_0\\\tilde q_1\end{pmatrix}.
\]
Since $M_{\kappa_{|\mu|/c^2}} M_0=
M_0M_{\kappa_{|\mu|/c^2}}$, we deduce, like for the asymptotic model, that $\mu<0$ should be  such that $\mathrm{det}(M_0M_{\kappa_{|\mu|/c^2}}-\mu\mathbb{I})=0$. This condition leads to 
\begin{align*}
  0&=(2-\tau \kappa_\gamma +\gamma c^2)^2-(\kappa_\gamma-2\tau)^2=(2-\tau \kappa_\gamma +\gamma c^2-\kappa_\gamma+2\tau)(2-\tau \kappa_\gamma +\gamma c^2+\kappa_\gamma-2\tau)\\
  &=\gamma c^2 (2(2-\tau\kappa_\gamma)+\gamma c^2)
\end{align*} 
where we set $\gamma=-\frac{\mu}{c^2}=\frac{|\mu|}{c^2}$.
 When $\tau=-1$ or $\tau=+1$ with $0<\kappa<2$, we have $2(2-\tau\kappa_\gamma)+\gamma c^2>0$ for any positive $\gamma$, hence there is no solution
to this equation:
 in these cases we have $N_n^-=0$.
%
If $\tau=+1$ and $\kappa>2$, 
it is thus required to make the non linear quantity 
\[F(\gamma)=\gamma-\ds\frac{2}{c^2}(\kappa_\gamma-2)
\]
vanishes. The function $F$ is continuous, increasing from $(0,\infty)$ to $(-\frac2{c^2}(\kappa-2),+\infty)$; 
hence there exists a unique $\gamma_c=-\frac{\mu_c}{c^2}>0$ such that $F(\gamma_c)=0$. 
Finally, we have to compute $(\mathbb KX, X)$. Since $X\in (\mathrm{Ker}(\mathscr{L}))^\perp$, $\mathcal P X=X$ and $(\mathbb KX, X)=(\tilde X,X)$. Using the equations for $(\tilde q_0,\tilde q_1)$ and $(\tilde p_0,\tilde p_1)$ together with $\gamma_c c^2=2(\kappa_{\gamma_c}-2)>0$, we deduce that the eigenvectors associated to $\mu_c$ are such that $\tilde q_1=-\tilde q_0$ and $\tilde p_1=-\tilde p_0$. On the one hand, choosing $\tilde q_0=1$ and $\tilde p_0=0$, we have 
\begin{align*}
  &\tilde X= \left(1,0,-1,0, \mathcal F^{-1}\left(-\ds\frac{\sqrt2|\xi|\widehat \sigma(\xi)}{\gamma_c+|\xi|^2}\right),0, \mathcal F^{-1}\left(\ds\frac{\sqrt2|\xi|\widehat \sigma(\xi)}{\gamma_c+|\xi|^2}\right),0\right)\\
  &X=\left(2-\kappa_{\gamma_c},0,\kappa_{\gamma_c}-2,0,\frac{1}{\sqrt2}\mathcal F^{-1}\left(-\ds\frac{|\xi|\widehat \sigma(\xi)}{\gamma_c+|\xi|^2}+\frac{\widehat \sigma(\xi)}{|\xi|}\right),0,\frac{1}{\sqrt2}\mathcal F^{-1}\left(\ds\frac{|\xi|\widehat \sigma(\xi)}{\gamma_c+|\xi|^2}-\frac{\widehat \sigma(\xi)}{|\xi|}\right),0\right)
\end{align*}
so that 
\begin{align*}
  (\tilde X, X)&=-2(\kappa_{\gamma_c} -2)+2\int_{\R^n}\left(\frac{|\xi|^2|\widehat \sigma(\xi)|^2}{(\gamma_c+|\xi|^2)^2}-\frac{|\widehat \sigma(\xi)|^2}{\gamma_c+|\xi|^2}\right)\ds\frac{\ud\xi}{(2\pi)^n}\\
  &= -\gamma_c c^2-2\gamma_c\int_{\R^n}\frac{|\widehat \sigma(\xi)|^2}{(\gamma_c+|\xi|^2)^2}\ds\frac{\ud\xi}{(2\pi)^n}<0.
\end{align*}
On the other hand, choosing $\tilde q_0=0$ and $\tilde p_0=1$, we have 
\begin{align*}
  &\tilde X= \left(0,1,0,-1, 0, -\frac{2\sqrt2}{c}\mathcal F^{-1}\left(\frac{\widehat{\sigma}(\xi)}{\gamma_c+|\xi|^2}\right),0,\frac{2\sqrt2}{c}\mathcal F^{-1}\left(\frac{\widehat{\sigma}(\xi)}{\gamma_c+|\xi|^2}\right)\right)\\
  &X=\left(0,2,0,-2,0,-\frac{\sqrt2}{c}\mathcal F^{-1}\left(\frac{\widehat{\sigma}(\xi)}{\gamma_c+|\xi|^2}\right),0,\frac{\sqrt2}{c}\mathcal F^{-1}\left(\frac{\widehat{\sigma}(\xi)}{\gamma_c+|\xi|^2}\right)\right)
\end{align*}
so that 
\begin{align*}
  (\tilde X, X)&=4+\frac{8}{c^2} \int_{\R^n}\frac{|\widehat \sigma(\xi)|^2}{(\gamma_c+|\xi|^2)^2}\ds\frac{\ud\xi}{(2\pi)^n}>0.
\end{align*}
We can conclude $N_n^-=1$.

When $\mu>0$, the symbol $\frac{\widehat \sigma(\xi)}{|\xi|^2-\mu/c^2}$
 has a singularity which is non square integrable; 
 this forces to set $p_0=p_1=0$, and $\tilde q_0=\tilde q_1=0$, so that $\varpi_0=\varpi_1=0$,
  and
 $\tilde \phi_0=\tilde \phi_1=0$. 
It implies $q_0=q_1=0$ and $\tilde p_0=\tilde p_1=0$; there is no non trivial solution of the generalized eigenvalue problem with $\mu>0$, that is $N_n^+=0.$
\\

For $\mu=0$, the equations reduce to 
\begin{align*}
  &\tau q_0 -q_1=0 & &\tau p_0 - p_1+\ds c\sqrt2 \ds\int_{\mathbb R^n}
  \sigma \varpi_0\ud z=0,\\
  & -q_0+\tau q_1=0, & &-p_0 +\tau  p_1+\tau\ds c\sqrt2\ds\int_{\mathbb R^n}
  \sigma \varpi_1\ud z=0\\
  &2c^2(-\Delta \varphi_0)
=0 & & 2c^2(-\Delta \varpi_0)+c\sqrt2\sigma p_0=0\\
&2 c^2(-\Delta \varphi_1)
=0 & & 2c^2(-\Delta \varpi_1)+\tau c\sqrt2\sigma p_1=0
\end{align*}
It yields $\varphi_0=\varphi_1=0$ and 
$\varpi_0=-\frac{1}{c\sqrt 2}(-\Delta)^{-1}\sigma p_0$, $\varpi_1=-\frac{\tau}{c\sqrt 2}(-\Delta)^{-1}\sigma p_1$, hence
the systems
\[
M_0\begin{pmatrix}q_0\\ q_1\end{pmatrix}=0,\qquad
M_\kappa\begin{pmatrix}p_0\\ p_1\end{pmatrix}=0.
\]
Solving these systems, we obtain $p_0=0=p_1$ and $q_1=\tau q_0$. As a consequence $X$ is proportional to $-\mathscr J X_0$. Reinterpreting $(\mathbb K X,X)=(\mathbb K(-\mathscr J X_0)|-\mathscr J X_0)$ as $(-\mathscr J X_0|Y_0)$ with $Y_0$ given in Theorem \ref{ThSpecSW}, we obtain $(\mathbb K(-\mathscr J X_0)|-\mathscr J X_0)<0$ and 
we conclude that $N^0_n=1$.

To sum up, we have the following 
\begin{align*}
  N^0_n=1,\ds N_n^+=0 \text{ and } N_n^-=\left\{\begin{aligned}
    0 & \text{ if } \tau=-1 \text{ or } \tau=1 \text{ and } \kappa <2,\\
    1 & \text{ if } \tau=1 \text{ and } \kappa >2.
  \end{aligned}\right..
\end{align*}

We remind the reader that the spectral stability means that  the spectrum of $\mathbb L$ is contained in $i\R$.  
To derive information about $\sigma(\mathbb L)$, we  use the counting argument of \cite[Theorem~1]{ChouPel} (see also \cite{LZ}) which asserts that
\[N^-_n+N^0_n+N^+_n+N_{C^+}=n(\mathscr L).\]
The presence of spectrally unstable directions corresponds to $N^-_n\neq 0$ or $N_{C^+}\neq 0$.
Gathering the obtained information, we infer that 
\begin{itemize}
  \item if $\tau =1$ and $\kappa<2$, $N_n^-=0$ and $N_{C^+}=n(\mathscr L)-1=0$;
  \item if $\tau =1$ and $\kappa>2$, $N_n^-=1$ and $N_{C^+}=n(\mathscr L)-1-1=0$;
  \item if $\tau =-1$, $N_n^-=0$ and $N_{C^+}=n(\mathscr L)-1=2$.
\end{itemize}
Accordingly, we conclude with the following claim.

\begin{proposition}\label{prop:spectralstabSW}
Suppose $0<\kappa<2$ and let $\tau=+1$. Then, the reference solution \eqref{SpecSol} of \eqref{S1}-\eqref{W1} is spectrally stable.
If $\tau=-1$ or $\tau=+1$ with $\kappa>2$, the reference solution \eqref{SpecSol} of \eqref{S1}-\eqref{W1} is spectrally unstable.
\end{proposition}

\begin{figure}[!h]
 \begin{subfigure}{0.48\textwidth}
  \includegraphics[height=5cm]{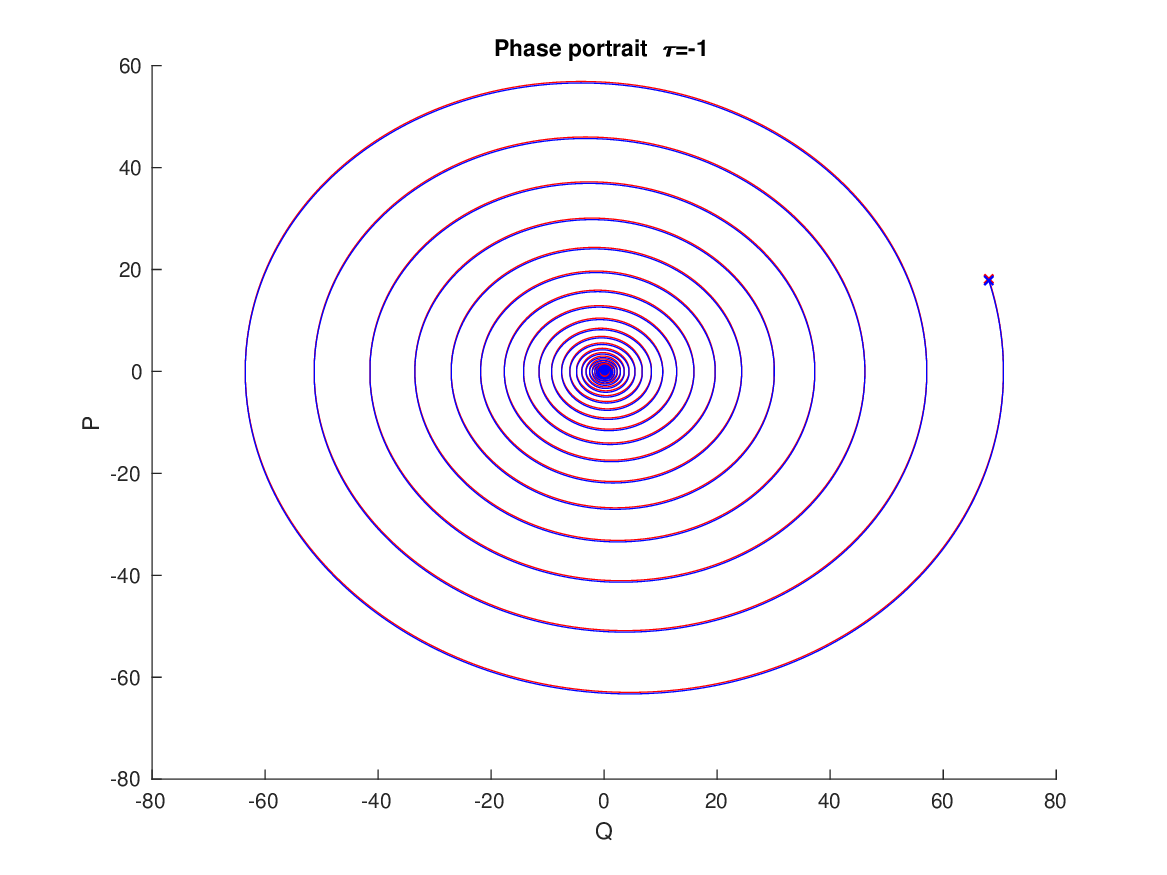}
  \caption{}
  \end{subfigure}
  \begin{subfigure}{0.48\textwidth}
   \includegraphics[height=5cm]{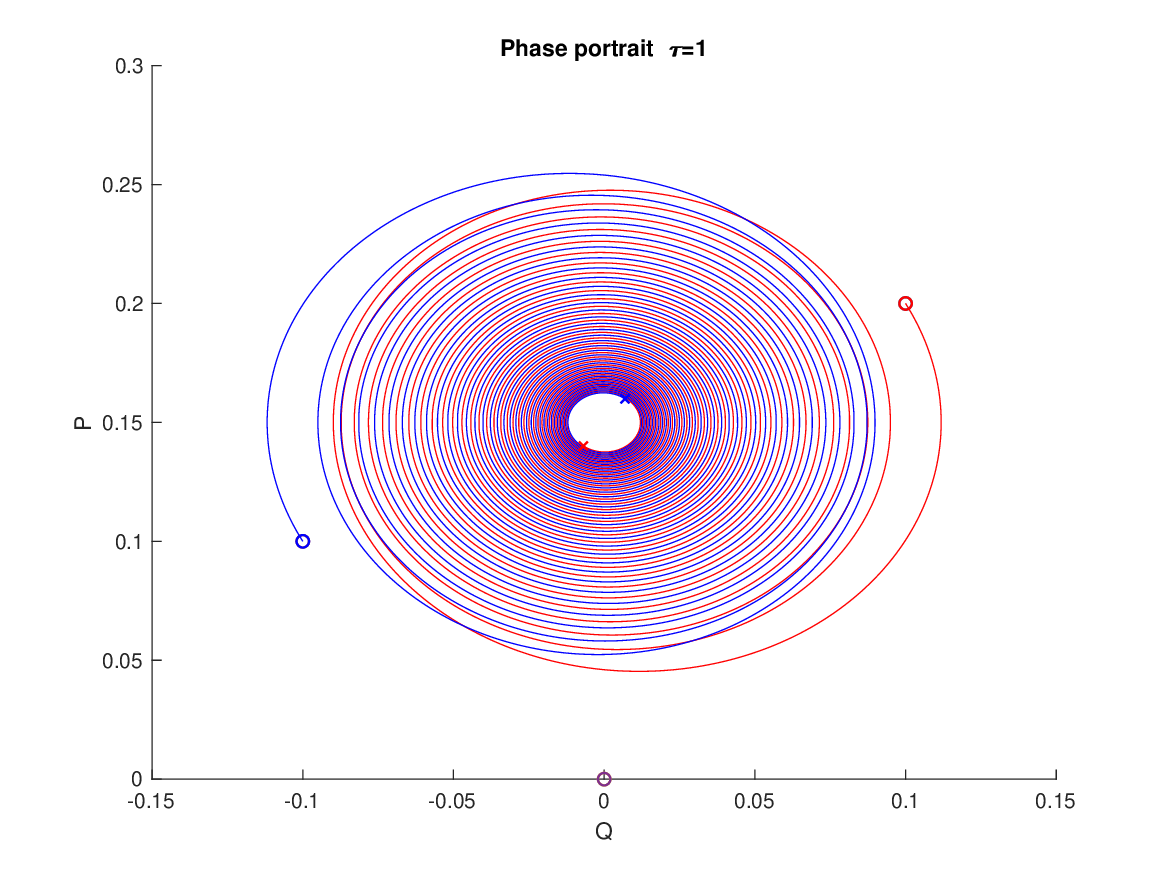}
   \caption{}
   \end{subfigure}
     \begin{subfigure}{0.52\textwidth}
   \includegraphics[height=5cm]{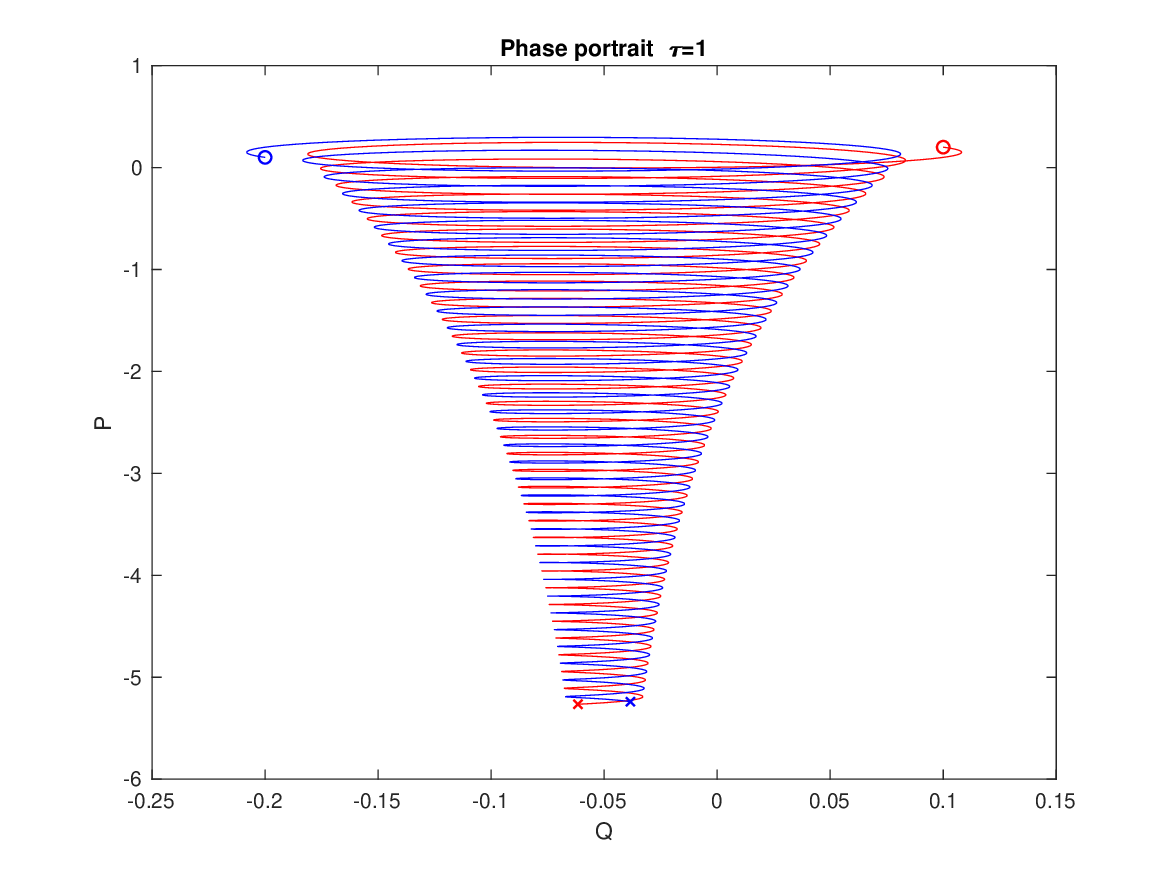}
    \caption{}
   \end{subfigure}
    \begin{subfigure}{0.48\textwidth}
   \includegraphics[height=5cm]{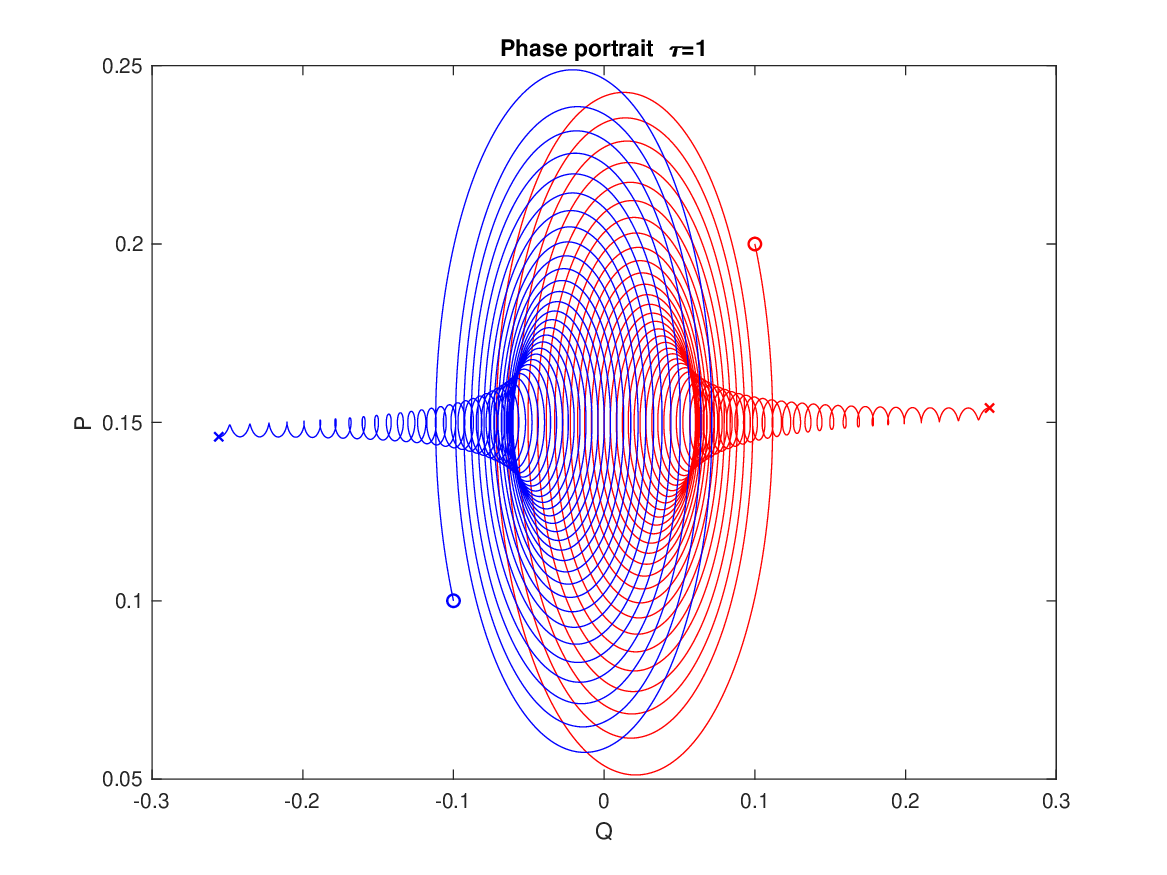}
    \caption{}
   \end{subfigure}
  \caption{Simulation of the  linearized coupled model: phase portrait at $T=100$, for $\kappa=1.1$ 
   with  $\tau=-1$ (a),
   with  $\tau=1$  for well-prepared initial data (b),
   with $\tau=1$  for ill-prepared initial data (c), 
   and with $\kappa=2.0688$, $\tau =1$ at $T=150$ for well-prepared initial data (d).
   The circled points indicate the initial state, the cross indicate the final state
  } 
  \label{NumSWLin}
  \end{figure}

This result is illustrated in Figure~\ref{NumSWLin}.
Inspired by the asymptotic problem, see Proposition~\ref{staLinHaEx} and Figure~\ref{stalin}, 
we guess that the linearized stabiltiy 
requires suitable orthogonality conditions.
Indeed, we can check  that $\mathrm{Ker}(\mathbb L)=\mathrm{Span}(0,1,0,\tau, \mathbf{0})$ 
and $\mathrm{Ker}(\mathbb L^*)=\mathrm{Span}(1,0,\tau,0, \mathbf{0})$.
In particular, denoting $\Psi=(1,0,\tau,0, \mathbf{0})$, we have $\frac{\ud}{\ud t}(X|\Psi)=0$, and in order to prevent 
grows of the linearized solution, we select initial data such that $(X_{\mathrm{init}}|\Psi)=0$, which reduces to 
$q_{\mathrm{init},0}+\tau q_{\mathrm{init},1}=0$.\\

 \begin{figure}[!h]
 \begin{center}
  \begin{subfigure}{0.48\textwidth}
  \includegraphics[height=4cm]{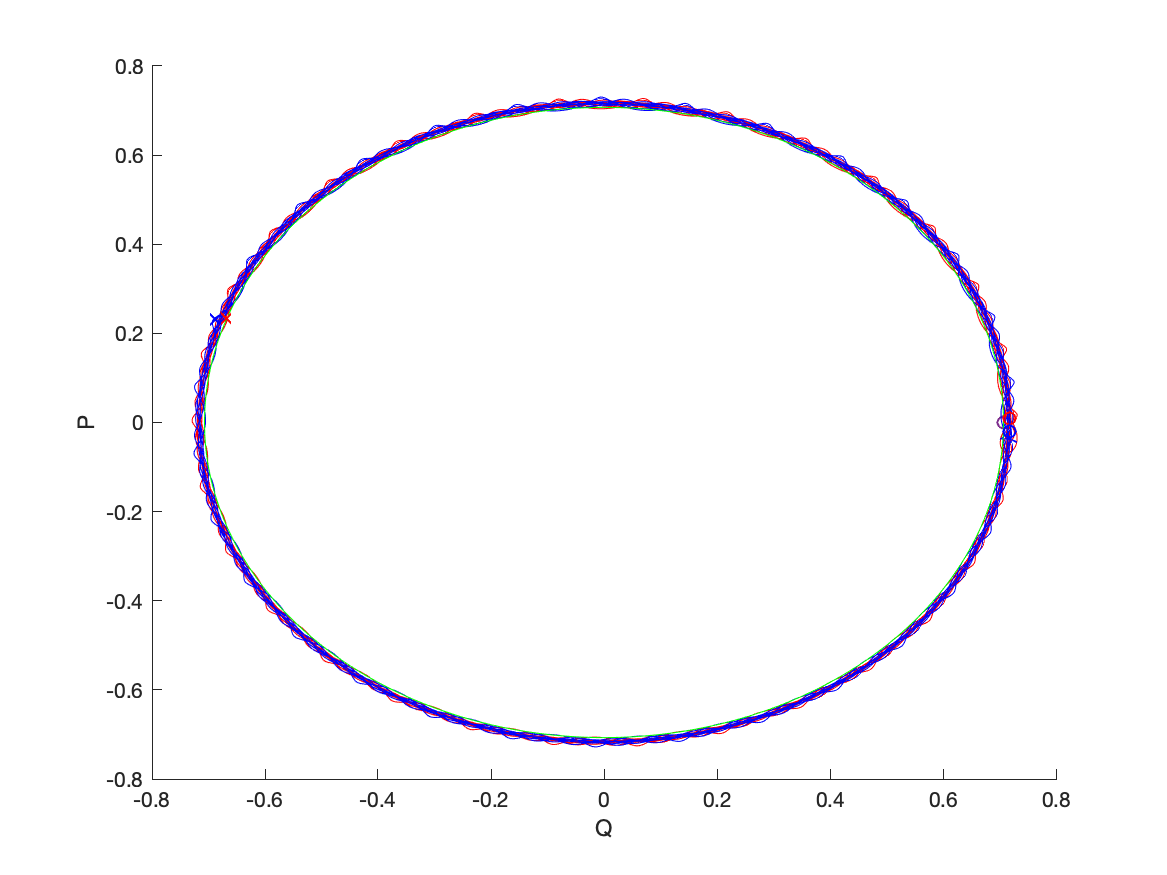}
   \caption{}
  \end{subfigure}
  \begin{subfigure}{0.48\textwidth}
   \includegraphics[height=4cm]{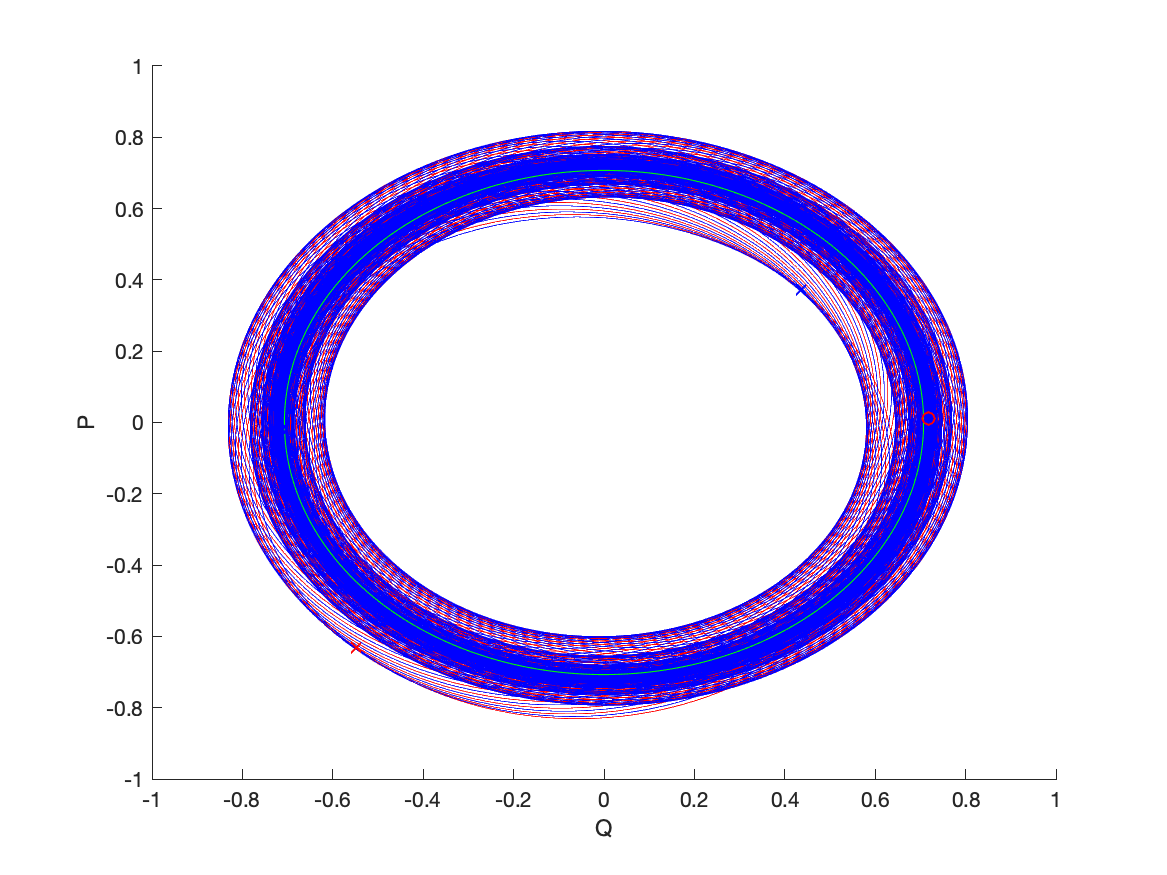}
    \caption{}
  \end{subfigure}
  \begin{subfigure}{0.48\textwidth}
    \includegraphics[height=4cm]{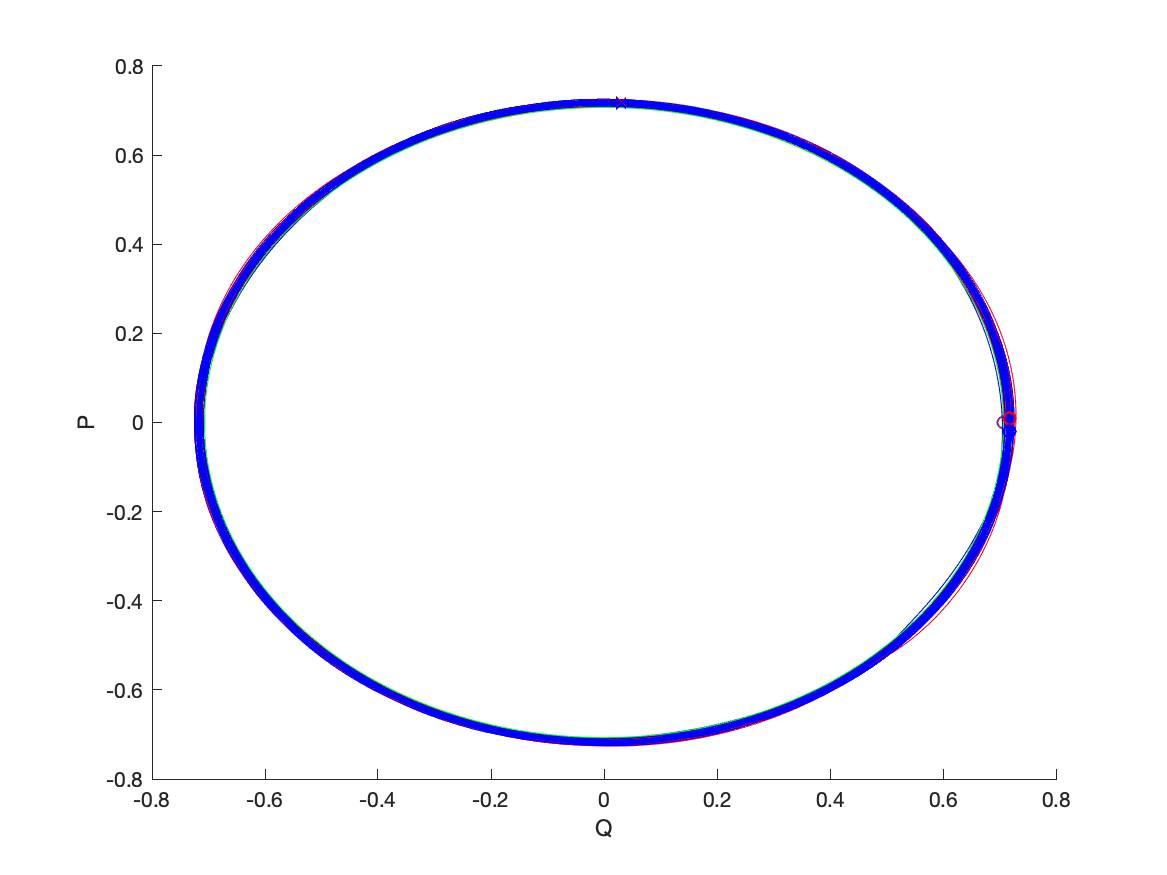}
     \caption{}
  \end{subfigure}
  \begin{subfigure}{0.48\textwidth}
    \includegraphics[height=4cm]{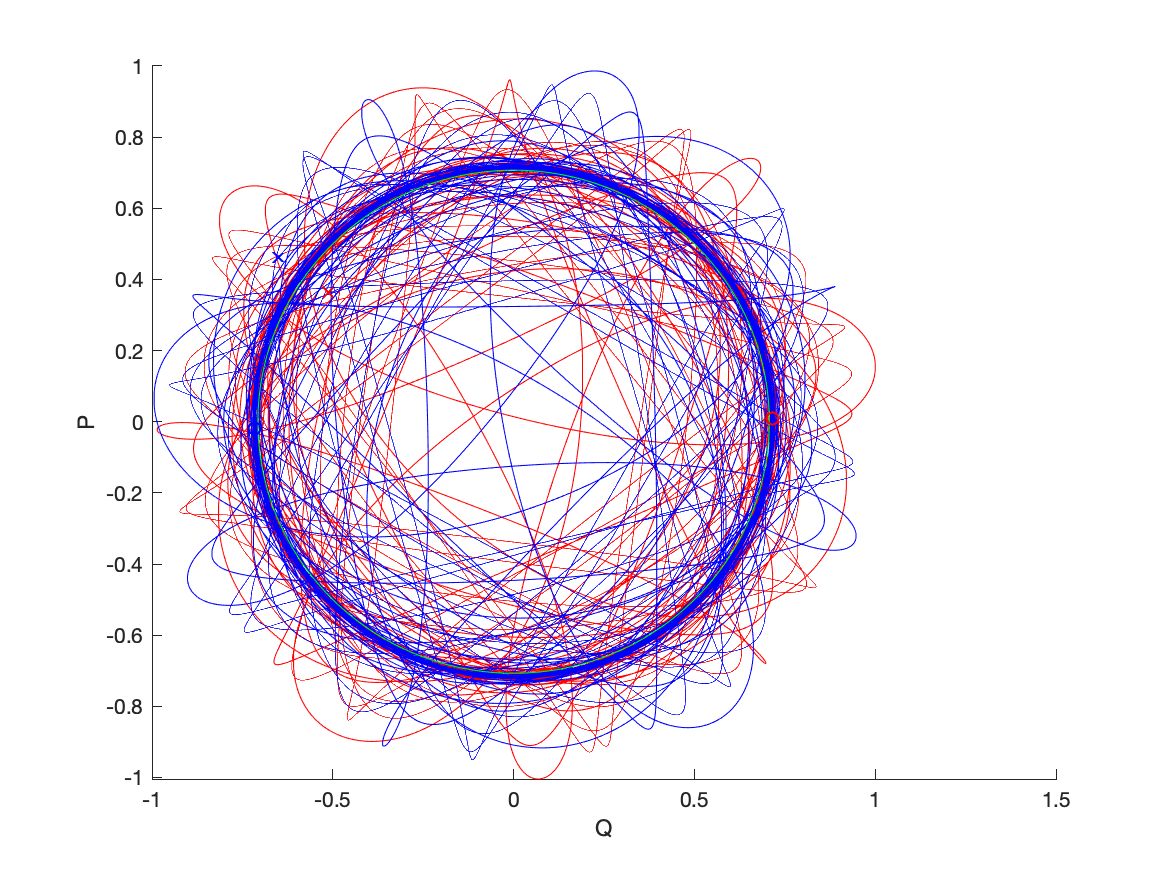}
 \caption{}
  \end{subfigure}
  \begin{subfigure}{0.48\textwidth}
    \includegraphics[height=4cm]{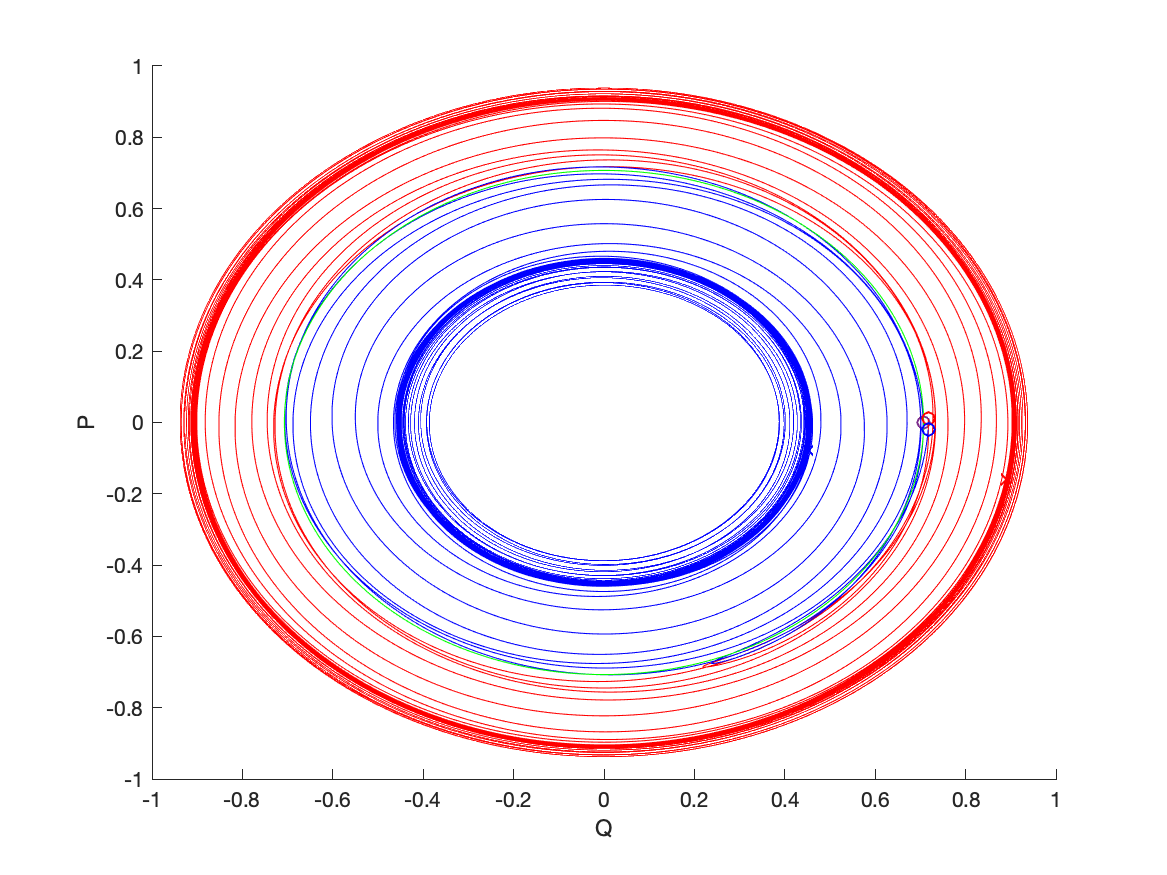} 
    \caption{}
  \end{subfigure}
  \begin{subfigure}{0.48\textwidth}
    \includegraphics[height=4cm]{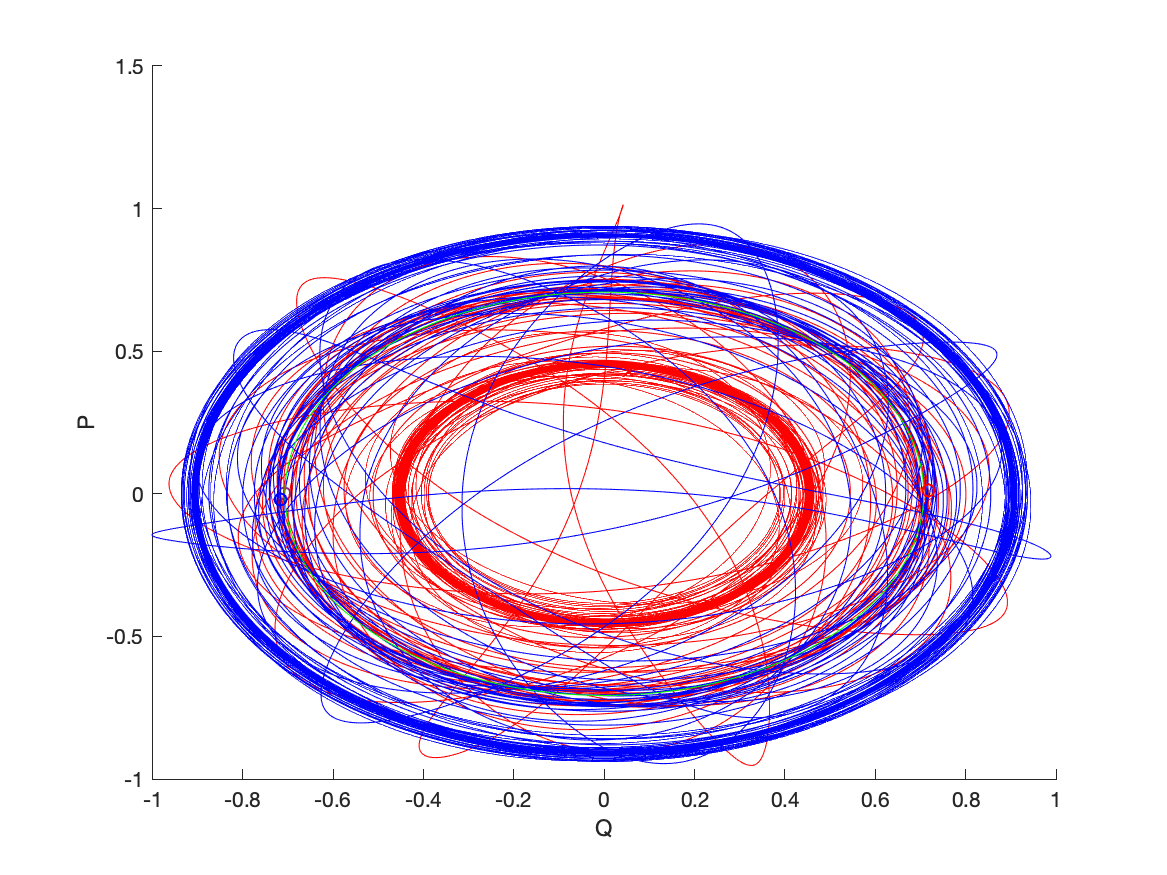}
    \caption{}
     \end{subfigure}
  \caption{Simulation of the  coupled model: phase portrait at $T=700$,  
   with  $\tau=+1$ (a, c, e),
   with  $\tau=-1$ (b, d, f), and several values of $\kappa$:  $\kappa=0.193$ for (a, b), $\kappa=1.58$ for (c, d), $\kappa=2.42$ (e, f).
   The circled points indicate the initial state, the cross indicate the final state
  } 
  \end{center}
  \end{figure}


When $\kappa>2$, a similar statement holds for the solution \eqref{extraSW}.

\begin{proposition}\label{prop:spectralstabSWextra}
Suppose $\kappa>2$. Then, the reference solution \eqref{SpecSolextra} of \eqref{S1}-\eqref{W1} is spectrally stable.
\end{proposition}

\noindent
{\bf Proof.}
As before, we are concerned with the  generalized 
eigenvalue problem
$\mathscr MX
=\mu \tilde X$, $ \mathscr L\tilde X=X$ with $X, \tilde X \in (\mathrm{Ker}(\mathscr L))^\perp$.
It now takes the form
\begin{align*}
  &A q_0 -q_1=\mu \tilde q_0 & & A p_0 - p_1+2\alpha c\ds\int_{\mathbb R^n}
  \sigma \varpi_0\ud z=\mu \tilde p_0,\\
  &-q_0+ Bq_1=\mu \tilde q_1,& &-p_0 +B  p_1+2\beta c \ds\int_{\mathbb R^n}
  \sigma \varpi_1\ud z=\mu \tilde p_1,\\
  &2c^2(-\Delta \varphi_0)
  =\mu \tilde \varphi_0& & 2c^2(-\Delta \varpi_0)+2 \alpha c \ds\sigma p_0=\mu \tilde \varpi_0,\\
  &2c^2(-\Delta \varphi_1)
=\mu\tilde \varphi_1, & &2c^2(- \Delta \varpi_1)+2 \beta c\ds \sigma p_1=\mu \tilde \varpi_1,
\end{align*}

coupled to 
\begin{align*}
  &A \tilde q_0 -\tilde q_1
 +\alpha \ds\int_{\mathbb R^n} \sigma (-\Delta)^{-1/2} \tilde \varphi_0\ud z=q_0,
 & & 
 A  \tilde p_0 -\tilde p_1=p_0
 \\
 &- \tilde q_0+ B \tilde q_1
 +\beta \ds\int_{\mathbb R^n} \sigma (-\Delta)^{-1/2} \tilde \varphi_1\ud z=q_1,
 & & 
- \tilde p_0+ B \tilde p_1=p_1,
 \\
 &\frac{\tilde \varphi_0}{2}+\alpha (-\Delta)^{-1/2}\sigma \tilde q_0=\varphi_0,
 & & 
 \frac{\tilde\varpi_0}{2}=\varpi_0,
 \\
 &\frac{\tilde \varphi_1}{2}+\beta (-\Delta)^{-1/2}\sigma \tilde q_1=\varphi_1,
 & & 
 \frac{\tilde\varpi_1}{2}=\varpi_1.
\end{align*}

As before the operator $(-\frac{\mu}{c^2}-\Delta)$ plays a crucial role. In particular, if $\mu>0$ it cannot be inverted so that the possibility to find non trivial solutions with $\mu>0$ is exluded. As a consequence, $N^+_n=0$.
When $\mu<0$, 
 we have
\begin{align*}
  &\left(-\Delta -\frac{\mu}{c^2}\right)\varpi_0=- \frac{\alpha}{c}  \ds\sigma p_0, & \left(-\Delta -\frac{\mu}{c^2}\right)\varpi_1=- \frac{\beta}{c}  \ds\sigma p_1,
\end{align*}
and 
\begin{align*}
  &\left(-\Delta -\frac{\mu}{c^2}\right)\tilde \varphi_0=- 2\alpha  \ds(-\Delta)^{1/2}\sigma \tilde q_0, & \left(-\Delta -\frac{\mu}{c^2}\right)\tilde \varphi_1=- 2\beta  \ds(-\Delta)^{1/2}\sigma \tilde q_1.
\end{align*}
Setting $\gamma=-\frac{\mu}{c^2}>0$ and $\kappa_\gamma$ as before leads to the following systems of equations
\begin{align*}
  &\begin{pmatrix}
    A-2\alpha^2 \kappa_\gamma & -1\\
    -1 & B-2\beta^2 \kappa_\gamma
  \end{pmatrix}
  \begin{pmatrix}
    p_0\\ p_1
  \end{pmatrix}
  =\mu \begin{pmatrix}
    \tilde p_0\\ \tilde p_1
  \end{pmatrix}, & \begin{pmatrix}
    A & -1\\
    -1 & B
  \end{pmatrix}
  \begin{pmatrix}
    \tilde p_0\\ \tilde p_1
  \end{pmatrix}
  =\begin{pmatrix}
    p_0\\ p_1
  \end{pmatrix}
\end{align*}
and 
\begin{align*}
  & \begin{pmatrix}
    A & -1\\
    -1 & B
  \end{pmatrix}
  \begin{pmatrix}
    q_0\\ q_1
  \end{pmatrix}
  =\mu \begin{pmatrix}
    \tilde q_0\\ \tilde q_1
  \end{pmatrix},& \begin{pmatrix}
    A-2\alpha^2 \kappa_\gamma & -1\\
    -1 & B-2\beta^2 \kappa_\gamma
  \end{pmatrix}
  \begin{pmatrix}
    \tilde q_0\\ \tilde q_1
  \end{pmatrix}
  =\begin{pmatrix}
     q_0\\ q_1
  \end{pmatrix}.
\end{align*}
As before, $\mu <0$ should be such that 
\begin{align*}
  \mathrm{det}\left(\begin{pmatrix}
    A & -1\\
    -1 & B
  \end{pmatrix}\begin{pmatrix}
    A-2\alpha^2 \kappa_\gamma & -1\\
    -1 & B-2\beta^2 \kappa_\gamma
  \end{pmatrix}-\mu\mathbb{I}\right)=0
\end{align*}
This condition is equivalent to
\begin{align*}
  0&=\mathrm{det}\left(\begin{pmatrix}
    A & -1\\
    -1 & B
  \end{pmatrix}\begin{pmatrix}
    A-2\frac{B}{\kappa} \kappa_\gamma & -1\\
    -1 & B-2\frac{A}{\kappa} \kappa_\gamma
  \end{pmatrix}-\mu\mathbb{I}\right)\\
  &= \mathrm{det}\left(\begin{pmatrix}
    A\kappa-2\frac{\kappa_\gamma}{\kappa}  & -\kappa+2A\frac{\kappa_\gamma}{\kappa} \\
    -\kappa+2B\frac{\kappa_\gamma}{\kappa} & B\kappa-2\frac{\kappa_\gamma}{\kappa} 
  \end{pmatrix}+\gamma c^2\mathbb{I}\right)\\
  &=\left(A\kappa-2\frac{\kappa_\gamma}{\kappa}+\gamma c^2\right)\left(B\kappa-2\frac{\kappa_\gamma}{\kappa}+\gamma c^2\right)-\left(2B\frac{\kappa_\gamma}{\kappa}-\kappa\right)\left(2A\frac{\kappa_\gamma}{\kappa}-\kappa\right)\\
  &= (A+B)\kappa\gamma c^2 -4\frac{\kappa_\gamma}{\kappa}\gamma c^2+(\gamma c^2)^2=\gamma c^2\left(\kappa^2 -4\frac{\kappa_\gamma}{\kappa}+\gamma c^2\right)\\
  &=\gamma c^2\left(\kappa^2-4 + 4\left(1 -\frac{\kappa_\gamma}{\kappa}\right)+\gamma c^2\right)
\end{align*}
where we use $A=\beta^2\kappa$, $B=\alpha^2 \kappa$ and $A+B=\kappa$. However, since $\kappa >2$ and $\kappa>\kappa_\gamma$ for any $\gamma>0$, we have
\begin{equation*}
  \gamma c^2\left(\kappa^2-4 + 4\left(1 -\frac{\kappa_\gamma}{\kappa}\right)+\gamma c^2\right)>0
\end{equation*}
so that $N^-_n=0$.

Finally, Theorem~\ref{ThSpec1SW} tells us that $N^0_n=1$, while $n(\mathscr L)=1$.
Applying the counting argument, we conclude that $N_{C^+}=0$.
\QED

\subsection{Orbital stability}

 To discuss the orbital stability of solutions to \eqref{S1}-\eqref{W1}, it would be useful to write the system in a more convenient way by means of the change of variables   
  \begin{equation*}
    u_j=q_j+ip_j,\ \varphi_j= (-\Delta)^{1/2}\psi_j,\ p_j=\frac{\partial_t\psi_j}{c}.
  \end{equation*}
Hence, \eqref{S1}-\eqref{W1} reads as
\begin{equation}
  \label{SWexp}
  \begin{aligned}
    & \partial_t q_0=p_0-p_1+p_0\int_{\R^n}(-\Delta)^{-1/2}\sigma \varphi_0 \ud z& &\partial_t \varphi_0= c(-\Delta)^{1/2}\varpi_0\\
    & \partial_t p_0=-q_0+q_1-q_0\int_{\R^n}(-\Delta)^{-1/2}\sigma \varphi_0 \ud z& &\partial_t \varpi_0=-c(-\Delta)^{1/2}\varphi_0-c\sigma(|q_0|^2+|p_0|^2)\\
    & \partial_t q_1=p_1-p_0+p_1\int_{\R^n}(-\Delta)^{-1/2}\sigma \varphi_1 \ud z& &\partial_t \varphi_1= c(-\Delta)^{1/2}\varpi_1\\
    & \partial_t p_1=-q_1+q_0-q_1\int_{\R^n}(-\Delta)^{-1/2}\sigma \varphi_1 \ud z& &\partial_t \varpi_1=-c(-\Delta)^{1/2}\varphi_1-c\sigma(|q_1|^2+|p_1|^2)
  \end{aligned}
\end{equation}
and it can be written as 
\begin{equation}
  \label{SW}
  \partial_t X=\mathscr J\nabla \mathscr H_{SW}(X).
\end{equation}
with $X=(S,W)\in \mathbb R^4\times (L^2(\mathbb R^n))^4$ as in \eqref{defXSW}, $\mathscr J$ defined by \eqref{defJ} and
\begin{align}
  \label{defHSWnew}
  \mathscr H_{SW}(X)=&\ \frac{|q_0-q_1|^2+|p_0-p_1|^2}{2}+\frac{1}{4}\int_{\R^n}(|\varpi_0|^2+|\varpi_1|^2+|\varphi_0|^2+|\varphi_1|^2)\ud z\nonumber\\
  &+\frac12 \int_{\R^n}(-\Delta)^{1/2}\sigma(\varphi_0(|q_0|^2+|p_0|^2)+\varphi_1(|q_1|^2+|p_1|^2))\ud z
\end{align}
Next, we denote by  
$F(S)=\frac{|X|^2}{2}=\frac{q_0^2+q_1^2+p_0^2+p_1^2}{2}$ and introduce the functional
\begin{equation*}
  \mathscr E (X)=\mathscr H_{SW}(X)+\omega F(S)
\end{equation*}
which is thus conserved by the dynamical system \eqref{SW}.
Let $X_*=(S_*,W_*)$ one of the special solutions described in subsection \ref{sec:solscoupled}. In particular, $S_*=(Q_{*0}, 0, Q_{*1},0)$ with 
\begin{align}
  \label{SpecSolSW}
  &\begin{pmatrix}
    Q_{*0}\\ Q_{*1}
  \end{pmatrix}= \frac1{\sqrt2}\begin{pmatrix}
    1\\ \tau
  \end{pmatrix},\ \omega=\frac{\kappa}{2}+\tau-1\text{ for }\eqref{SpecSol}\\
  \label{SpecSolSWextra}
  &\begin{pmatrix}
    Q_{*0}\\ Q_{*1}
  \end{pmatrix}= \begin{pmatrix}
    \alpha\\ \beta
  \end{pmatrix},\ \omega=\kappa-1\text{ for }\eqref{SpecSolextra}
\end{align}
and $W_*=(\varphi_{0*},0,\varphi_{1*}, 0)$ where $\varphi_{*j}=-|Q_{*j}|^2(-\Delta)^{-1/2}\sigma$.

Adapting the argument used for the asymptotic model,  
we consider the level set
  \[\mathscr S=\{X=(S,W),\ F(S)=F(S_*)=1/2\}.
  \]
and its tangent set given by
  \[
  T\mathscr S=\{(S,W), \ \nabla F(S_*)\cdot S=0\}.\]
  Note that $(S,W)=(q_0,p_0,q_1,p_1,W)\in T\mathscr S$ if and only if $Q_{*0}q_0+ Q_{*1}q_1=0$.
The orbit associated to $X_*$ is given by
\[\mathscr O=\Big\{
  (S_\theta,-|Q_{*0}|^2(-\Delta)^{-1/2}\sigma,0,-|Q_{*1}|^2(-\Delta)^{-1/2}\sigma,0),\
 S_\theta= R(\theta)S_*,\  \theta\in \mathbb R\Big\}.\]
  and 
$(T\mathscr O)^\perp$ is made of $(S,W)$ with $S=(q_0,p_0,q_1,p_1)$ such that $Q_{*0}p_0+Q_{*1} p_1=0$.

\begin{rmk}
In contrast to 
the observation made for the asymptotic problem in Remark~\ref{phase_inv},
and to a common property of Hamiltonian systems,  here the phase invariance property holds in a restricted sense:
the energy $\mathscr H(X)$, and $\mathscr E(X)$ as well,  is left unchanged when changing $X=(S,W)$ into $(R(\theta)S,W)$, 
where the rotation $R(\theta)$ acts only on a part of the variables.
\end{rmk}

  The Euler-Lagrange relation for the coupled problem can be reformulated as 
  \begin{equation}\label{0gradESW}\nabla\mathscr E(X_*)=0\end{equation}
  and $\mathscr L$, defined in \eqref{defL} or \eqref{defL1}, corresponds to the Hessian of $\mathscr E$ evaluated at $X_*$ given by \eqref{SpecSolSW} or \eqref{SpecSolSWextra} respectively.
  We wish to establish a coercivity estimate, on a certain subspace,
  for the quadratic form $ X\mapsto D^2\mathscr{E}(X_*)(X,X)$. This is a crucial property for establishing the orbital stability. A straightforward computation gives, for any  $X\in T\mathscr S\cap (T\mathscr O)^\perp$,
  \begin{align*}
    D^2\mathscr{E}(X_*)(X,X)=&\,\frac{Q_{*1}}{Q_{*0}}q_0^2 - q_1q_0 +2 Q_{*0}q_0 \int_{\R^n}(-\Delta)^{-1/2}\sigma \varphi_0\ud z+\frac{Q_{*1}}{Q_{*0}}p_0^2 - p_1p_0\\
    &+\frac{Q_{*0}}{Q_{*1}}q_1^2 - q_0q_1 +2 Q_{*1}q_1 \int_{\R^n}(-\Delta)^{-1/2}\sigma \varphi_1\ud z+\frac{Q_{*0}}{Q_{*1}}p_1^2 - p_0p_1\\
    &+\frac12(\|\varphi_0\|_{L^2(\mathbb R^n)}^2+\|\varphi_1\|_{L^2(\mathbb R^n)}^2)+ \frac12(\|\varpi_0\|_{L^2(\mathbb R^n)}^2+\|\varpi_1\|_{L^2(\mathbb R^n)}^2)
  \end{align*}
   where we use the fact that $(1-|Q_{*0}|^2\kappa+\omega)Q_{*0}=Q_{*1}$ and $(1-|Q_{*1}|^2\kappa+\omega)Q_{*1}=Q_{*0}$. Now, since $X\in T\mathscr S\cap (T\mathscr O)^\perp$, we have $Q_{*0}q_0+Q_{*1}q_1=0=Q_{*0}p_0+Q_{*1}p_1$, so that 
   \begin{align*}
    D^2\mathscr{E}(X_*)(X,X)=&\,\left(\frac{1}{Q_{*0}Q_{*1}}\right)(q_0^2+q_1^2) +2 Q_{*0}q_0 \int_{\R^n}(-\Delta)^{-1/2}\sigma \varphi_0\ud z\\
    &+2 Q_{*1}q_1 \int_{\R^n}(-\Delta)^{-1/2}\sigma \varphi_1\ud z+\left(\frac{1}{Q_{*0}Q_{*1}}\right)(p_0^2+p_1^2)\\
    &+\frac12(\|\varphi_0\|_{L^2(\mathbb R^n)}^2+\|\varphi_1\|_{L^2(\mathbb R^n)}^2)+ \frac12(\|\varpi_0\|_{L^2(\mathbb R^n)}^2+\|\varpi_1\|_{L^2(\mathbb R^n)}^2)
  \end{align*}

  By virtue of the Cauchy-Schwarz inequality, 
  \begin{equation*}
    2|Q_{*j}||q_j|\left|\int_{\R^n}(-\Delta)^{-1/2}\sigma \varphi_j\ud z\right|\le 2|Q_{*j}||q_j| \sqrt\kappa\|\varphi_j\|_{L^2}\le \frac{\kappa}{\epsilon} |Q_{*j}|^2 q_j^2 +\epsilon \|\varphi_j\|_{L^2}^2
  \end{equation*}
  for any $\epsilon>0$. Therefore, 
  \begin{align}
    \label{coercHessiangenSW}
    D^2\mathscr{E}(X_*)(X,X)\ge &\,\left(\frac{1}{Q_{*0}Q_{*1}}-\frac{\kappa}{\epsilon}|Q_{0*}|^2\right)q_0^2+ \left(\frac{1}{Q_{*0}Q_{*1}}-\frac{\kappa}{\epsilon}|Q_{1*}|^2\right)q_1^2 \nonumber\\
    &+\left(\frac{1}{Q_{*0}Q_{*1}}\right)(p_0^2+p_1^2)+\left(\frac12-\epsilon\right)(\|\varphi_0\|_{L^2(\mathbb R^n)}^2+\|\varphi_1\|_{L^2(\mathbb R^n)}^2)\nonumber\\&+ \frac12(\|\varpi_0\|_{L^2(\mathbb R^n)}^2+\|\varpi_1\|_{L^2(\mathbb R^n)}^2).
  \end{align}
  Assume that $X_*$ is such that $(Q_{*0}, Q_{*1})$ is given by \eqref{SpecSolSW}. Then Proposition \ref{prop:spectralstabSW} implies that $X_*$ is spectrally stable only for $\tau=1$ and $0<\kappa<2$. In this case, \eqref{coercHessiangenSW} reads as
  \begin{align*}
    D^2\mathscr{E}(X_*)(X,X)\ge &\,\left(2-\frac{\kappa}{2\epsilon}\right)(q_0^2+ q_1^2)+2(p_0^2+p_1^2)+\left(\frac12-\epsilon\right)(\|\varphi_0\|_{L^2(\mathbb R^n)}^2+\|\varphi_1\|_{L^2(\mathbb R^n)}^2)\nonumber\\&+ \frac12(\|\varpi_0\|_{L^2(\mathbb R^n)}^2+\|\varpi_1\|_{L^2(\mathbb R^n)}^2)\ge C(\epsilon)\|X\|^2.
  \end{align*}
  with $C(\epsilon)=\min\big\{\left(2-\frac{\kappa}{2\epsilon}\right), \left(\frac12-\epsilon\right)\big\}$. Note that $C(\epsilon)>0$ provided $\epsilon$ is chosen such that $\frac{\kappa}{4}<\epsilon<\frac12$. This leads to the orbital stability of $X_*$ given by \eqref{SpecSolSW} when $\tau=1$ and $0<\kappa<2$.

  Note that if $\tau=1$ and $\kappa>2$ or $\tau=-1$, then the quadratic form $X\mapsto D^2\mathscr{E}(X_*)(X,X)$
   has no definite sign on $T\mathscr S\cap (T\mathscr O)^\perp$.

 Next, let $\kappa>2$ and let $X_*$ be such that $(Q_{*0}, Q_{*1})$ is given by \eqref{SpecSolSWextra}. Then Proposition~\ref{prop:spectralstabSWextra} implies that $X_*$ is always spectrally stable. In this case, \eqref{coercHessiangenSW} reads as
  \begin{align*}
    D^2\mathscr{E}(X_*)(X,X)\ge &\,\left(\kappa-\frac{\kappa}{\epsilon}\alpha^2\right)q_0^2+ \left(\kappa-\frac{\kappa}{\epsilon}\beta^2\right)q_1^2+\kappa(p_0^2+p_1^2)\nonumber\\&+\left(\frac12-\epsilon\right)(\|\varphi_0\|_{L^2(\mathbb R^n)}^2+\|\varphi_1\|_{L^2(\mathbb R^n)}^2)+ \frac12(\|\varpi_0\|_{L^2(\mathbb R^n)}^2+\|\varpi_1\|_{L^2(\mathbb R^n)}^2)
  \end{align*}
  where we use that $\alpha\beta=\frac1{\kappa}$. Next, we consider separately the cases $\tau=1$ and $\tau=-1$.

  If $\tau=1$, we write $q_1=-\frac{\alpha}{\beta}q_0$ so that 
  \begin{align*}
    D^2\mathscr{E}(X_*)(X,X)\ge &\,\kappa\left(\frac{1}{\beta^2}-\frac{2}{\epsilon}\alpha^2\right)q_0^2+\kappa(p_0^2+p_1^2)\nonumber\\&+\left(\frac12-\epsilon\right)(\|\varphi_0\|_{L^2(\mathbb R^n)}^2+\|\varphi_1\|_{L^2(\mathbb R^n)}^2)+ \frac12(\|\varpi_0\|_{L^2(\mathbb R^n)}^2+\|\varpi_1\|_{L^2(\mathbb R^n)}^2)\\
    =& \,\frac{\kappa}{\beta^2\epsilon}\left(\epsilon-\frac{2}{\kappa^2}\right)\left(1-\frac{\alpha^2}{\beta^2}\right)q_0^2+\frac{\kappa}{\beta^2\epsilon}\left(\epsilon-\frac{2}{\kappa^2}\right)q_1^2+\kappa(p_0^2+p_1^2)\nonumber\\&+\left(\frac12-\epsilon\right)(\|\varphi_0\|_{L^2(\mathbb R^n)}^2+\|\varphi_1\|_{L^2(\mathbb R^n)}^2)+ \frac12(\|\varpi_0\|_{L^2(\mathbb R^n)}^2+\|\varpi_1\|_{L^2(\mathbb R^n)}^2).
  \end{align*}
  By choosing $\frac{2}{\kappa^2}<\epsilon<\frac{1}{2}$ which is possible since $\kappa>2$ and since $\frac{\alpha^2}{\beta^2}<1$, we obtain $D^2\mathscr{E}(X_*)(X,X)\ge C(\epsilon)\|X\|^2$ with $C(\epsilon)>0$ and for any $X\in T\mathscr S\cap (T\mathscr O)^\perp$. 
  
  If $\tau=-1$, we write $q_0=-\frac{\beta}{\alpha}q_1$ so that 
  \begin{align*}
    D^2\mathscr{E}(X_*)(X,X)\ge & \,\frac{\kappa}{\alpha^2\epsilon}\left(\epsilon-\frac{2}{\kappa^2}\right)q_0^2+\frac{\kappa}{\alpha^2\epsilon}\left(\epsilon-\frac{2}{\kappa^2}\right)\left(1-\frac{\beta^2}{\alpha^2}\right)q_1^2+\kappa(p_0^2+p_1^2)\nonumber\\&+\left(\frac12-\epsilon\right)(\|\varphi_0\|_{L^2(\mathbb R^n)}^2+\|\varphi_1\|_{L^2(\mathbb R^n)}^2)+ \frac12(\|\varpi_0\|_{L^2(\mathbb R^n)}^2+\|\varpi_1\|_{L^2(\mathbb R^n)}^2).\\
  \end{align*}
  and we can conclude as above. 

\QED

\subsection{Instability}
\label{Sec:insta}

 In this section we study the nonlinear instability of the solution $X_*$ given by \eqref{SpecSolSW} whenever $\tau=-1$ or $\tau=1$ and $\kappa>2$, \emph{i.e.} whenever $X_*$ is spectrally unstable. To this goal, we use again the same change of variables as in the previous section so that \eqref{S1}-\eqref{W1} reads as \eqref{SWexp}.
Note that the reasoning of \cite{Maeda}, as described in Section~\ref{insta_As}, can be applied only in the case $\tau=1$ and $\kappa>2$,
the Morse index of $\mathscr L$ being larger than 2 when $\tau=-1$.
Hence, we are going to apply the general result described in \cite{ShSt}  to treat both cases at the same time.
The instability analysis is of different nature than in Section~\ref{insta_As}.
In Section~\ref{insta_As}, the method of \cite{Maeda} relies on the spectral property of the self-adjoint operator $\mathscr L$; it shows a linear growth of the perturbation by using the energy conservation, but it requires a strong assumption on the dimension of the eigenspace 
 of unstable directions. Here, the arguments of \cite{ShSt}, which has been extended to various 
 type of non linear Schr\"odinger equation in \cite{CCO, GeOh},
 uses the fact that $\mathbb L$ admits an eigenvalue with a positive real part.
 This property can be deduced from the counting argument.
 As pointed out in \cite{ShSt}, we show that 
 the perturbation, starting from the worst linearly unstable
  direction, exits a certain ball in a time which scales like the logarithm 
  of the inverse of the size of the initial perturbation.
 
We start by observing that the linearized operator $\mathbb L$ satisfies
\[
( \mathbb LX|X)
=-\ds\frac1{\sqrt2}\ds\int_{\mathbb R^n}(-\Delta)^{-1/2}\sigma(\varphi_0p_0+\tau \varphi_1p_1)\ud z
-c\sqrt 2\ds\int_{\mathbb R^n}\sigma(\varpi_0q_0+\tau \varpi_1q_1)\ud z.\]
The Cauchy-Schwarz inequality yields
\[|( \mathbb LX|X)|\leq 2(\sqrt{\kappa/2} +c\sqrt2\|\sigma\|_{L^2(\mathbb R^n)}) \|X\|^2.\]
As it will be detailed below, the operator $\lambda -\mathbb L$ is onto for sufficiently large (real part of) $\lambda$'s.
Accordingly, we can apply Lumer-Phillips' theorem \cite[Th.~12.22]{Ren} to the linearized equation 
\[\partial_ t X=\mathbb L X.\]
It can be formulated as the existence of the semi-group $t\mapsto e^{\mathbb Lt}$, which satisfies 
the continuity estimate: there exists $\Lambda>0$ 
such that  for any $t\geq 0$, $\|e^{\mathbb Lt}\|\leq e^{\Lambda t}$.
For further purposes, we denote 
\[K_0=\sup\big\{\|e^{\mathbb Lt}\|, 0\leq t\leq 1\big\}.\]

Then, we express the problem by considering the evolution of a perturbation of $X_*$ 
according to the dynamical system \eqref{SWexp}.  More precisely, we set 
\begin{align*}
 \begin{pmatrix}
  \tilde q_j\\ \tilde p_j
 \end{pmatrix} = R(\omega t)\begin{pmatrix}
  Q_{*j}+q_j\\ p_j
 \end{pmatrix},\ \tilde \varphi_j= \varphi_{*j}+\varphi_j,\ \tilde \varphi_j=\varpi_j.  
\end{align*}
From \eqref{SWexp}, we deduce that
the perturbation  $Y=(q_0,p_0,q_1,p_1,\varphi_0,\varpi_0,\varphi_1,\varpi_1)$ satisfies
\[
\partial_t Y=\mathbb LY+F(Y)\]
where the nonlinear remainder reads
\[
F(Y)=
\begin{pmatrix}
p_0\ds\int_{\mathbb R^n} \sigma(-\Delta)^{-1/2}\varphi_0\ud z
\\
-q_0\ds\int_{\mathbb R^n} \sigma(-\Delta)^{-1/2}\varphi_0\ud z
\\
p_1\ds\int_{\mathbb R^n} \sigma(-\Delta)^{-1/2}\varphi_1\ud z
\\
-q_1\ds\int_{\mathbb R^n} \sigma(-\Delta)^{-1/2}\phi_1\ud z
\\0
\\
- c\sigma(|p_0|^2+|q_0|^2)
\\
0
\\
-c\sigma(|p_1|^2+|q_1|^2)
\end{pmatrix}.\]
The orbital stability of the solution  \eqref{SpecSolSW} to \eqref{SWexp}
is rephrased in the orbital stability of $0$ for this problem. More precisely, we shall obtain the critical estimates by using 
 the integral formulation
\begin{equation}\label{Duhamel}
Y(t)=
e^{\mathbb Lt} Y_{\mathrm{init}} + \ds\int_0^t e^{\mathbb L(t-s)}F(Y(s))\ud s
\end{equation}
of the problem.
The application of the reasonings in \cite{ShSt} relies on the following estimate
\begin{lemma}\label{EstNL}
There exists $C_1>0$ such that for any $X$, we have 
$|F(X)|\leq C_1|X|^2$.
\end{lemma}
 
 \noindent
 {\bf Proof.}
 In order to estimate $F(Y)$, we make the following quantities appear
 \[
 (|p_j|^2+|q_j|^2)\left( \ds\int_{\mathbb R^n} (-\Delta)^{-1/2}\sigma\phi_j\ud z\right)^2
 \leq \kappa(|p_j|^2+|q_j|^2)\|\varphi_j\|^2_{L^2(\mathbb R^n)}\]
 and
 \[(|p_j|^2+|q_j|^2)^2
 \ds\int_{\mathbb R^n}|\sigma|^2\ud z =\|\sigma\|^2_{L^2(\mathbb R^n)}
 (|p_j|^2+|q_j|^2)^2.\]
 It leads to the asserted conclusion with $C_1=2(\sqrt\kappa+c\|\sigma\|_{L^2(\mathbb R^n)})$.

 \QED
  
Next, we need the following information on the spectrum of $\mathbb L$.  
\begin{proposition}\label{pr_Pruss}
  $\sigma(e^{\mathbb L})=e^{\sigma(\mathbb L)}$.
\end{proposition}

This statement strengthens the embedding 
 $\exp(\sigma(\mathbb L))\subset \sigma(e^{\mathbb L})$ which  always holds.
 It is not a prerequisite but it simplifies the argument, see   \cite{ShSt}.
 According to 
Gearhart-Greiner-Herbst-Pr\"uss Spectral Mapping Theorem, see e.g.
\cite[Prop.~1]{Pruss} (in fact, we use the criterion in the same form as  in \cite[Section~2]{Gesz}), the proof relies on a uniform estimate on the resolvent $(\lambda-\mathbb L)^{-1}$, as $\mathrm{Im}(\lambda)\to \pm \infty $ with $\mathrm{Re}(\lambda) \neq 0$ fixed, that we are going to establish.
We denote $$\mathbb H=\mathbb R^4\times(L^2(\mathbb R^n))^4$$
endowed with the norm
\[
\|X\|_{\mathbb H}=\sqrt{
|q_0|^2+|p_0|^2+|q_1|^2+|p_1|^2+\|\varphi_0\|_{L^2(\mathbb R^n)}^2
+\|\varpi_0\|^2_{L^2(\mathbb R^n)}
+\|\varphi_1\|_{L^2(\mathbb R^n)}^2+\|\varpi_1\|^2_{L^2(\mathbb R^n)}
}
.\]
Let $\lambda\in \mathbb C\setminus\{0\}$ and for a given  data $X'$, 
we consider the equation\[(\lambda-\mathbb L)X=X',\]
that is
\[\begin{array}{l}
\lambda q_0-\tau p_0+p_1=q'_0,
\\
\lambda p_0+\tau q_0-q_1+\ds\frac{1}{\sqrt2}\ds\int_{\mathbb R^n} (-\Delta)^{-1/2}\sigma \varphi_0\ud z=p'_0,
\\
\lambda q_1+p_0-\tau p_1=q'_1,
\\
\lambda p_1-q_0+\tau q_1+\ds\frac{\tau}{\sqrt2} \ds\int_{\mathbb R^n} (-\Delta)^{-1/2}\sigma \varphi_1\ud z=p'_1,
\\
\lambda \varphi_0-c(-\Delta)^{1/2}\varpi_0=  \varphi_0',
\\
\lambda \varpi_0+c(-\Delta)^{1/2}\varphi_0+c\sqrt 2\sigma q_0=  \varpi_0',
\\
\lambda \varphi_1-c(-\Delta)^{1/2}\varpi_1=  \varphi_1',
\\
\lambda \varpi_1+c(-\Delta)^{1/2}\varphi_1+\tau c\sqrt 2\sigma q_1=  \varpi_1'.
\end{array}
\]
Therefore, we get 
\[
\varpi_0=\ds\frac{(-\Delta)^{-1/2}}{c}(\lambda\varphi_0-\varphi_0'),\qquad
\varpi_1=\ds\frac{(-\Delta)^{-1/2}}{c}(\lambda\varphi_1-\varphi_1'),\]
which allows us to write
\[\begin{array}{l}
\Big(\ds\frac{\lambda^2}{c^2}-\Delta\Big)\varphi_0= \ds\frac{\lambda}{c^2} \varphi_0'
+(-\Delta)^{1/2}\varpi_0'-\sqrt 2(-\Delta)^{1/2}\sigma q_0,
\\[.3cm]
\Big(\ds\frac{\lambda^2}{c^2}-\Delta\Big)\varphi_1= \ds\frac{\lambda}{c^2} \varphi_1'
+(-\Delta)^{1/2}\varpi_1' - \tau \sqrt2(-\Delta)^{1/2}\sigma q_1.
\end{array}\]
We solve these equations by means of Fourier transform. Note that this makes the symbol $\big(\frac{\lambda^2}{c^2}+c^2\xi^2\big)$ appear. However, it does not vanish out of the axis $i\mathbb R$. 
Hence, we still can use the function
\[z\in \mathbb C\setminus i\mathbb R \longmapsto \kappa_z=\ds\int_{\mathbb R^n} \ds\frac{|\widehat \sigma(\xi)|^2}{z^2+|\xi|^2}\ds\frac{\ud \xi}{(2\pi)^n}\]
As consequence, we arrive at the reduced system:
\[\begin{array}{l}
\lambda q_0-\tau p_0+p_1=q'_0,
\\
\lambda p_0+\tau q_0-q_1-\kappa_{\lambda^2/c^2}q_0=S_0,
\\
\lambda q_1+p_0-\tau p_1=q'_1,
\\
\lambda p_1-q_0+\tau q_1 - \kappa_{\lambda^2/c^2}q_1=S_1,
\end{array}\]
where we have set 
\begin{equation}\label{datarhs}
\begin{array}{l}
S_0=p'_0-\ds\frac{1}{\sqrt 2}\ds\int_{\mathbb R^n} \ds\frac{\widehat  \sigma(\xi) \widehat {\varpi_0'}(\xi)}{\lambda^2/c^2+|\xi|^2}
\ds\frac{\ud \xi}{(2\pi)^n}
-\ds\frac{\lambda}{\sqrt 2c^2}\ds\int_{\mathbb R^n} \ds\frac{\widehat  \sigma(\xi) \widehat {\varphi_0'}(\xi)}{(\lambda^2/c^2+|\xi|^2)|\xi|}
\ds\frac{\ud \xi}{(2\pi)^n},
\\[.3cm]
S_1=p'_1-\ds\frac{\tau}{\sqrt 2}\ds\int_{\mathbb R^n} \ds\frac{\widehat  \sigma(\xi) \widehat {\varpi_1'}(\xi)}{\lambda^2/c^2+|\xi|^2}
\ds\frac{\ud \xi}{(2\pi)^n}
-\tau \ds\frac{\lambda}{\sqrt 2c^2}\ds\int_{\mathbb R^n} \ds\frac{\widehat  \sigma(\xi) \widehat {\varphi_1'}(\xi)}{(\lambda^2/c^2+|\xi|^2)|\xi|}
\ds\frac{\ud \xi}{(2\pi)^n}.
\end{array}\end{equation}
Since 
$$\lambda(q_0+\tau q_1)=q'_0+\tau q'_1,$$
we obtain
$$\begin{array}{l}
\lambda p_0 +2\tau  q_0-\kappa_{\lambda^2/c^2}q_0=S_0+\tau\ds\frac{(q'_0+\tau q'_1)}{\lambda},
\\[.3cm]
\lambda p_1 -2  q_0
+\tau \kappa_{\lambda^2/c^2}q_0=S_1-
\tau(\tau-\kappa_{\lambda^2/c^2})\ds\frac{(q'_0+\tau q'_1)}{\lambda}
.
\end{array}
$$
It eventually yields
\begin{equation}\label{resolvent}
(\lambda^2 + 4- 2\tau \kappa_{\lambda^2/c^2})q_0
=\lambda q'_0
- (S_1-\tau S_0)+(2-\tau\kappa_{\lambda^2/c^2})\frac{(q_0'+\tau q_1')}{\lambda}
\end{equation}
which already explains (when setting $X'=0$) the relation \eqref{nllam} for studying the eigenvalues of $\mathbb L$.
Next, we are going to use the following elementary claim.

\begin{lemma}
\label{cp}
Let $\lambda=a+ib\in \mathbb C$, with $a$ and $ b$ reals.  If $|b|\geq \sqrt 3|a|$, then, for any $\epsilon\geq 0$ we have
$$
\Big|\ds\frac{1}{\lambda^2+\epsilon}\Big|\leq \ds\frac {\sqrt2}{|\lambda|^2},\qquad
\Big|\ds\frac{\lambda}{\lambda^2+\epsilon}\Big|\leq \ds\frac {\sqrt2}{|\lambda|}.$$
\end{lemma}

\noindent
{\bf Proof.} 
We write $\lambda=re^{i\theta}$ with $r=\sqrt{a^2+b^2}$, so that 
$$
\Big|\ds\frac{1}{\lambda^2+\epsilon}\Big|=\Big|\ds\frac{1}{e^{i\theta} r^2+e^{-i\theta}\epsilon}\Big|,\qquad 
\Big|\ds\frac{\lambda}{\lambda^2+\epsilon}\Big|=\Big|\ds\frac{1}{e^{i\theta} r+e^{-i\theta}\epsilon/r}\Big|
.$$
Now, we 
re-organize \[\begin{array}{l}
|e^{i\theta} r^2+e^{-i\theta}\epsilon|^2=r^4+\epsilon^2+2r^2\epsilon\cos(2\theta)
=(r^2- \epsilon)^2+ 4r^2\epsilon \cos^2(\theta)
\\
\hspace*{4cm}\geq \ds\frac{(r^2-\epsilon)^2}{2}+\ds\frac{r^4+\epsilon^2}{2}\geq \ds\frac{r^4}{2}
,
\\
|e^{i\theta} r+e^{-i\theta}\epsilon/r|^2 = r^2+\ds\frac{\epsilon^2}{r^2}+2\epsilon\cos(2\theta)
=\Big(r-\ds\frac\epsilon r\Big)^2+4\epsilon\cos^2(\theta)
\\
\hspace*{4cm}\geq 
\ds\frac{(r-\epsilon/ r)^2}{2}+\ds\frac{r^2+\epsilon^2/r^2}{2}
\geq \ds\frac{r^2}{2}.
\end{array}\]
where we have used that, by assumptions on $a,b$, $|\cos(\theta)|=\frac{|a|}{\sqrt{a^2+b^2}}\leq \frac12$. It thus implies
\[\Big|\ds\frac{1}{\lambda^2+\epsilon}\Big|\leq  \ds\frac{\sqrt 2}{r^2}, \qquad  
\Big|\ds\frac{\lambda}{\lambda^2+\epsilon}\Big|\leq \ds\frac{\sqrt 2}{r}.\]
\QED

\noindent
This allows us to estimate the resolvent $(\lambda-\mathbb L)^{-1}$. 

\begin{lemma}
Let $\lambda=a+ib\in \mathbb C$, with $a\neq 0$, $|b|\geq \sqrt 3|a|$.
Then, there exists a constant $C_a>0$ such that  the quantities $S_0, S_1$ in \eqref{datarhs} satisfy
\[
|S_j|\leq C_a \|X'\|_{\mathbb H}.
\]
\end{lemma}

\noindent
{\bf Proof.} 
The only difficulty is to estimate the integrals involving $\sigma$. Owing to Lemma~\ref{cp}
and the Cauchy-Schwarz inequality, we obtain 
\[\begin{array}{lll}
\Big|
\ds\frac{\lambda}{c^2}\ds\int_{\mathbb R^n} \ds\frac{\widehat  \sigma(\xi) \widehat {\varphi_0'}(\xi)}{(\lambda^2/c^2+|\xi|^2)|\xi|}\ds\frac{\ud \xi}{(2\pi)^n} \Big|
&\leq& \ds\frac{\sqrt2}{|\lambda|}\ds\int_{\mathbb R^n} \ds\frac{|\widehat  \sigma(\xi)|\ | \widehat {\varphi_0'}(\xi)|}{|\xi|}\ds\frac{\ud \xi}{(2\pi)^n} 
\\
&\leq&
\ds\frac{\sqrt2}{|\lambda|}\left(\ds\int_{\mathbb R^n} \ds\frac{|\widehat  \sigma(\xi)|^2}{|\xi|^2}\ds\frac{\ud \xi}{(2\pi)^n} 
\right)^{1/2}
\left(\ds\int_{\mathbb R^n}|\widehat  \varphi_0'(\xi)|^2\ds\frac{\ud \xi}{(2\pi)^n} 
\right)^{1/2}
\\[.3cm]
&\leq&
\ds\frac{\sqrt{2\kappa} }{|\lambda|}\|\varphi'_0\|_{L^2(\mathbb R^n)}
\leq 
\ds\frac{\sqrt\kappa }{\sqrt 2|a|}\|\varphi'_0\|_{L^2(\mathbb R^n)}.
\end{array}\]
Similarly, we get 
\[\begin{array}{lll}
\Big|
\ds\frac{1}{c^2}\ds\int_{\mathbb R^n} \ds\frac{\widehat  \sigma(\xi) \widehat {\varpi_0'}(\xi)}{\lambda^2/c^2+|\xi|^2}\ds\frac{\ud \xi}{(2\pi)^n} \Big|
&\leq& \ds\frac{\sqrt 2}{|\lambda|^2}\ds\int_{\mathbb R^n}
|\widehat  \sigma(\xi)|\ | \widehat {\varpi_0'}(\xi)|
\ds\frac{\ud \xi}{(2\pi)^n} 
\\[.3cm]
&\leq&
\ds\frac{\sqrt 2}{|\lambda|^2}\|\sigma\|_{L^2(\mathbb R^n)}
\|\varpi'_0\|_{L^2(\mathbb R^n)}
\leq 
\ds\frac{\sqrt 2\|\sigma\|_{L^2(\mathbb R^n)}  }{4|a|^2}\|\varpi'_0\|_{L^2(\mathbb R^n)}.
\end{array}\]
\QED

\noindent
By direct inspection, Lemma~\ref{cp} also yields the following estimate.

\begin{lemma}\label{klc}
Let $\lambda=a+ib\in \mathbb C$, with $a\neq 0$ and $|b|\geq \sqrt 3|a|$.
Then, we have 
$$|
\kappa_{\lambda^2/c^2}|
\leq c^2\frac{\|\sigma\|_{L^2(\mathbb R^n)}^2}{2\sqrt{2}|a|^2}.$$
\end{lemma}

Let $\lambda=a+ib\in \mathbb C$. By virtue of Lemma~\ref{klc},
when $b$ is large enough, 
$\lambda^2 + 4-2\kappa_{\lambda^2/c^2}$ does not vanish.
We can therefore obtain $q_0$ from the data $X'$ with \eqref{resolvent}.
Moreover, as $b\to \infty$, with $a\neq 0$ fixed, 
$\frac{1}{\lambda^2 + 4-2\kappa_{\lambda^2/c^2}}$,
 and $\frac{\lambda}{\lambda^2 + 4-2\kappa_{\lambda^2/c^2}}$ both tend to 0.
 We conclude that we can find some $r>0$ and $M>0$ 
 (depending on $a$, $c$, $\sigma$) such that for any $b\in \mathbb R$, $|b|\geq r$, 
we have $\|X\|_{\mathbb H}=\|(a+ib-\mathbb L)^{-1}X'\|_{\mathbb H}
\leq M\|X'\|_{\mathbb H}$.
This justifies Proposition~\ref{pr_Pruss}.
\QED

 In case of spectral instability, $\mathbb L$ admits eigenvalues with positive real value.
 There is only a finite number of such eigenvalues (as indicated by the counting argument).
 Since, $\exp(\sigma(\mathbb L))\subset \sigma(e^{\mathbb L})$ we thus already know that the spectral radius of 
 $e^{\mathbb L}$ is larger than 1. In fact, we can use the identity in Proposition~\ref{pr_Pruss}.
  Let us denote 
  \[
 \lambda_*=a_*+ib_*\in \sigma(\mathbb L),\quad
 a_*=\sup\big\{\mathrm{Re}(\lambda),\ \lambda\in \sigma(\mathbb L)\big\}>0.\]
  Of course, for any $t\geq 0$, we have 
 $|e^{\lambda_* t}|=e^{a_*t}$ and the spectral radius of $e^{\mathbb L}$ is $e^{a_*}>1$, see \cite{GeOh} for more details.
We are going to use the following claim.

\begin{lemma}\cite[Lemma~2 \& Lemma~3]{ShSt}
There exists a constant $K_1$, such that for any $t\geq 0$, there holds
$\|e^{t\mathbb L}\|_{\mathscr L(\mathbb H)}\leq K_1e^{3a_*t/2}$.\label{LShSt}
\end{lemma}

Let us define $\epsilon>0$ such that 
\[\ds\frac{4 K_1(1+C_1)^2}{a_*}\epsilon<1
\]
with $C_1$ and $K_1$ defined in Lemma \ref{EstNL} and \ref{LShSt} respectively.
Then, pick an arbitrary $0<\delta <\epsilon$ and set
$$T_\epsilon=\ds\frac{1}{a_*}\ln\Big(\ds\frac\epsilon \delta\Big)
$$
Let $Y_*$ be a normalized eigenvector of $\mathbb L$ associated to $\lambda_*$:
\[\mathbb LY_*=\lambda_*Y_*,\qquad \|Y_*\|_{\mathbb H}=1.\]
It will serve to define the initial perturbation that leads to instability:
we start from the perturbation
\[
Y\big|_{t=0}=\delta Y_*,\]
which has thus an arbitrarily small norm.
As a matter of fact, \eqref{Duhamel} becomes
\[
Y(t)
=\delta e^{\lambda_*t} Y_{*} + \ds\int_0^t e^{\mathbb L(t-s)}F(Y(s))\ud s.\]
We 
are going to
contradict the orbital stability by showing that  $Y(T_\epsilon)$
is at a distance larger than  
$\kappa \epsilon$, for a certain constant $\kappa>0$, to the orbit $\mathscr O$.
Let $$\tilde T_\epsilon
=\ds\sup\big\{t\in [0,T_\epsilon],\ \|Y(s)\|\leq (1+C_1)\delta e^{a_*s} \textrm{ for } 0\leq s\leq t\big\}
\in (0,T_\epsilon].$$ 
The Duhamel formula \eqref{Duhamel} yields
\[
\|Y(t)\|\leq \delta e^{a_*t}+\ds\int_0^t  K_1 e^{3a_*(t-s)/2}C_1\|Y(s)\|^2 \ud s
\]
by using Lemma~\ref{EstNL} and~\ref{LShSt}.
Therefore 
\[
\begin{array}{lll}
\|Y(t)\|&\leq& \delta e^{a_*t}+K_1C_1 (1+C_1)^2 \delta^2 \ds\int_0^t e^{3a_*(t-s)/2} e^{2a_*s}\ud s
\\[.3cm]
&\leq&\delta e^{a_*t}+K_1C_1 (1+C_1)^2 \delta^2 \ds\frac{2e^{2a_*t}}{a_*}
\\[.3cm]
&\leq& \delta e^{a_*t}\left(1+\ds\frac{2K_1C_1(1+C_1)^2}{a_*}\delta e^{a_*T_\epsilon}\right)
\leq  \delta e^{a_*t}\left(1+C_1\ds\frac{2K_1(1+C_1)^2}{a_*}\epsilon\right)
\end{array}\]
holds for any $t\in [0,\tilde T_\epsilon]\subset [0, T_\epsilon]$.
Hence, $\epsilon$ is chosen small enough so that 
this implies 
\[
\|Y(t)\|
< \left(1+\ds\frac{C_1}{2}\right)\delta e^{a_*t},\]
which would contradict the definition of $\tilde T_\epsilon$ if $\tilde T_\epsilon< T_\epsilon$.
We deduce that 
\[\|Y(t)\|\leq (1+C_1)\delta e^{a_*t}\leq (1+C_1)\epsilon\]
holds for any $t\in [0,T_\epsilon]$.
Owing to this estimate, we go back to the Duhamel formula 
and we obtain, for $0\leq t\leq T_\epsilon$,
\begin{equation}\label{est_Duh_key}\begin{array}{lll}
\|Y(t)-\delta e^{\lambda_* t}Y_*\|
&\leq& 
\ds\int_0^t |e^{\mathbb L(t-s)}F(Y(s))|\ud s
\leq 
\ds\int_0^t e^{3a_*(t-s)/2} K_1C_1 \|Y(s)\|^2\ud s
\\
&\leq&
K_1C_1 (1+C_1)^2\delta^2 \ds\int_0^t e^{3a_*(t-s)/2} e^{2a_*s}\ud s
\leq\ds\frac{2K_1C_1 (1+C_1)^2}{a_*}\delta^2 e^{2a_*t}
\\
&\leq& \ds\frac{2K_1C_1 (1+C_1)^2}{a_*}\delta^2 e^{2a_*T_\epsilon}
=
 \ds\frac{2K_1C_1 (1+C_1)^2}{a_*}\epsilon^2
.
\end{array}
\end{equation}
We distinguish the components of the solution 
$X_*=(S_*,W_*)$, $Y(t)=(\tilde S(t),\tilde W(t))$ and
$X(t)=(S(t),W(t))=(R(\omega t)(S_*+\tilde S(t)),W_*+\tilde W(t))
$. We wish to evaluate 
\[\begin{array}{lll}
\Xi_\epsilon&=&\ds\inf_\theta\|X(T_\epsilon)-(R(\theta)S_*,W_*)\|
\\&=&
\ds\inf_\theta\|(R(\omega T_\epsilon)(S_*+\tilde S(T_\epsilon)),W_*+\tilde W(T_\epsilon))-(R(\theta)S_*,W_*)\|
\\
&=&
\ds\inf_\theta\|(S_*+\tilde S(T_\epsilon),W_*+\tilde W(T_\epsilon))-(R(-\omega T_\epsilon)R(\theta)S_*,W_*)\|
\\&=&
\ds\inf_{\theta'}\| Y(T_\epsilon)+X_*-(R(\theta')S_*,W_*)\|
.\end{array}\]
Let $\theta_\epsilon$ denote the phase which reaches this infimum:
$$\Xi_\epsilon=\| Y(T_\epsilon)+X_*-(R(\theta_\epsilon)S_*,W_*)\|.$$
We observe that
\[
\Xi_\epsilon\leq \ds\inf_{\theta'}\big(\| Y(T_\epsilon)\|+\|X_*-(R(\theta')S_*,W_*)\|\big)
\leq \| Y(T_\epsilon)\|\leq (1+C_1)\epsilon.\]
Next, we have
\[\|X_*-(R(\theta_\epsilon)S_*,W_*)\|\leq 
\Xi_\epsilon+\| Y(T_\epsilon)\|\leq 2(1+C_1)\epsilon,\]
which implies that $\lim_{\epsilon\to 0}\theta_\epsilon=0$.
Hence, a basic Taylor expansion tells us that
\[
X_*-(R(\theta_\epsilon)S_*,W_*)=(-\theta_\epsilon \mathscr J_S S_*,0)+\epsilon r_\epsilon,\qquad 
\lim_{\epsilon\to 0}\|r_\epsilon\|=0.\]
Now, we are going to use the following splitting of the initial perturbation
$$Y_*=\big(Y_*|(\mathscr J_S S_*,0) \big)\ds\frac{(\mathscr J_S S_*,0)}{\|(\mathscr J_S S_*,0)\|^2}+Y_*^\perp,\qquad 
\big(Y_*^\perp|(\mathscr J_S S_*,0)\big)=0.$$
The Cauchy-Schwarz inequality yields 
\[\begin{array}{lll}
\Xi_\epsilon \|Y_*^\perp\| &\geq& \big|\big(Y(T_\epsilon)+X_*-(R(\theta_\epsilon)S_*,W_*)| Y_*^\perp\big)\big|
\\
&\geq& \big|
\delta e^{\lambda_*T_\epsilon}\big(Y_*
| Y_*^\perp\big)+\big(Y(T_\epsilon)-\delta e^{\lambda_*T_\epsilon}Y_*
| Y_*^\perp\big)
-\theta_\epsilon \underbrace{\big((\mathscr J_S S_*,0)|Y_*^\perp\big)}_{=0}+\epsilon \big(r_\epsilon|Y_*^\perp\big)\big|
\end{array}\]
Possibly at the price of choosing a smaller $\epsilon$, coming back to 
\eqref{est_Duh_key}, we can make both quantities
\[
\big|\big(Y(T_\epsilon)-\delta e^{\lambda_*T_\epsilon}Y_*
| Y_*^\perp\big)\big|\leq 
\big\|\big(Y(T_\epsilon)-\delta e^{\lambda_*T_\epsilon}Y_*\big\|
\| Y_*^\perp\|
\textrm{ and } 
\epsilon \big|\big(r_\epsilon|Y_*^\perp\big)\big|
\leq \epsilon \|r_\epsilon\| \|Y_*^\perp\|\]
smaller than 
$\frac\epsilon 4\|Y_*^\perp\|^2.$
It follows that 
\[\begin{array}{lll}
\Xi_\epsilon \|Y_*^\perp\|
& \geq &  \big|
\delta e^{\lambda_*T_\epsilon}\big(Y_*
| Y_*^\perp\big)\big|
-\big|\big(Y(T_\epsilon)-\delta e^{\lambda_*T_\epsilon}Y_*\big)
| Y_*^\perp\big)\big|
-\epsilon \big|\big(r_\epsilon|Y_*^\perp\big)\big|
\\[.3cm]
& \geq &  
\delta e^{a_*T_\epsilon}\|Y_*^\perp\|^2 - 
\ds\frac {\epsilon}2 \|Y_*^\perp\|^2=
\ds\frac\epsilon2\|Y_*^\perp\|^2.
\end{array}\]
%
%
%
This estimate is meaningful provided $Y^\perp_*\neq 0$.
This is indeed the case because $\mathscr J_S S_*=\frac1{\sqrt 2}(0, -1,0,-\tau)$ and we can check that
$(\mathscr J_S S_*,0)$ lies in $\mathrm{Ker}(\mathbb L)$
while $Y_*\in \mathrm{Ker}(\mathbb L-\lambda_*)$,  with $\lambda_*\neq 0$. 
\QED

\appendix
\section{Proof of $L^2$ and energy conservation properties}
\label{App_Cons}

The three models 
can be cast under the general form 
\[
i\ds\frac{\ud}{\ud t}\begin{pmatrix}
u_0\\
u_1
\end{pmatrix}=
\begin{pmatrix}
A_{0} & -1 \\
-1 & A_1
\end{pmatrix}
\begin{pmatrix}
u_0\\
u_1
\end{pmatrix}\]
where $A_0=A_1=1$ for 
\eqref{0cou}, $A_0=1-\kappa |u_0|^2 $,  $A_1=1-\kappa |u_1|^2$ for \eqref{Hartree},
and 
$A_0=1+\int \sigma \psi_0\ud z $,  $A_1=1+\int \sigma \psi_1\ud z $
for  \eqref{S1}-\eqref{W1}. In any case, $A_0$ and $A_1$ are real.
Therefore, we obtain 
\[\begin{array}{lll}
\ds\frac{\ud}{\ud t}(|u_0|^2 + |u_1|^2)
&=&\ds\frac {\overline {u_0}}i (A_0 u_0-u_1) -\ds\frac {u_0}i (A_0 \overline {u_0}-\overline {u_1}) 
+\ds\frac {\overline {u_1} }i (A_1 u_1-u_0) -\ds\frac {u_1}i (A_1 \overline {u_1}-\overline {u_0})
\\[.4cm]
&=&
\ds\frac{A_0}i(\overline {u_0}u_0-u_0\overline {u_0})
+\ds\frac{A_1}i(\overline {u_1}u_1-u_1\overline {u_1})
+\ds\frac1i(
-\overline {u_0}u_1
+\overline {u_1}u_0
-u_0\overline {u_1}+u_1\overline {u_0}
)\\[.4cm]
&=&0,
\end{array}\]
which proves the conservation of $|u_0|^2 + |u_1|^2$.

Moreover, we have
\[\begin{array}{lll}
\ds\frac12\ds\frac{\ud}{\ud t}|u_0 - u_1|^2
& = &-\ds\frac12 \ds\frac{\ud}{\ud t}(u_0\overline {u_1}+ u_1\overline {u_0})
\\[.4cm]
& =&
-\ds\frac12\left(\ds\frac {\overline {u_1}}i (A_0 u_0-u_1) -\ds\frac {u_1}i (A_0 \overline {u_0}-\overline {u_1}) 
+\ds\frac {\overline {u_0} }i (A_1 u_1-u_0) -\ds\frac {u_0}i (A_1 \overline {u_1}-\overline {u_0})\right)
\\[.4cm]
& =&
-\ds\frac1{2i}\left(A_0(\overline {u_1}u_0-u_1\overline {u_0})
+A_1(\overline {u_0}u_1-u_0\overline {u_1})\right)
=-(A_0-A_1)\mathrm{Im}(u_0\overline{u_1}).
\end{array}\]
For \eqref{Hartree}, this 
combines to
\[\begin{array}{lll}
\ds\frac\kappa4 \ds\frac{\ud}{\ud t}(|u_0|^4+|u_1|^4)
& =&
\ds\frac{\kappa|u_0|^2}{2}\left(
\ds\frac {\overline {u_0}}i (A_0 u_0-u_1) -\ds\frac {u_0}i (A_0 \overline {u_0}-\overline {u_1}) 
\right)
\\[.4cm]
&&+\ds\frac{\kappa|u_1|^2}{2}\left(
\ds\frac {\overline {u_1} }i (A_1 u_1-u_0) -\ds\frac {u_1}i (A_1 \overline {u_1}-\overline {u_0})
\right)
\\[.4cm]
&=&
-\ds\frac{\kappa|u_0|^2}{2i}(\overline {u_0}u_1-u_0\overline {u_1})
-
\ds\frac{\kappa|u_1|^2}{2i}(\overline {u_1}u_0-u_1\overline {u_0})
\\[.4cm]
&=&
\kappa(|u_0|^2-|u_1|^2) \mathrm{Im}(u_0\overline{u_1})
=-(A_0-A_1)\mathrm{Im}(u_0\overline{u_1}),
\end{array}\]
so that  \eqref{EnerHa} holds.
For \eqref{S1}-\eqref{W1}, we also compute
the energy of the vibrational field
\[\begin{array}{l}
\ds\frac12 \ds\frac{\ud}{\ud t}
\ds\int_{\mathbb R^n}
\Big(
\ds\frac{1}{c^2}
(|\partial_t \psi_0|^2+|\partial_t \psi_1|^2) +|\nabla\psi_0|^2+|\nabla\psi_1|^2\Big)\ud z
\\[.4cm]
=
\ds\int_{\mathbb R^n}\left\{
\Big(
\ds\frac{1}{c^2}
\partial_{tt}^2 \psi_0-\Delta \psi_0\Big)\partial_t\psi_0
+\Big(
\ds\frac{1}{c^2}\partial_{tt}^2 \psi_1-\Delta \psi_1\Big)\partial_t\psi_1\right\}\ud z
\\[.4cm]
=
-\ds\int_{\mathbb R^n}\sigma
\big(
|u_0|^2\partial_t\psi_0
+|u_1|^2\partial_t\psi_1\big)\ud z.
\end{array}\]
Finally, we compute the evolution of the  interaction energy
\[
\begin{array}{lll}
\ds\frac{\ud}{\ud t}
\ds\int_{\mathbb R^n}\sigma (\psi_0|u_0|^2+\psi_1|u_1|^2)
\ud z
&=&
\ds\int_{\mathbb R^n}\sigma
\big(
|u_0|^2\partial_t\psi_0
+|u_1|^2\partial_t\psi_1\big)\ud z
\\[.4cm]
&&+
\ds\int_{\mathbb R^n}\sigma
\psi_0\left(\ds\frac {\overline {u_0}}i (A_0 u_0-u_1) -\ds\frac {u_0}i (A_0 \overline {u_0}-\overline {u_1})\right)
\ud z
\\[.4cm]
&&+
\ds\int_{\mathbb R^n}\sigma
\psi_1\left(\ds\frac {\overline {u_1}}i (A_1 u_1-u_0) -\ds\frac {u_1}i (A_1 \overline {u_1}-\overline {u_0})\right)
\ud z\\[.4cm]
&=&
\ds\int_{\mathbb R^n}\sigma
\big(
|u_0|^2\partial_t\psi_0
+|u_1|^2\partial_t\psi_1\big)\ud z
\\[.4cm]
&&-
\ds\frac1i\ds\int_{\mathbb R^n}\sigma\Big(
\psi_0(\overline {u_0}u_1-\overline {u_1}u_0)
+
\psi_1(\overline {u_1}u_0-\overline {u_0}u_1)
\Big)
\ud z
\\[.4cm]
&=&
\ds\int_{\mathbb R^n}\sigma
\big(
|u_0|^2\partial_t\psi_0
+|u_1|^2\partial_t\psi_1\big)\ud z
+
2(A_0-A_1)\mathrm{Im}(u_0\overline{u_1})
.
\end{array}\]
Gathering these identities, we arrive at \eqref{EnerSW}.

\section{Proof of Lemma~\ref{Plem}.}
\label{App:Plem}

We extend $\Sigma$ by 0 on $(-\infty,0)$ and we assume that $\Sigma$ is supported in $[-R+\sqrt \mu,R+\sqrt \mu]$, for some $0<R<\infty$. 
Extending the discussion to a function with fast decay at infinity follows from a standard density argument.
We start by defining the principal value
\[
\mathrm{P.V.}\ds\int_{0}^{+\infty} \ds\frac{\Sigma(r)}{(r-\sqrt \mu)(r+\sqrt \mu)}\ud r
=\lim_{\epsilon\to 0}\ds\int_{-R+\sqrt \mu}^{R+\sqrt \mu}
 \mathbf 1_{|r-\sqrt \mu|\geq \epsilon}\ds\frac{\Sigma(r)}{(r-\sqrt \mu)(r+\sqrt \mu)}\ud r.
\]
We decompose
\[
\ds\frac{1}{(r-\sqrt \mu)(r+\sqrt \mu)}=
 \ds\frac{1}{2\sqrt\mu(r-\sqrt \mu)}- \ds\frac{1}{2\sqrt \mu(r+\sqrt \mu)}.\]
There is no difficulty in handling the last  term by means of the Lebesgue's dominated convergence theorem 
(bearing in mind that $\Sigma$
is supported on $[0,\infty)$) and we obtain
\[\begin{array}{l}
\ds\lim_{\epsilon\to 0}\ds\int_{-R+\sqrt \mu}^{R+\sqrt \mu} \mathbf 1_{|r-\sqrt \mu|\geq \epsilon} 
 \ds\frac{\Sigma(r)}{2\sqrt\mu(r+\sqrt \mu)}\ud r
=
\ds\int_{-R+\sqrt \mu}^{R+\sqrt \mu}
\ds\frac{\Sigma(r)}{2\sqrt \mu(r+\sqrt \mu)}\ud r
.
\end{array}\]
Next, we make use of parity so that for any $\epsilon>0$
\[
\ds\int_{-R+\sqrt \mu}^{R+\sqrt \mu}
\mathbf 1_{|r-\sqrt \mu|\geq \epsilon}\ds\frac{\ud r}{r-\sqrt \mu}
=0.\]
Hence, we rewrite
\[\begin{array}{lll}
\ds\ds\lim_{\epsilon\to 0}\ds\int_{-R+\sqrt \mu}^{R+\sqrt \mu} \mathbf 1_{|r-\sqrt \mu|\geq \epsilon}
\ds\frac{\Sigma(r)}{r-\sqrt \mu}\ud r&=&\ds\lim_{\epsilon\to 0}\ds\int_{-R+\sqrt \mu}^{R+\sqrt \mu} \mathbf 1_{|r-\sqrt \mu|\geq \epsilon}
\ds\frac{\Sigma(r)-\Sigma (\sqrt \mu)}{r-\sqrt \mu}\ud r
\\[.4cm]
&=&\ds\int_{-R+\sqrt \mu}^{R+\sqrt \mu} 
\ds\frac{\Sigma(r)-\Sigma (\sqrt \mu)}{r-\sqrt \mu}\ud r
\end{array}\]
which is well defined since the integrand is bounded by $\|\Sigma'\|_{L^\infty}$ and the integral is over a bounded domain.
We conclude that 
\[\begin{array}{l}\mathrm{P.V.}\ds\int_{0}^{+\infty} \ds\frac{\Sigma(r)}{(r-\sqrt \mu)(r+\sqrt \mu)}\ud r
\\[.4cm]
=
-\ds\int_{-R+\sqrt \mu}^{R+\sqrt \mu}
\ds\frac{\Sigma(r)}{2\mu(r+\sqrt \mu)}\ud r
+
\ds\int_{-R+\sqrt \mu}^{R+\sqrt \mu} 
\ds\frac{\Sigma(r)-\Sigma (\sqrt \mu)}{2\sqrt\mu(r-\sqrt \mu)}\ud r.
\end{array}\]
We split $P(-\mu,B)$ into its real and imaginary parts; it yields
\[ P(-\mu,B)=\ds\int_{0}^\infty \ds\frac{ (r^2-\mu)\Sigma(r)}{B^2+ (r^2-\mu)^2}\ud r
-i B\ds\int_{0}^\infty \ds\frac{\Sigma(r)}{B^2+ (r^2-\mu)^2}  \ud r
.
\]
With $B>0$, and the change of variable $r'=\frac{r-\sqrt \mu}{B}$, the imaginary parts recasts as
\[
\ds\int_{0}^\infty \ds\frac{\Sigma(r)}{1+ ((r-\sqrt\mu)/B)^2(r+\sqrt\mu)^2}  \ds\frac{\ud r}{B}
=
\ds\int_{-\sqrt\mu/B}^\infty \ds\frac{\Sigma(Br'+\sqrt\mu)}{1+ r'^2(Br'+2\sqrt\mu)^2}  \ud r'.\]
A direct application of the Lebesgue's dominated convegence theorem shows that it tends to 
\[
\Sigma(\sqrt\mu)\ds\int_{-\infty}^\infty \ds\frac{1}{1+ 4\mu r'^2}  \ud r'
=
\ds\frac{\pi\Sigma(\sqrt\mu)}{2\sqrt \mu}
\]
as $B\to 0^+$.
Taking the limit $B\to 0^-$, changes the sign of this expression.

We proceed in two steps to handle the real part.
Let $\epsilon>0$ and compute
\[
\begin{array}{l}
\ds\int_{|r-\sqrt \mu|\leq \epsilon} \ds\frac{ (r^2-\mu)\Sigma(r)}{B^2+ (r^2-\mu)^2}\ud r
-\ds\int_{|r-\sqrt \mu|\leq \epsilon} \ds\frac{ 2\sqrt\mu (r-\sqrt\mu)\Sigma(r)}{B^2+ 4\mu(r-\sqrt\mu)^2}\ud r
\\[.4cm]
\qquad
=
\ds\int_{|r-\sqrt \mu|\leq \epsilon}\Sigma(r)(r-\sqrt\mu)
\\
\qquad\qquad\times \ds\frac{
B^2(r+\sqrt\mu -2\sqrt \mu)+
4\mu(r-\sqrt\mu)^2(r+\sqrt\mu)-2\sqrt\mu(r-\sqrt\mu)^2(r+\sqrt\mu)^2 }
{(B^2+ (r^2-\mu)^2)(B^2+ 4\mu(r-\sqrt\mu)^2)}\ud r
\\[.4cm]
\qquad
=
\ds\int_{|r-\sqrt \mu|\leq \epsilon}\Sigma(r)(r-\sqrt\mu)^2 
\ds\frac{B^2
+(r-\sqrt\mu)(r+\sqrt\mu)(4\mu-2\sqrt\mu(r+\sqrt\mu))}{(B^2+ (r^2-\mu)^2)(B^2+ 4\mu(r-\sqrt\mu)^2)}\ud r
\\[.4cm]
\qquad
=
\ds\int_{|r-\sqrt \mu|\leq \epsilon}\Sigma(r)(r-\sqrt\mu)^2 \ds\frac{
B^2
-2\sqrt\mu(r-\sqrt\mu)^2(r+\sqrt\mu)
}{(B^2+ (r^2-\mu)^2)(B^2+ 4\mu(r-\sqrt\mu)^2)}\ud r.
\end{array}\]
This difference is dominated by
\[\begin{array}{l}
\ds\int_{|r-\sqrt \mu|\leq \epsilon}\Sigma(r) 
\left(
\ds\frac{B^2(r-\sqrt\mu)^2 }{B^24\mu(r-\sqrt\mu)^2}
+\ds\frac{2\sqrt\mu(r-\sqrt\mu)^4(r+\sqrt\mu)}{
4\mu(r-\sqrt\mu)^4(r+\sqrt\mu)^2}
\right)
\ud r
\\[.4cm]
\qquad
=
\ds\int_{|r-\sqrt \mu|\leq \epsilon}\Sigma(r)
\left(\ds\frac{1 }{4\mu}+\ds\frac{1}{
2\sqrt\mu(r+\sqrt\mu)}
\right)
\ud r
.\end{array}\]
Pick $\delta>0$. Since $r\mapsto \Sigma(r)\big(\frac{1 }{4\mu}+\frac{1}{
2\sqrt\mu(r+\sqrt\mu)}
\big)$ is integrable over $(0,\infty)$, 
this quantity can be made $\leq \delta$ by choosing $\epsilon$ small enough.
Possibly at the price of reducing $\epsilon$, we also suppose that
\[
\left|
\ds\int_{0}^{+\infty} \ds\frac{\Sigma(r)}{2\sqrt\mu(r+\sqrt \mu)}\ud r
-\ds\int_{|r-\sqrt \mu|\geq \epsilon} \ds\frac{\Sigma(r)}{2\sqrt\mu(r+\sqrt\mu)}\ud r
\right|\leq \delta
\]
holds.
Having disposed of this preliminary,  we write
\[\begin{array}{lll}\ds\int \ds\frac{ (r^2-\mu)\Sigma(r)}{B^2+ (r^2-\mu)^2}\ud r
&=&\ds\int_{|r-\sqrt \mu|\leq \epsilon} \ds\frac{ (r^2-\mu)\Sigma(r)}{B^2+ (r^2-\mu)^2}\ud r
+
\ds\int_{|r-\sqrt \mu|\geq \epsilon} \ds\frac{ (r^2-\mu)\Sigma(r)}{B^2+ (r^2-\mu)^2}\ud r
\\[.4cm]
&=&\left(\ds\int_{|r-\sqrt \mu|\leq \epsilon} \ds\frac{ (r^2-\mu)\Sigma(r)}{B^2+ (r^2-\mu)^2}\ud r
-\ds\int_{|r-\sqrt \mu|\leq \epsilon} \ds\frac{ 2\sqrt\mu(r-\sqrt\mu)\Sigma(r)}{B^2+ 4\mu(r-\sqrt\mu)^2}\ud r\right)
\\[.4cm]
&&+
\ds\int_{|r-\sqrt \mu|\leq \epsilon} \ds\frac{ 2\sqrt\mu(r-\sqrt\mu)\Sigma(r)}{B^2+ 4\mu(r-\sqrt\mu)^2}\ud r
+
\ds\int_{|r-\sqrt \mu|\geq \epsilon} \ds\frac{ (r^2-\mu)\Sigma(r)}{B^2+ (r^2-\mu)^2}\ud r
\end{array}\]
where the first term can be made $\leq \delta$, uniformly with respect to $B$, \emph{i.e.}
\[
\ds\sup_{B\in \mathbb R}\left|
\ds\int_{|r-\sqrt \mu|\leq \epsilon} \ds\frac{ (r^2-\mu)\Sigma(r)}{B^2+ (r^2-\mu)^2}\ud r
-\ds\int_{|r-\sqrt \mu|\leq \epsilon} \ds\frac{ 2\sqrt\mu(r-\sqrt\mu)\Sigma(r)}{B^2+ 4\mu(r-\sqrt\mu)^2}\ud r
\right|\leq \delta.\]
The limit of the last integral is identified by applying Lebesgue's dominated convergence theorem
\[\begin{array}{l}
\ds\lim_{B\to 0}
\ds\int_{|r-\sqrt \mu|\geq \epsilon} \ds\frac{ (r^2-\mu)\Sigma(r)}{B^2+ (r^2-\mu)^2}\ud r
=
\ds\int_{|r-\sqrt \mu|\geq \epsilon} \ds\frac{ \Sigma(r)}{ (r^2-\mu)}\ud r
\\[.4cm]
\qquad
=
\ds\int_{|r-\sqrt \mu|\geq \epsilon} \ds\frac{ \Sigma(r)}{ 2\sqrt\mu(r-\sqrt\mu)}\ud r-
\ds\int_{|r-\sqrt \mu|\geq \epsilon} \ds\frac{ \Sigma(r)}{ 2\sqrt\mu(r+\sqrt\mu)}\ud r
\\[.4cm]
\qquad=
\ds\int_{-R+\sqrt \mu}^{R+\sqrt\mu} \mathbf 1_{|r-\sqrt \mu|\geq \epsilon}\ds\frac{ \Sigma(r)-\Sigma(\sqrt\mu)}{ 2\sqrt\mu(r-\sqrt\mu)}\ud r-
\ds\int_{|r-\sqrt \mu|\geq \epsilon} \ds\frac{ \Sigma(r)}{ 2\sqrt\mu(r+\sqrt\mu)}\ud r
.
\end{array}\]
Finally, the second term can be recast as
\begin{equation*}
  \ds\int_{-R+\sqrt \mu}^{R+\sqrt\mu} \mathbf 1_{|r-\sqrt \mu|\leq \epsilon} \ds\frac{ 2\sqrt\mu(r-\sqrt\mu)(\Sigma(r)-\Sigma(\sqrt\mu))}{B^2+ 4\mu(r-\sqrt\mu)^2}\ud r
\end{equation*}
so that, in the limit as $B$ goes to $0$, we get
\begin{equation*}
  \ds\int_{-R+\sqrt \mu}^{R+\sqrt\mu} \mathbf 1_{|r-\sqrt \mu|\leq \epsilon} \ds\frac{ \Sigma(r)-\Sigma(\sqrt\mu)}{2\sqrt\mu(r-\sqrt\mu)}\ud r.
\end{equation*}

Therefore, we get
\[\begin{array}{l}
\ds\lim_{B\to 0}\left(\ds\int_{|r-\sqrt \mu|\leq \epsilon} \ds\frac{ 2\sqrt\mu(r-\sqrt\mu)\Sigma(r)}{B^2+ 4\mu(r-\sqrt\mu)^2}\ud r
+
\ds\int_{|r-\sqrt \mu|\geq \epsilon} \ds\frac{ (r^2-\mu)\Sigma(r)}{B^2+ (r^2-\mu)^2}\ud r\right)
\\[.4cm]
\qquad\qquad=
\ds\int_{-R+\sqrt\mu}^{R+\sqrt \mu} \ds\frac{ \Sigma(r)-\Sigma(\sqrt\mu)}{ 2\sqrt\mu(r-\sqrt\mu)}\ud r-
\ds\int_{|r-\sqrt \mu|\geq \epsilon} \ds\frac{ \Sigma(r)}{ 2\sqrt\mu(r+\sqrt\mu)}\ud r,
\end{array}\]
which is close, up to $\delta$, to $\mathrm{P.V.}\int_0^\infty \frac{\Sigma(r)}{(r-\sqrt\mu)(r+\sqrt\mu)}\ud r$. As a consequence, we conclude that, for any $\delta>0$, we can exhibit $B(\delta)>0$ small enough so that
\[
\left|
\ds\int \ds\frac{ (r^2-\mu)\Sigma(r)}{B^2+ (r^2-\mu)^2}\ud r-
\mathrm{P.V.}\int_0^\infty \frac{\Sigma(r)}{(r-\sqrt\mu)(r+\sqrt\mu)}\ud r\right|\leq2\delta
\]
holds for any $0<|B|<B(\delta)$.
\QED

%
%

\bibliography{SchroPW2}
\bibliographystyle{plain}

\end{document}